\documentclass[a4paper, 10pt]{article}
\usepackage[T2A]{fontenc}
\usepackage[english]{babel}

\usepackage{amsmath, amsthm}
\usepackage{amssymb}
\usepackage{amsfonts}
\usepackage{srcltx, dsfont, color, graphicx}

\theoremstyle{plain}
\newtheorem{theorem}{Theorem}[section]
\newtheorem{lemma}{Lemma}[section]

\numberwithin{equation}{section}

\def\tht{\theta}
\def\Om{\Omega}
\def\om{\omega}
\def\e{\varepsilon}
\def\g{\gamma}
\def\G{\Gamma}
\def\l{\lambda}
\def\p{\partial}
\def\D{\Delta}
\def\a{\alpha}
\def\b{\beta}
\def\Ups{\Upsilon}

\def\d{\delta}
\def\L{\Lambda}
\def\z{\zeta}

\def\vp{\varphi}

\def\vr{\varrho}
\def\vt{\vartheta}
\def\vk{\varkappa}

\def\Ho{\mathring{W}}

\def\hf{\mathfrak{h}}

\def\la{\langle}
\def\ra{\rangle}

\def\iu{\mathrm{i}}

\def\rA{\mathrm{A}}

\def\cL{\mathcal{L}}

\def\sE{\mathsf{E}}

\def\Ho{\mathring{W}}

\DeclareMathOperator{\RE}{Re}

\DeclareMathOperator{\mes}{mes}

\DeclareMathOperator{\dist}{dist}
\DeclareMathOperator{\supp}{supp}

\DeclareMathOperator{\dvr}{div}

\allowdisplaybreaks

\begin{document}

\title{Operator estimates for non-periodically perforated domains with Dirichlet and nonlinear Robin conditions: strange term}

\author{D.I. Borisov$^{1,2,3}$
}

\date{\empty}

\maketitle

{\small
    \begin{quote}
1) Institute of Mathematics, Ufa Federal Research Center, Russian Academy of Sciences,  Chernyshevsky str. 112, Ufa, Russia, 450008
\\
2) Bashkir State  University, Zaki Validi str. 32, Ufa, Russia, 450076
\\
3) University of Hradec Kr\'alov\'e
62, Rokitansk\'eho, Hradec Kr\'alov\'e 50003, Czech Republic
\\
Emails: borisovdi@yandex.ru
\end{quote}

{\small
\begin{quote}
\noindent{\bf Abstract.}  We consider a boundary value problem for a general second order linear equation in a domain with a fine perforation. The latter is made by small cavities; both the shapes of the cavities and their distribution are arbitrary. The boundaries of the cavities are subject either to a Dirichlet or a nonlinear Robin condition. On the perforation, certain rather weak conditions are imposed to ensure that under the homogenization we obtain a similar problem in a non-perforated domain with an additional potential in the equation usually called a strange term. Our main results state  the convergence of the solution of the perturbed problem to that of the homogenized one in $W_2^1$- and $L_2$-norms uniformly in $L_2$-norm of the right hand side in the equation. The estimates for the convergence rates are established and their order sharpness is discussed.

\medskip


\medskip

 	
\end{quote}
}

\section{Introduction}

Nowadays a new direction in the homogenization theory devoted to so-called operator estimates is quite intensively developed. In contrast to classical results in the homogenization theory on strong and weak convergence of the solutions, here the studies are aimed on proving the norm resolvent convergence and obtaining estimates for the convergence rates; the latter are often called operator estimates.
Recently, such results were obtained in few papers for problems in domains with fine perforation distributed along entire domain. Problems in such perforated domains are classical in the homogenization theory, see, for instance,  \cite{MK}, \cite{MK1}, \cite{OIS}, \cite{Sha}, \cite{ZKO}, and the references therein, and there are many results describing the convergence in a strong or weak sense in $L_2$ and $W_2^1$ for fixed right hand sides in the equations and boundary conditions. In   \cite{ChDR}, \cite{Khra}, \cite{Khra2}, \cite{Past}, \cite{Sus}, \cite{Zhi}  the classical results were improved and operator estimates were established for several cases of periodic and almost periodic perforation in arbitrary domains. The case of the Neumann condition was addressed  in \cite{Sus}, \cite{Past}, \cite{Zhi} and the sizes  of the cavities were of the same order as the distances between them and the perforation was purely periodic. In \cite{Khra2}, on the boundaries of the cavities  the Dirichlet condition was imposed and the   sizes of these cavities were assumed  to satisfy certain relation with respect to the size of the periodicity cell. All cavities were of the same shapes up to an arbitrary rotation and its location in the periodicity cell was also quite arbitrary. In \cite{Khra},  \cite{ChDR} the perforation was pure periodic and it was made by small balls with the Dirichlet or Neumann \cite{ChDR} or Robin \cite{Khra}  condition on the  boundaries.  The main results of the cited papers were the formulation of the homogenized problems  and various operator estimates; their order sharpness was not established.

A non-periodic perforation was studied in \cite{Post}. Here the domain was a manifold with a perforation made by arbitrary cavities with the Dirichlet or Neumann condition and the operator was the Laplacian. The main results were again homogenized problems and operator estimates and they were established under the validity of certain local upper bounds for $L_2$-norm in terms of $W_2^1$-norms. And these bounds were the main tools in proving the convergence and operator estimates. Then several cases of possible homogenized problems were addressed and as examples, it was shown that the developed scheme worked for perforation by small balls.

We also mention several recent papers on operator estimates for domain perforated along a given manifold  \cite{MSB21}, \cite{AA22}, \cite{PRSE16}. The perforation was  non-periodic and formed by arbitrary cavities and distribution. The homogenized problems were classified and a series of operator estimates was  established. In some cases these estimates turned out to be order sharp.

In this paper we study a  boundary value problem for a  linear second order elliptic equation in a perforated domain. The differential expression is general, involves complex-valued varying coefficients
and is not formally symmetric. The perforation is arbitrary and non-periodic and is assumed to satisfy natural geometric conditions.
On the boundaries of the cavities we impose the Dirichlet or a nonlinear Robin condition; both types of conditions can be simultaneously present on different cavities. Then we impose additional rather weak conditions on the perforation to describe the case, when the homogenization produces a so-called strange term, namely, when in the homogenized equation an additional potential appears. Our main results states the convergence of the perturbed solution to the homogenized one in $L_2$- and $W_2^1$-norms uniformly in $L_2$-norm of the right hand side in the equation. The estimates for the convergence rates are also proved and some terms in these estimates are shown to be order sharp. An important feature of our results is that our assumptions are rather weak and do not apriori require any local estimates like  in \cite{Post}. Instead of this  we prove that similar estimates are guaranteed by our assumptions. One more advantage of our study is that we can deal with a nonlinear Robin condition.

In conclusion we mention that a similar problem was studied in a very recent paper \cite{BK} but in the situation, when the solution to the perturbed problem vanishes as the perforation becomes finer. Operator estimates in such case were obtained and the convergence rates were shown to be order sharp.

\section{Problem and main results}\label{sec2}

\subsection{Formulation of problem}

Let $x=(x_1,\ldots,x_n)$ be Cartesian coordinates in $\mathds{R}^n$ and $\Om$ be  an arbitrary  domain in $\mathds{R}^n$; if its boundary is non-empty, we suppose that its smoothness is $C^2$. The domain $\Om$ can be  both  bounded or unbounded.  In this domain, we choose a family of points $M_k^\e$, $k\in\mathds{M}^\e$, where $\e$ is a small positive parameter and $\mathds{M}^\e$ is some at most countable set of indices. We also choose  a family of bounded non-empty domains $\om_{k,\e}\subset \mathds{R}^d$, $k\in\mathds{M}^\e$, with $C^2$-boundaries. Then we define
\begin{equation*}
\om_k^\e:=\big\{x:\, (x - M_k^\e)\e^{-1}\eta^{-1}(\e)\in \om_{k,\e}\big\}, \quad k\in\mathds{M}^\e,\qquad \tht^\e:=\bigcup\limits_{k\in\mathds{M}^\e} \om_k^\e,
\end{equation*}
where $\eta=\eta(\e)$ is some function obeying $0<\eta(\e)\leqslant 1$.
We shall formulate rigorously the assumptions on the cavities $\om_{k,\e}$ later, now we just say that they are assumed to be approximately of the same size (but not the shapes!) and there is a minimal distance between the points $M_k^\e$, which ensures that the domains $\om_k^\e$ are mutually disjoint.

By means of the domains $\om_k^\e$ we introduce a perforation of the domain $\Om$ as   $\Om^\e:=\Om\setminus\tht^\e$. In the perforated domain $\Om^\e$ we consider a  boundary value problem for an elliptic equation with the coefficients $A_{ij}=A_{ij}(x)$, $A_j=A_j(x)$, $A_0=A_0(x)$ defined in the non-perforated domain $\Om$, which are supposed to satisfy the conditions
\begin{gather}\label{2.5}
 A_{ij}\in W_\infty^1(\Om), \qquad A_j,\, A_0\in L_\infty(\Om),
\\
 A_{ij}=A_{ji}, \qquad
\sum\limits_{i,j=1}^{n} A_{ij}(x)\xi_i\overline{\xi_j}\geqslant c_0\sum\limits_{j=1}^{n} |\xi_j|^2,\qquad x\in\Om,\quad \xi_i\in\mathds{C},\label{2.5a}
\end{gather}
where $c_0>0$ is some fixed constant independent of $\xi$ and $x$. The functions $A_{ij}$ are real-valued, while the functions $A_j$ and $A_0$ are complex-valued.

The boundaries of the cavities $\om_k^\e$ are subject to either the Dirichlet condition or a nonlinear Robin condition. In order to introduce them, we first partition arbitrarily the set $\tht^\e$:
\begin{equation*}
\tht_D^\e:=\bigcup\limits_{k\in\mathds{M}_D^\e} \om_k^\e,\qquad \tht_R^\e:=\bigcup\limits_{k\in\mathds{M}_R^\e} \om_k^\e, \qquad \mathds{M}_D^\e\cup \mathds{M}_R^\e=\mathds{M}^\e, \qquad \mathds{M}_D^\e\cap \mathds{M}_R^\e=\emptyset.
\end{equation*}
For $x\in\p\tht_R^\e$ and $u\in\mathds{C}$ by  $a^\e=a^\e(x,u)$  we denote
a measurable complex-valued function, which will serve as a nonlinear term in the Robin condition; the main assumptions about this function will be formulated later.

The main object of our study is the following boundary value problem:
\begin{equation}\label{2.7}
(\cL-\l)u=f\quad\text{in}\quad \Om^\e,\qquad u=0\quad\text{on} \quad \p\Om\cup\p\tht_D^\e,\qquad  \frac{\p u}{\p\boldsymbol{\nu}} + a^\e(x,u)=0\quad\text{on}\quad \p\tht_R^\e.
\end{equation}
Here  $\cL$ and  $\frac{\p\ }{\p\boldsymbol{\nu}}$
are  a differential expression and a conormal derivative:
\begin{gather}
\nonumber
\cL:=-\dvr \rA \nabla  + \sum\limits_{j=1}^{n} A_j\frac{\p\ }{\p x_j}  + A_0,\qquad \frac{\p\ }{\p\boldsymbol{\nu}} =\nu\cdot\rA\nabla,
\\
\rA(x):=\begin{pmatrix}
A_{11}(x) & \ldots & A_{1n}(x)
\\
\vdots && \vdots
\\
A_{n1}(x) & \ldots & A_{nn}(x)
\end{pmatrix},
\nonumber
\end{gather}
$f\in L_2(\Om^\e)$ is an arbitrary function, $\l\in\mathds{C}$ is a fixed constant, $\nu$ is the unit normal to $\p\tht^\e_R$ directed inside $\tht_R^\e$.

Our main aim  is analyze the behavior of a generalized solution to problem (\ref{2.7})  as $\e\to+0$. Namely, we address two questions:  how does a homogenized problem (\ref{2.7}) read and whether the operator estimates can be established and if so, what are the corresponding convergence rates? It is very well known that the homogenized problem depends very much on the distribution and shapes of the cavities as well as on their sizes and the distances between them. In this paper consider the case, when a so-called strange term appears and we ensure such situation by a few assumptions on the cavities and   the nonlinearity $a^\e$ in the Robin condition. All of them will be formulated later, now we just say that the sizes of cavities, controlled by the function $\eta(\e)$, depends on the small parameter $\e$, governing the distances between the holes, as follows:
\begin{equation}\label{2.14}
\e^{-2}\eta^{n-2}(\e)\vk^{-1}(\e) \to \g,\qquad \e\to+0,
\end{equation}
where  $\g$ is a non-negative constant and
\begin{equation*}
\vk(\e):=|\ln\eta(\e)|+1\quad\text{as}\quad n=2,\qquad
\vk(\e):=1\quad\text{as}\quad n\geqslant 3.\nonumber
\end{equation*}
Convergence (\ref{2.14}) is the key point guaranteeing the appearance of the strange term in our model. Namely, we show that under our assumptions the homogenized problem reads as
\begin{equation}\label{2.11}
(\cL+\g\Ups\b-\l)u_0=f\quad\text{in}\quad \Om,\qquad u_0=0\quad\text{on} \quad \p\Om,
\end{equation}
where $\b\in L_\infty(\Om)$ is some function determined by the shapes and distribution of the cavities $\om_k^\e$ and $\Ups$ is just a fixed function:
\begin{equation*}
\Ups(x):=\sqrt{\det\rA(x)} \mes_{n-1}\p B_1(0), \qquad x\in\Om;
\end{equation*}
hereinafter by   $B_r(M)$ we denote a ball in $\mathds{R}^n$ of a radius $r$ centered at a point $M$, while $\mes_{n-1}$ denotes the $(n-1)$-dimensional measure on surfaces.  We observe that in view of ellipticity condition (\ref{2.5a})  the matrix $\rA(x)$ is symmetric, positive and bounded uniformly in $x$ and this is why the function $\Ups$ is well-defined. The assumed smoothness of the functions $A_{ij}$ implies that $\Ups\in W_\infty^1(\Om)$.

\subsection{Main assumptions}

In this subsection we formulate our main assumptions. We begin with a geometric assumption on the cavities $\om_k^\e$. In the vicinity of the boundaries $\p\om_{k,\e}$ we define a local variable  $\tau$ being the distance measured along the normal vector $\nu$ to $\p\om_{k,\e}$.

\begin{enumerate}
\def\theenumi{{A\arabic{enumi}}}
\item\label{A1} The points $M_k^\e$ and the domains $\om_{k,\e}$ obey the conditions
\begin{equation}
B_{R_1}(y_{k,\e})\subseteq \om_{k,\e}\subset  B_{R_2}(0), \quad
B_{\e R_3}(M_k^\e)\cap B_{\e R_3}(M_j^\e)=\emptyset, \quad \dist(M_k^\e,\p\Om)\geqslant R_3\e,
\quad k\ne j,\quad k,j\in\mathds{M}^\e,
\label{2.2}
\end{equation}
where $y_{k,\e}$ are some points, and $R_1<R_2<R_3$ are some fixed constants independent of $\e$, $\eta$, $k$ and $j$. The sets $B_{R_2}(0)\setminus \om_{k,\e}$ are connected. For each $k\in\mathds{M}_R^\e$ there exist local variables $s$ on $\p\om_{k,\e}$ such that the variables $(\tau,s)$ are  well-defined at least on $\{x\in\mathds{R}^n:\, \dist(x,\p\om_{k,\e})\leqslant \tau_0\}\subseteq B_{R_2}(0)$, where $\tau_0$ is a fixed constant independent of $k\in\mathds{M}^\e$ and $\e$. The Jacobians corresponding to passing from variables $x$ to $(\tau,s)$ are separated from zero and bounded from above uniformly in $\e$, $k\in\mathds{M}_R^\e$ and $x$ as $\dist(x,\p\om_{k,\e})\leqslant \tau_0$. The derivatives of $x$ with respect to $(\tau,s)$ and of $(\tau,s)$ with respect to $x$ up to the second order
are bounded uniformly in $\e$, $k\in\mathds{M}_R^\e$ and $x$ as $\dist(x,\p\om_{k,\e})\leqslant \tau_0$.
\end{enumerate}

\noindent The first relation in (\ref{2.2}) means that all domains $\om_{k,\e}$ are approximately of the same sizes: we can inscribe a fixed ball of the radius $R_1$ inside each domain, which in its turn is contained in a fixed ball $B_{R_2}(0)$. The second condition in (\ref{2.2}) guarantees that each two neighbouring cavities do not intersect and there is a minimal distance $2R_3$ between each two neighbouring points $M_k^\e$, while the third condition says that the cavities are not too close to the boundary of $\Om$, see Figure~\ref{fig1}. The connectedness of the domains $B_{R_2}(0)\setminus \om_{k,\e}$ is also a natural condition meaning that the perforation produces no new isolated connected components in the domain $\Om$. The rest of Assumption~\ref{A1} postulates a regularity of the boundaries $\p\om_{k,\e}$ uniformly in $k$ and $\e$.

Our second assumptions concern the function $a^\e$. We first suppose that
\begin{equation}\label{2.6}
\begin{aligned}
&
\RE (a^\e(x,u_1)-a^\e(x,u_2))\overline{(u_1-u_2)}
\geqslant - \mu_0(\e)
|u_1-u_2|^2,\qquad
\\
&
|a^\e(x,u_1)-a^\e(x,u_2)|\leqslant c_1 |u_1-u_2|,\qquad a^\e(x,0)=0,
\end{aligned}
\end{equation}
where $c_1$ is some constant independent of $\e$, $x\in\p\tht_R^\e$ and $u_1, u_2\in\mathds{C}$ and $\mu_0(\e)$ is some nonnegative function such that
\begin{gather}\label{2.10}
\e\eta(\e)\vk(\e)\mu_0(\e)\to+0,\quad \e\to+0.
\end{gather}
These conditions describe the class of admissible nonlinearities and in view of the technique we use, they guarantee the unique solvability of problem (\ref{2.7}).

\begin{figure}
\begin{center}
\includegraphics[scale=0.4]{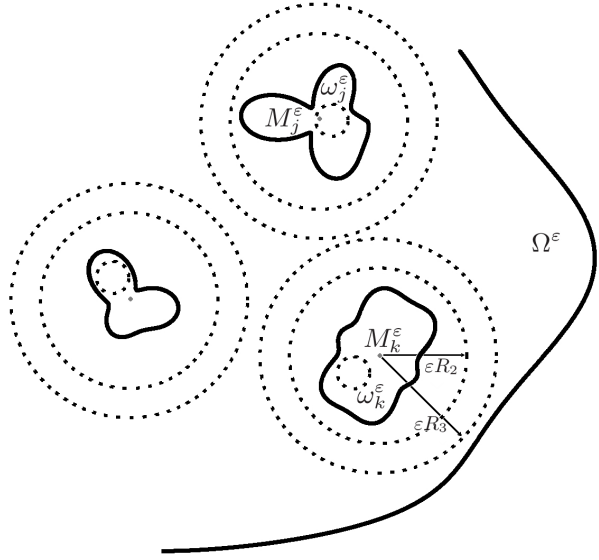}

\caption{{\small Domains $\om_k^\e$ (indicated by solid black lines) and points $M_k^\e$ (indicated by gray color). Dotted lines show (rescaled) balls $B_{\e R_1}(M_{k,\e}+\e y_{k,\e})$ and $B_{\e R_2}(M_k^\e)$ from the first condition in (\ref{2.2})
and balls $B_{\e R_3}(M_k^\e)$ from the second condition.}}\label{fig1}
\end{center}
\end{figure}

A more important assumption for $a^\e$ is as follows; it is needed to ensure that the homogenized problem is indeed (\ref{2.11}).

\begin{enumerate}
\def\theenumi{{A\arabic{enumi}}}\setcounter{enumi}{1}

\item\label{A6} The set $\mathds{M}_R^\e$ is partitioned into two disjoint subsets $\mathds{M}_{R,1}^\e$ and $\mathds{M}_{R,2}^\e$ obeying  the following   conditions:
\begin{align}
&\RE a^\e(x,u)\overline{u}\geqslant \mu_1(\e) |u|^2, \quad x\in\p\om_k^\e,\quad u\in\mathds{C},\quad k\in\mathds{M}_{R,1}^\e,
 \label{2.25}
 \\
 &\label{2.22}
a^\e(x,u)=\e^{-1}\eta^{-1}(\e) b_k\big((x-M_k^\e)\e^{-1}\eta^{-1},\e\big) u +\tilde{a}_k^\e(x,u),  \quad x\in\p\om_k^\e,\quad u\in\mathds{C},\quad k\in\mathds{M}_{R,2}^\e,
\end{align}
where  $\mu_1=\mu_1(\e)$ is a fixed function independent of $k\in\mathds{M}_{R,1}^\e$, the functions
 $b_k(\xi,\e)$, $\xi\in\p\om_k^\e$,
and $\tilde{a}_k^\e(x,u)$, $(x,u)\in\p\om_k^\e\times \mathds{C}$, are  complex-valued, the functions $b_k(\xi,\e)$ belong to $C^1(\p\om_{k,\e})$, while $\tilde{a}_k^\e(x,u)$ are measurable in $x$ and $u$ for each $\e$ and
\begin{align}\label{2.24}
&\e\eta(\e)\vk^{-1}(\e)\mu_1(\e)\to+\infty,\qquad \e\to+0,
\\
&
\begin{aligned}
&\|b_k(\,\cdot\,,\e)\|_{C^1(\p\om_{k,\e})}\leqslant c_3,
\qquad \RE b_k(\xi,\e)\geqslant c_2,\quad \xi\in\p\om_{k,\e},
\\
& |\tilde{a}_k^\e(x,u)|\leqslant \mu_2(\e)|u|,\quad x\in\p\om_k^\e, \quad u\in\mathds{C},
\qquad \e\eta(\e)\vk(\e)\mu_2(\e)\to+0,\quad \e\to+0,
\end{aligned}
\label{2.26}
\end{align}
$c_2$, $c_3$ are some positive constants independent of $\e$, $\xi$ and $k$ and $\mu_2$ is some function independent of $k$.
\end{enumerate}

This condition says that we deal with two main types of the nonlinear Robin condition. The first  is imposed  for $k\in\mathds{M}_{R,1}^\e$ and here the nonlinear term $a^\e$ is sign definite and large in the sense of inequality (\ref{2.25}) and convergence (\ref{2.24}). These conditions ensure that the corresponding cavities, for $k\in\mathds{M}_{R,1}^\e$, behave similar to ones with the Dirichlet condition for $k\in\mathds{M}_D^\e$. Namely, the traces of the function $u_\e$ on $\p\om_k^\e$ for $k\in\mathds{M}_{R,1}^\e$ tend to zero as $\e\to+0$.

The second type of the Robin condition is imposed for $k\in\mathds{M}_{R,2}^\e$ and here the function $a^\e$ is linear in the leading term as it is described by (\ref{2.22}), (\ref{2.26}). The coefficient  $\e^{-1}\eta^{-1}$ at the functions $b_k$ indicates the minimal growth of the linear term in the Robin condition; a faster growth is also allowed   since the functions $b_k$ can additionally depend on $\e$.

All cavities with both the Dirichlet and the Robin conditions contribute to the function $\b$ in the strange term in (\ref{2.11}). The contribution of each cavity is made via certain constants, which are related with the following boundary value problems:
\begin{align}\label{2.19}
&\dvr_\xi \rA(M_k^\e)\nabla_\xi X_{k,\e}=0\quad\text{in}\quad\mathds{R}^n\setminus\overline{\om_{k,\e}}, \qquad k\in\mathds{M}^\e,
\\
&X_{k,\e}=0\quad\text{on}\quad\p\om_{k,\e},\qquad k\in\mathds{M}_D^\e\cup \mathds{M}_{R,1}^\e,\label{2.27}
\\
&\nu_\xi\cdot\rA(M_k^\e)\nabla_\xi X_{k,\e}+b_k(\xi,\e) X_{k,\e}=0\quad\text{on}\quad\p\om_{k,\e},\qquad k\in \mathds{M}_{R,2}^\e,\label{2.28}
\\
& \label{2.20}
X_{k,\e}(\xi)=\left\{
\begin{aligned}
&1+K_{k,\e} |\rA^{-\frac{1}{2}}(M_k^\e)\xi|^{-n+2} +
O\big(|\xi|^{-n+1}\big), && \xi\to\infty\quad \text{as}\quad n\geqslant 3,
\\
&\ln |\rA^{-\frac{1}{2}}(M_k^\e)\xi|+K_{k,\e}+O\big(|\xi|^{-1}\big), && \xi\to\infty\quad \text{as}\quad n=2,
\end{aligned}
\right.
\end{align}
where $K_{k,\e}$ are some constants and $\nu_\xi$ is the unit normal to to $\p\om_{k,\e}$ directed inside $\om_{k,\e}$. We shall show in Lemma~\ref{lm6.1} that these problems are uniquely solvable and have classical solutions belonging to $C^2(\mathds{R}^d\setminus\om_{k,\e})\cap C^\infty(\mathds{R}^d\setminus\overline{\om_{k,\e}})$. As $n\geqslant 3$, the constants $K_{k,\e}$ can be treated as certain capacities of the cavities $\om_{k,\e}$.

We introduce an auxiliary function:
\begin{equation}\label{2.9}
\begin{aligned}
&\b_\e(x):=\frac{(2-n)K_{k,\e}}{R_3^n \mes_n B_1(0)}\quad && \text{on}\quad B_{\e R_3}(M_k^\e),\quad k\in\mathds{M}^\e,\quad\text{as}\quad n\geqslant 3,
\\
&\b_\e(x):=\frac{1}{R_3^2 \mes_n B_1(0)} && \text{on}\quad B_{\e R_3}(M_k^\e),\quad k\in\mathds{M}^\e,\quad\text{as}\quad n=2,
\\
& \b_\e(x):=0 && \text{on}\quad \Om\setminus \bigcup\limits_{k\in\mathds{M}^\e}
B_{\e R_3}(M_k^\e),
 \end{aligned}
\end{equation}
where $\mes_n$ stands for the Lebesgue measure in $\mathds{R}^n$.
It will be shown, see Lemma~\ref{lm4.2}, that the constants $K_{k,\e}$ for $n\geqslant 3$ are bounded uniformly in $k$ and $\e$. Then it follows from the above definition that the family of functions  $\b_\e$ belongs to $L_\infty(\Om)$ and is   bounded uniformly in $\e$ in this space.  Our third assumption says that the function $\b_\e$ converges to some limit $\b$ as $\e\to+0$ in an appropriate space of multipliers; this is how the strange term in (\ref{2.11}) appears and how the cavities contributes to this term. The mentioned space of multipliers is denoted by $\mathfrak{M}$ and this is the space of the functions $F$ defined on $\Om$ such that for each $u\in \Ho_2^1(\Om)\cap W_2^2(\Om)$ the function $F u$ is a continuous antilinear functional on $\Ho_2^1(\Om)$; here $\Ho_2^1(\Om)$ is the space of the functions from $W_2^1(\Om)$ with the zero trace on $\p\Om$.  The norm in $\mathfrak{M}$ is introduced as
\begin{equation}\label{2.3}
\|F\|_{\mathfrak{M}}
 =
\sup\limits_{\substack{u\in \Ho_2^1(\Om)\cap W_2^2(\Om)
\\
v\in \Ho_2^1(\Om)
}}  \frac{|\la F u,v \ra|
}{\|u\|_{W_2^2(\Om)}\|v\|_{W_2^1(\Om)}},
\end{equation}
where $\la F u,v\ra$ stands for the action of the functional $F u$ on a function $v$. The space $L_\infty(\Om)$ is a subset of $\mathfrak{M}$ due to to the identity $\la F u,v \ra=(F u,v)_{L_2(\Om)}$ for $F\in L_\infty(\Om)$.

Our third assumption reads as follows.

\begin{enumerate}
\def\theenumi{{A\arabic{enumi}}}
\setcounter{enumi}{2}

\item\label{A7} The family of the functions $\b^\e$ converges in $\mathfrak{M}$.
\end{enumerate}

This assumption means that there exists   a function $\b\in \mathfrak{M}$ such that $\|\b^\e-\b\|_{\mathfrak{M}}\to0$ as $\e\to+0$. Since the function $\b_\e$ is determined by the distribution of the points $M_k^\e$ and also by the shapes of the cavities as $n\geqslant 3$, this assumption describes the class of   non-periodic perforations, which we can consider.  We shall discuss the convergence in the space $\mathfrak{M}$ as well as possible examples of the perforations in a separate Section~\ref{secAS7}. In particular, it will be shown that the limit $\b$ is necessary an element of the space $L_\infty(\Om)$.

\subsection{Main results}

Here we formulate our main results. They involve a special boundary corrector generated by the cavities and this corrector is introduced as follows. We denote
\begin{equation*}
E_{k,\e}:=\big\{\xi: \ |\rA^{-\frac{1}{2}}(M_k^\e)\xi|<R_4\eta^{-1} \big\},
\qquad E_k^\e:=\big\{x: \ |\rA^{-\frac{1}{2}}(M_k^\e)(x-M_k^\e)|<R_4\e  \big\},
\end{equation*}
where $R_4>0$ is some fixed positive constant independent of $\e$, $k$ and $\eta$ such that
\begin{equation*}
\om_k^\e\subset B_{\e\eta R_2}(M_k^\e)\subset E_k^\e \subset B_{\e R_3}(M_k^\e).
\end{equation*}
Convergence (\ref{2.14}) yields that $\eta(\e)\to0$ as $\e\to+0$ and in view of Assumption~\ref{A1} and conditions (\ref{2.5}), (\ref{2.5a})   such constant $R_4$ obviously exists. We consider one more family of boundary value problem similar to (\ref{2.19}), (\ref{2.27}), (\ref{2.28}), (\ref{2.20}):
\begin{align}\label{2.34}
&\dvr_\xi \rA(M_k^\e)\nabla_\xi Z_{k,\e}=0\quad\text{in} \quad E_{k,\e}\setminus\overline{\om_{k,\e}},
\qquad k\in\mathds{M}^\e,
\\
&Z_{k,\e}=\left\{
\begin{aligned}
&1+K_{k,\e} R_4^{-n+2}\eta^{n+2} &&\text{as}\quad n\geqslant 3,
\\
&|\ln\eta|+\ln R_4 + K_{k,\e} &&\text{as}\quad n=2,
\end{aligned}\qquad\text{on}\quad\p E_{k,\e},
\right. \label{2.35}
\end{align}
with boundary conditions (\ref{2.27}), (\ref{2.28}) on $\p\om_{k,\e}$. We shall show in Lemma~\ref{lm3.2} that  these problems are uniquely solvable and possess classical solutions belonging to $C^2(\overline{E_{k,\e}}\setminus\om_{k,\e})\cap C^\infty(E_{k,\e}\setminus\overline{\om_{k,\e}})$. The aforementioned corrector is introduced as follows:
\begin{align}\label{2.31}
&\Xi_\e(x):=\left\{
\begin{aligned}
& \frac{Z_{k,\e}\big((x-M_k^\e)\e^{-1}\eta^{-1}\big)}
{1+K_{k,\e} R_4^{-n+2}\eta^{n-2}} &&\text{in}\quad E_k^\e\setminus\om_k^\e,\quad k\in\mathds{M}^\e,\quad n\geqslant 3,
\\
& \frac{Z_{k,\e}\big((x-M_k^\e)\e^{-1}\eta^{-1}\big)}
{|\ln\eta|+\ln R_4+K_{k,\e}} &&\text{in}\quad E_k^\e\setminus\om_k^\e,\quad k\in\mathds{M}^\e,\quad n=2,
\\
& \hphantom{X_{k,\e}\big(R_4\eta^{-1}(}1 &&\text{in}\quad \Om\setminus \bigcup\limits_{k\in\mathds{M}^\e} E_k^\e.
\end{aligned}
\right.
\end{align}

Now we are in position to formulate our main result.

\begin{theorem}\label{th1}
Let Assumption~\ref{A1} and (\ref{2.14}) be satisfied. In the case  $\mathds{M}_R^\e\ne\emptyset$ suppose also that Assumption~\ref{A6}
holds true and if $\g\ne0$, let Assumption~\ref{A7} hold  as well. Then there exists a fixed $\l_0\in\mathds{R}$ independent of $\e$ such that as $\RE\l\leqslant \l_0$, problems (\ref{2.7}), (\ref{2.11}) are uniquely solvable for each $f\in L_2(\Om)$ and the solutions  satisfy the estimates:
\begin{align}\label{2.16}
&\|u_\e-u_0\Xi_\e\|_{W_2^1(\Om^\e)}  \leqslant  C\Big(\e +\big|\eta^{n-2}\e^{-2}\vk^{-1}-\g\big| +
\g\|\b_\e-\b\|_{\mathfrak{M}} +(\e\eta\vk\mu_1)^{-\frac{1}{2}} + \e\eta\vk\mu_2
\Big)\|f\|_{L_2(\Om)},
\\
\label{2.17}
&\|u_\e-u_0\|_{L_2(\Om^\e)}\leqslant C \Big(\e +\big|\eta^{n-2}\e^{-2}\vk^{-1}-\g\big| +
\g\|\b_\e-\b\|_{\mathfrak{M}}+(\e\eta\vk\mu_1)^{-\frac{1}{2}}   + \e\eta\vk\mu_2
\Big) \|f\|_{L_2(\Om)},
\end{align}
where $C$ is some constant independent of $\e$ and $f$. If the set $\mathds{M}_R^\e$  is empty, then the terms $(\e\eta\mu_1\vk)^{-\frac{1}{2}}$  and $\e\eta\mu_2$ can be omitted in the above estimates. The terms
$\e$, $\big|\eta^{n-2}\e^{-2}\vk^{-1}-\g\big|$, $\|\b_\e-\b\|_{\mathfrak{M}}$ and $\e\eta\vk\mu_2$ in estimate  (\ref{2.16}) are order sharp. The terms  $\big|\eta^{n-2}\e^{-2}\vk^{-1}-\g\big|$ and $\|\b_\e-\b\|_{\mathfrak{M}}$
in (\ref{2.17}) are order sharp.
\end{theorem}

In a particular case $\g=0$ the assumptions on the perforation can be weakened; this case is treated in the following theorem.

\begin{theorem}\label{th2}
Let condition (\ref{2.14}) hold with $\g=0$ and Assumption~\ref{A1}
be satisfied. Then there exists a fixed $\l_0\in\mathds{R}$ independent of $\e$ such that as $\RE\l\leqslant \l_0$, problems (\ref{2.7}) and (\ref{2.11}) with $\b=0$ are uniquely solvable for each $f\in L_2(\Om)$ and the solutions  satisfy the estimates
\begin{equation}\label{2.21}
\|u_\e-u_0\|_{W_2^1(\Om^\e)}  \leqslant  C\Big(\e+\e^{-1}\eta^{\frac{n}{2}-1} \vk^{-\frac{1}{2}} + \e^{\frac{1}{2}}\eta^\frac{1}{2}
\Big)\|f\|_{L_2(\Om)},
\end{equation}
and if, in addition, $A_j\in W_\infty^1(\Om)$, then
\begin{equation}
\|u_\e-u_0\|_{L_2(\Om^\e)}  \leqslant  C\Big(\e^2+\e^{-2}\eta^{n-2} \vk^{-1} +\e\eta   \Big)\|f\|_{L_2(\Om)},\label{2.23}
\end{equation}
where $C$ are some constants independent of $\e$ and $f$.  The term  $\eta^{\frac{n}{2}-1}\e^{-1}\vk^{-\frac{1}{2}}$ in (\ref{2.21}) is order sharp.
\end{theorem}

Let us discuss briefly  the problem and main result. There are several main features of our problem. The first is that we consider a general perforation of a rather arbitrary non-periodic structure.  Assumption~\ref{A1} is  very natural and rather weak. While the second condition in (\ref{2.2}) describes the minimal distance between the points $M_k^\e$,  at the same time there are no upper bound for these distances and they can be uniformly bounded from below or even growing as $\e$ goes to zero. In particular, this means that our results also applies to the case of finitely many small cavities separated by fixed distances.

The second feature is that the boundaries of the cavities can be subject either to the Dirichlet condition or to the nonlinear Robin condition; both types of these conditions can be simultaneously present on the boundaries of the cavities. The structure of the Robin condition is described
by Assumption~\ref{A6} and these the only serious restrictions. As we have already said, conditions~(\ref{2.6}),~(\ref{2.10}) are needed only to ensure the unique solvability of problem (\ref{2.7}) and they hold immediately once we deal with the classical linear Robin condition, that is, as $a^\e(x,u)=a_\e(x)u$ with an appropriate function $a_\e(x)$.  The third feature of our model is that we consider a general second linear elliptic equation and the differential expression $\cL$ is not supposed to be formally symmetric. The coefficients $A_j$  and $A_0$ are allowed to be complex-valued.

Our main theorem states that under the above discussed conditions, the homogenized problem is (\ref{2.11}) and the convergence in $W_2^1(\Om^\e)$ and $L_2(\Om^\e)$ holds uniformly in the right hand side $f$; the estimates for the convergence rates are our main results. In the case when the Robin condition is present and is linear, these are operator estimates describing the norm resolvent convergence of the perturbed operator to the homogenized one. Estimate (\ref{2.16}) says that the solution $u_\e$ to problem (\ref{2.7}) can be approximated by $u_0\Xi_\e$ in $W_2^1(\Om^\e)$ and the estimate for the convergence rate is provided. The corrector $\Xi_\e$  can be omitted and then we have a similar result but only in  $L_2(\Om^\e)$ with the same convergence rate. The terms involving $\mu_1$ and $\mu_2$ are generated only due to the presence of the nonlinear Robin condition. If the set $\mathds{M}_R^\e$ is empty, they can be removed from the estimates.

It is also shown that all terms except for $(\e\eta\vk\mu_1)^{-\frac{1}{2}}$ are order sharp in (\ref{2.16}). In particular, this implies that the term $\|\b_\e-\b\|_{\mathfrak{M}}$ can not be omitted and hence, the same concerns Assumption~\ref{A7}.  In estimate (\ref{2.17}), two terms in the convergence rate are also order sharp.   It is unclear to us whether other terms in (\ref{2.16}), (\ref{2.17}) are also order sharp or the estimate could be improved by using some additional techniques. This question remained open.

As $\g=0$, it is possible to omit Assumptions~\ref{A6},~\ref{A7} and to prove similar results only under Assumption~\ref{A1}, see estimates (\ref{2.21}), (\ref{2.23}). We stress that in (\ref{2.21}) the corrector is absent in comparison with (\ref{2.16}) but the price we pay for this and for omitting additional assumptions is a worse convergence rate. However, estimating then in $L_2(\Om^\e)$-norm, the order of the convergence rate can be improved twice, see (\ref{2.23}). The term $\e^{-1}\eta^{\frac{n}{2}-1} \vk^{-\frac{1}{2}} $ is shown to be order sharp in (\ref{2.21}). The sharpness of the other terms in (\ref{2.21}) and of all terms in (\ref{2.23}) remains an open question.

We observe that the sharpness of the term $\e^{-1}\eta^{\frac{n}{2}-1} \vk^{-\frac{1}{2}} $ does not contradicts the sharpness of a similar term in (\ref{2.16}) since in (\ref{2.21}) we have omitted the corrector. We also stress that if $\g\ne0$, then by omitting the corrector in  (\ref{2.16}) we destroy the convergence in $W_2^1(\Om^\e)$, namely, it turns out that $\|u_\e-u_0\|_{W_2^1(\Om^\e)}$ is just of order $O(1)$ once $\g\ne0$.

\section{Auxiliary lemmata}

In this section we provide of series of auxiliary lemmata, which will be employed then in the proofs of Theorem~\ref{th1},~\ref{th2}.

\begin{lemma}\label{lm3.1}
Under Assumption~\ref{A1}
for all $k\in\mathds{M}^\e$ and all $u\in W_2^1(B_{\e R_3}(M_k^\e)\setminus\om_k^\e)$ the estimates
\begin{align}\label{3.6}
&\|u\|_{L_2(\p\om_k^\e)}^2 \leqslant  C\Big(\e\eta \vk\|\nabla u\|_{L_2(B_{\e R_3}(M_k^\e)\setminus\om_k^\e)}^2 +\e^{-1}\eta^{n-1}
\|u\|_{L_2(B_{\e R_3}(M_k^\e)\setminus B_{\e R_2}(M_k^\e))}^2\Big),
\\
&\label{3.71}
\|u\|_{L_2(B_{\e\eta R_2}(M_k^\e)\setminus\om_k^\e)}^2\leqslant C \Big(\e\eta\|u\|_{L_2(\p\om_k^\e)}^2 +
\e^2\eta^2\|\nabla u\|_{L_2(B_{\e R_3}(M_k^\e)\setminus\om_k^\e)}^2
\Big),
\end{align}
hold, where $C$ are constants independent of $k$, $\e$, $\eta$ and $u$.  If, in addition,
\begin{equation}\label{3.73}
\int\limits_{B_{\e\eta R_2}(M_k^\e)\setminus \om_k^\e} u\,dx=0,
\end{equation}
then the estimate
\begin{equation}\label{3.72}
\|u\|_{L_2(B_{\e\eta R_2}(M_k^\e)\setminus\om_k^\e)}^2\leqslant C\e^2\eta^2 \|\nabla u\|_{L_2(B_{\e R_3}(M_k^\e)\setminus\om_k^\e)}^2
\end{equation}
holds,  where $C$ is a constant independent of $k$, $\e$, $\eta$ and $u$.
\end{lemma}

\begin{proof}
Inequality (\ref{3.6}) was proved in  \cite[Lm. 3.6]{BK}, while inequality (\ref{3.72}) was established in \cite[Lm. 3.5]{BK}. Given an arbitrary $k\in\mathds{M}^\e$ and $u\in W_2^1(B_{\e\eta R_2}(M_k^\e)\setminus\om_k^\e)$, we denote
\begin{equation}\label{3.78}
\la u\ra_k:=\frac{1}{|B_{\e\eta R_3}(M_k^\e)\setminus\om_k^\e|} \int\limits_{B_{\e\eta R_3}(M_k^\e)\setminus \om_k^\e} u\,dx, \qquad
u_k:=u-\la u \ra_k.
\end{equation}
The function $u_k$ satisfies condition (\ref{3.73}) and by (\ref{3.72})  we have
\begin{equation*}
\|u_k\|_{L_2(B_{\e\eta R_3}(M_k^\e)\setminus\om_{k,\e})}^2 \leqslant C\e^2\eta^2 \|\nabla u\|_{L_2(B_{\e\eta R_3}(M_k^\e) \setminus\om_{k,\e})}^2;
\end{equation*}
hereinafter in the proof we denote by $C$ inessential constants independent of $\e$, $\eta$, $k$ and $v$. By \cite[Lm. 3.5]{BK} we also have:
\begin{equation}\label{3.74}
\|u_k\|_{L_2(\p\om_k^\e)}^2\leqslant C\e\eta \|\nabla u\|_{L_2(B_{\e\eta R_3}\setminus\om_k^\e)}^2.
\end{equation}

Employing the first condition in (\ref{2.2}), we argue as follows:
\begin{align*}
\int\limits_{\p\om_k^\e}   |u|^2 \,ds= & |\la u\ra_k|^2 \mes_{n-1}\p\om_k^\e
+ 2 \RE \overline{\la u\ra_k} \int\limits_{\p\om_k^\e}
u_k\,ds + \int\limits_{\p\om_k^\e}  |u_k|^2 \,ds
\\
\geqslant &   (R_1\e\eta)^{n-1}|\la u\ra_k|^2  \mes_{n-1} \p B_1(0) -2|\la u\ra_k| \int\limits_{\p\om_k^\e} |u_k|\,ds
\\
\geqslant & \frac{1}{2}(R_1\e\eta)^{n-1} |\la u\ra_k|^2 \mes_{n-1} \p B_1(0) -C(\e\eta)^{-n+1}\Bigg(\int\limits_{\p\om_k^\e}   |u_k|\,ds\Bigg)^2
\\
\geqslant & \frac{1}{2}(R_1\e\eta)^{n-1}  |\la u\ra_k|^2 \mes_{n-1} \p B_1(0) -C
\|u_k\|_{L_2(\p\om_k^\e)}^2.
\end{align*}
Hence, by (\ref{3.74}),
\begin{equation*}
(\e\eta)^n |\la u \ra_k|^2 \leqslant C \Big( \e\eta \|u\|_{L_2(\p\om_k^\e)}^2 \,ds + \e^2\eta^2 \|\nabla u\|_{L_2(B_{\e\eta R_3}\setminus\om_k^\e)}^2\Big).
\end{equation*}
Using this estimate and (\ref{3.74}), we obtain:
\begin{align*}
\|u\|_{L_2(B_{\e\eta R_2}(M_k^\e)\setminus\om_k^\e)}^2\leqslant  & 2 |\la u\ra_k |^2 \mes_n (B_{\e\eta R_2}(M_k^\e)\setminus\om_k^\e) + 2 \|u_k\|_{L_2(B_{\e\eta R_3}\setminus\om_k^\e)}^2
\\
\leqslant  & C \Big( \e\eta\|u\|_{L_2(\p\om_k^\e)}^2  + \e^2\eta^2 \|\nabla u\|_{L_2(B_{\e\eta R_3}\setminus\om_k^\e)}^2\Big).
\end{align*}
This proves (\ref{3.71}).
\end{proof}

We note that under  convergence (\ref{2.14}) we have
\begin{equation}\label{3.23}
\eta^{n-1}=\e^2\eta\vk\cdot \e^{-2}\eta^{n-2}\vk^{-1}\leqslant C \e^2\eta \vk,
\end{equation}
where $C$ is some fixed constant independent of $\e$ and $\eta$. Then   estimate (\ref{3.6}) can be rewritten as
\begin{equation}\label{3.22}
\|u\|_{L_2(\p\om_k^\e)}^2 \leqslant  C \e\eta \vk \| u\|_{W_2^1(B_{\e R_3}(M_k^\e)\setminus\om_k^\e)}^2.
\end{equation}

\begin{lemma}\label{lm3.6}
Under Assumption~\ref{A1}  for all $k\in\mathds{M}_R^\e$ and all $u\in W_2^1(B_{\e R_3}(M_k^\e)\setminus\om_k^\e)$ the estimate
\begin{equation*}
\|u\|_{L_2(B_{\e\eta R_3}(M_k^\e)\setminus\om_k^\e)}^2 \leqslant C  \Big(\e^2\eta^2 \vk\|\nabla u\|_{L_2(B_{\e R_3}(M_k^\e)\setminus\om_k^\e)}^2 + \eta^n
\|u\|_{L_2(B_{\e R_3}(M_k^\e)\setminus\om_k^\e)}^2\Big)
\end{equation*}
holds with a constant $C$ independent of $k$, $\e$ and $u$.
\end{lemma}

\begin{proof}
For each $u\in W_2^1(B_{\e R_3}(M_k^\e)\setminus\om_k^\e)$ we integrate by parts as follows:
\begin{align*}
 \|u\|_{L_2(B_{\e\eta R_3}(M_k^\e)\setminus\om_k^\e)}^2 = & \frac{1}{n}\int\limits_{B_{\e\eta R_3}(M_k^\e)\setminus\om_k^\e} |u|^2\dvr x\,dx
\\
= &\frac{\e\eta R_3}{n} \int\limits_{\p B_{\e\eta R_3}(M_k^\e)} |u|^2\,ds + \frac{1}{n}\int\limits_{\p\om_k^\e} |u|^2 x\cdot \nu \,ds - \frac{1}{n}\int\limits_{B_{\e\eta R_3}(M_k^\e)\setminus\om_k^\e} x\cdot \nabla |u|^2\,dx.
\end{align*}
We estimate the integral over $B_{\e\eta R_3}(M_k^\e)\setminus\om_k^\e$ in this inequality as
\begin{align*}
 \Bigg| \int\limits_{B_{\e\eta R_3}(M_k^\e)\setminus\om_k^\e} x\cdot \nabla |u|^2\,dx \Bigg| \leqslant & 2\e\eta R_3\int\limits_{B_{\e\eta R_3}(M_k^\e)\setminus\om_k^\e} |u||\nabla u|\,dx
\\
\leqslant & \|u\|_{L_2(B_{\e\eta R_3}(M_k^\e)\setminus\om_k^\e)}^2 + \e^2\eta^2 R_3^2  \|\nabla u\|_{L_2(B_{\e\eta R_3}(M_k^\e)\setminus\om_k^\e)}^2.
\end{align*}
The integrals over $\p B_{\e\eta R_3}(M_k^\e)$ and $\p\om_k^\e$ are estimated by means of inequality (\ref{3.6})  and this finally  completes the proof.
\end{proof}

\begin{lemma}\label{lm3.3}
Under  Assumption~\ref{A1}
for all $k\in\mathds{M}_R^\e$ and all $u\in W_2^1(B_{\e  R_3}(M_k^\e)\setminus\om_k^\e)$ obeying the identity
\begin{equation}\label{3.1}
\int\limits_{B_{\e   R_3}(M_k^\e)\setminus\om_k^\e} u(x)\,dx=0
\end{equation}
the estimate
\begin{equation}\label{3.2a}
\|u\|_{L_2(B_{\e  R_3}(M_k^\e)\setminus\om_k^\e)}^2\leqslant C\e^2  \|\nabla u\|_{L_2(B_{\e   R_3}(M_k^\e)\setminus\om_k^\e)}^2,
\end{equation}
holds, where $C$ is a constant independent of the parameters $k$, $\e$, $\eta$ and the function $u$.
\end{lemma}

\begin{proof}
We first consider the Neumann Laplacian on the annulus $B_{R_3}(0)\setminus B_{\eta R_3}(0)$; we denote it by $\mathcal{D}^N_\eta$. Since $\eta\to+0$, by \cite[Thm. 1.2]{PMA18},  this operator converges in the norm resolvent sense to the Neumann Laplacian on $B_{R_3}(0)$, which we denote by $\mathcal{D}^N_0$. Namely, the estimate holds:
\begin{equation*}
\|(\mathcal{D}^N_\eta-\iu)^{-1}g-(\mathcal{D}^N_0-\iu)^{-1}g\|_{W_2^1(B_{R_3}(0)\setminus B_{\eta R_3}(0))}\leqslant C\eta^\frac{1}{2}\|g\|_{L_2(B_{R_3}(0))}
\end{equation*}
for each $g\in L_2(B_{R_3}(0))$ with a constant $C$ independent of $\eta$ and $g$. Employing this estimate and proceeding as in \cite[Sect. 7]{PRSE16}, we easily see that the spectrum of  $\mathcal{D}^N_\eta$ converges to that of $\mathcal{D}^N_0$ as $\eta\to+0$ and hence, the second eigenvalue of the operator $\mathcal{D}^N_\eta$ is positive and separated from zero uniformly in $\eta$. Therefore, by the minimax principle, for each $\tilde{u}\in W_2^1(B_{R_3}(0)\setminus B_{\eta R_3}(0))$
obeying the condition
\begin{equation*}
\int\limits_{B_{R_3}(0)\setminus B_{\eta R_3}(0)} \tilde{u}\,d\xi=0
\end{equation*}
we have the estimate
\begin{equation*}
\frac{\|\nabla \tilde{u}\|_{L_2(B_{R_3}(0)\setminus B_{\eta R_3}(0))}^2}{\|\tilde{u}\|_{L_2(B_{R_3}(0)\setminus B_{\eta R_3}(0))}^2}\geqslant C,
\end{equation*}
where $C$ is a positive constant independent of $\tilde{u}$ and $\eta$. Given then an arbitrary function $u\in W_2^1(B_{\e R_3}(M_k^\e)\setminus B_{\e\eta R_3}(M_k^\e))$ such that
\begin{equation}\label{3.8}
\int\limits_{B_{\e R_3}(M_k^\e)\setminus B_{\e\eta R_3}(M_k^\e)} u\,dx=0
\end{equation}
 and applying the above inequality to $\tilde{u}(\xi):=u(M_k^\e+\e\xi)$, we obtain:
\begin{equation}\label{3.7}
\|u\|_{L_2(B_{\e R_3}(M_k^\e)\setminus B_{\e\eta R_3}(M_k^\e))}^2\leqslant C\e^2
\|\nabla u\|_{L_2(B_{\e R_3}(M_k^\e)\setminus B_{\e\eta R_3}(M_k^\e))}^2,
\end{equation}
where $C$ is a positive constant independent of $u$, $k$, $\e$ and $\eta$.

Let $u\in W_2^1(B_{\e R_3}(M_k^\e)\setminus\om_k^\e)$ be an arbitrary function obeying condition (\ref{3.1}). In $B_{\e R_3}(M_k^\e)\setminus B_{\e\eta R_3}(M_k^\e)$ we represent it as
\begin{equation}\label{3.9}
u=u_\dag +u_\ddag,\qquad u_\dag:=\frac{1}{\mes_n B_{\e R_3}(M_k^\e)\setminus B_{\e\eta R_3}(M_k^\e)} \int\limits_{B_{\e R_3}(M_k^\e)\setminus B_{\e\eta R_3}(M_k^\e)} u\,dx,\qquad u_\ddag:=u-u_\dag.
\end{equation}
The function $u_\ddag$ obviously satisfies condition (\ref{3.8}) and hence, inequality (\ref{3.7}), while due to (\ref{3.1}) and the Cauchy-Schwarz inequality for the constant $u_\dag$ we have
\begin{equation*}
|u_\dag|^2=\Bigg|\frac{1}{\mes_n B_{\e R_3}(M_k^\e)\setminus B_{\e\eta R_3}(M_k^\e)}
\int\limits_{B_{\e\eta R_3}(M_k^\e)\setminus\om_k^\e} u\,dx\Bigg|^2\leqslant C\eta^n\e^{-n} \|u\|_{L_2(B_{\e\eta R_3}(M_k^\e)\setminus\om_k^\e)}^2,
\end{equation*}
where $C$ is a constant independent of $u$, $k$, $\e$ and $\eta$. Then by the above estimate, inequality  (\ref{3.7}) for $u_\ddag$, Lemma~\ref{lm3.6} and the convergence $\eta\to0$ we have:
\begin{align*}
\|u\|_{L_2(B_{\e R_3}(M_k^\e)\setminus\om_k^\e)}^2\leqslant & C\big( \e^n |u_\dag|^2 +\|u_\ddag\|_{L_2(B_{\e R_3}(M_k^\e)\setminus B_{\e\eta R_3}(M_k^\e))}^2\big)+\|u\|_{L_2(B_{\e\eta R_3}(M_k^\e)\setminus\om_k^\e)}^2
\\
 \leqslant & C
\|u\|_{L_2(B_{\e\eta R_3}(M_k^\e)\setminus\om_k^\e)}^2 + C \e^2 \|\nabla u\|_{L_2(B_{\e R_3}(M_k^\e)\setminus B_{\e\eta R_3}(M_k^\e))}^2
\\
 \leqslant & C\eta^n
\|u\|_{L_2(B_{\e\eta R_3}(M_k^\e)\setminus\om_k^\e)}^2 + C \e^2 \|\nabla u\|_{L_2(B_{\e R_3}(M_k^\e)\setminus \om_k^\e)}^2,
\end{align*}
where $C$ are some fixed constants independent of $u$, $k$, $\e$ and $\eta$. The obtained estimates imply (\ref{3.2a}).  The proof is complete.
\end{proof}

We recall that a generalized solution to problem (\ref{2.7})
 is  a function $u\in \Ho_2^1(\Om^\e,\p\Om\cup\p\tht_D^\e)$  satisfying the integral identity
\begin{equation}\label{2.18}
\hf(u,v) + (a^\e(\,\cdot\,,u),v)_{L_2(\p\tht_R^\e)}-\l(u_\e,v)_{L_2(\Om^\e)}=(f,v)_{L_2(\Om^\e)}
\end{equation}
for each $v\in \Ho_2^1(\Om^\e,\p\Om\cup\p\tht_D^\e)$, where
\begin{align*}
\hf(u,v):=
(\rA \nabla u,\nabla v)_{L_2(\Om^\e)}
 + \sum\limits_{j=1}^{n} \left(A_j \frac{\p u}{\p x_j}, v\right)_{L_2(\Om^\e)}
  + (A_0 u,v)_{L_2(\Om^\e)}.
\end{align*}
A generalized solution to problem (\ref{2.11}) is defined in a similar way. The next lemma ensures the unique solvability of problems (\ref{2.7}), (\ref{2.11}).

\begin{lemma}\label{lm4.1}
Under Assumption~\ref{A1} and also under Assumption~\ref{A6} if $\mathds{M}_R^\e\ne\emptyset$, there exists $\l_0\in\mathds{R}$ independent of $\e$ such that for $\RE\l<\l_0$ problems (\ref{2.7}), (\ref{2.11}) are uniquely solvable  for each $f\in L_2(\Om^\e)$ and for their solutions the estimates
\begin{align}\label{2.12}
&\|u_\e\|_{W_2^1(\Om^\e)}\leqslant C(\l) \|f\|_{L_2(\Om^\e)},
\\
&\|u_0\|_{W_2^2(\Om)}\leqslant C(\l)\|f\|_{L_2(\Om)} \label{2.15}
\end{align}
hold,
where $C(\l)$ are some constants independent of $\e$ and $f$.
\end{lemma}

\begin{proof}
Assumption~\ref{A6} guarantees that  estimate (2.5) in \cite{BK} is satisfied and hence, by Lemma~3.7 in \cite{BK},
there exists $\l_0\in\mathds{R}$ independent of $\e$ such that for $\RE\l<\l_0$ problem (\ref{2.7}) is solvable in $\Ho_2^1(\Om^\e,\p\Om\cup\p\tht_D^\e)$ for each $f\in L_2(\Om^\e)$. This is why we just need to check the uniqueness of the solution.

Supposing that there are two solutions $u_\e$ and $\tilde{u}_\e$ for some $f$, the difference $\hat{u}_\e:=u_\e-\tilde{u}_\e$ then solves the boundary value problem
\begin{equation*}
(\cL-\l)\hat{u}_\e=0\quad\text{in}\quad \Om^\e,\qquad \hat{u}_\e=0\quad\text{on} \quad \p\Om\cup\p\tht_D^\e,\qquad  \frac{\p \hat{u}_\e}{\p\boldsymbol{\nu}} + a^\e(x,u_\e)-a^\e(x,\tilde{u}_\e)=0\quad\text{on}\quad \p\tht_R^\e.
\end{equation*}
Writing  the corresponding integral identity with $\hat{u}_\e$ as the test function, we immediately get:
\begin{equation}\label{3.3}
\hf(\hat{u}_\e,\hat{u}_\e) - \l \|\hat{u}_\e\|_{L_2(\Om^\e)}^2 + \big(a^\e(\,\cdot\,,u_\e)-a^\e(\,\cdot\,,\tilde{u}_\e), \hat{u}_\e\big)_{L_2(\p\tht_R^\e)}=0.
\end{equation}
Thanks to the lower bound in the two-sided inequality in (\ref{2.6}) and also to (\ref{3.6}),  the second term in the above identity satisfies the estimate:
\begin{equation}\label{3.4}
\begin{aligned}
\RE\big(a^\e(\,\cdot\,,u_\e)-a^\e(\,\cdot\,,\tilde{u}_\e), \hat{u}_\e\big)_{L_2(\p\tht_R^\e)}\geqslant & -\mu_0(\e) \|u_1-u_2\|_{L_2(\p\tht_R^\e)}^2
\\
\geqslant & - C\mu_0 \sum\limits_{k\in\mathds{M}_R^\e} \Big(\e\eta \vk\|\nabla u\|_{L_2(B_{\e R_3}(M_k^\e)\setminus\om_k^\e)}^2
\\
&\hphantom{- C\mu_0 \sum\limits_{k\in\mathds{M}_R^\e} \Big(}+\e^{-1}\eta^{n-1}
\|u\|_{L_2(B_{\e R_3}(M_k^\e)\setminus B_{\e R_2}(M_k^\e))}^2\Big)
\\
\geqslant & C \Big(\mu_0\e\eta \vk\|\nabla u\|_{L_2(\Om^\e)}^2 +\mu_0\e^{-1}\eta^{n-1}
\|u\|_{L_2(\Om^\e)}^2\Big).
\end{aligned}
\end{equation}
By (\ref{2.10}) and   (\ref{2.14}) we have:
\begin{equation}\label{3.5}
\mu_0\e^{-1}\eta^{n-1}=\mu_0\e\eta \vk\cdot \vk^{-2} \cdot \e^{-2} \eta^{n-2} \vk^{-1} \to +0,\qquad \e\to+0.
\end{equation}
Conditions (\ref{2.5}), (\ref{2.5a})  and the Cauchy-Schwarz inequality imply that
\begin{equation}\label{4.3}
\RE \hf(u,u)\geqslant \frac{3c_0}{4} \|\nabla u\|_{L_2(\Om^\e)}^2 - C \|u\|_{L_2(\Om^\e)}^2
\end{equation}
for all $u\in W_2^1(\Om^\e)$, where $C$ is some absolute constant independent of $\e$ and $u\in W_2^1(\Om^\e)$. This estimate and (\ref{3.3}), (\ref{3.4}), (\ref{3.5}) then yield
\begin{align*}
0=&\RE\hf(\hat{u}_\e,\hat{u}_\e) - \RE \l \|\hat{u}_\e\|_{L_2(\Om^\e)}^2
+ \RE \big(a^\e(\,\cdot\,,u_\e)-a^\e(\,\cdot\,,\tilde{u}_\e), \hat{u}_\e\big)_{L_2(\p\tht_R^\e)}
\\
\geqslant & \frac{c_0}{2} \|\nabla \hat{u}_\e\|_{L_2(\Om^\e)}^2 -(\RE \l + C) \|\hat{u}_\e\|_{L_2(\Om^\e)}^2
\end{align*}
for $\e$ small enough, where $C$ is some fixed constant independent of $\e$ and $\hat{u}_\e$.  Hence, as $\RE \l<-C-1$, we necessarily have $\hat{u}_\e=0$ and this proves the uniqueness of solution to problem (\ref{2.7}). Writing the integral identity corresponding to problem (\ref{2.7}) with $u_\e$ as the test function
 and proceeding as above, we prove easily estimate (\ref{2.12}).

Problem (\ref{2.11}) can be treated as a resolvent equation for the operator generated by the differential expression $\cL+\g\Ups\b$ in $L_2(\Om)$ subject to the Dirichlet condition. Conditions~(\ref{2.5}),~(\ref{2.5a}) imply easily that such operator is $m$-sectorial and this is why the unique solvability is just a standard fact from the theory of $m$-sectorial operators. This operator is bounded as that from $W_2^2(\Om)\cap \Ho_2^1(\Om)$ into $L_2(\Om)$ and this is why, by the Banach theorem, its resolvent is a bounded operator from $L_2(\Om)$ into $W_2^2(\Om)\cap \Ho_2^1(\Om)$. This implies estimate (\ref{2.15}). The proof is complete.
\end{proof}

\section{Properties of corrector}

In this section we prove the solvability of problems (\ref{2.19}), (\ref{2.27}), (\ref{2.28}), (\ref{2.20}), (\ref{2.34}), (\ref{2.35}) and study certain properties of their solutions. Throughout this section we suppose that Assumption~\ref{A1} is satisfied.

The main point in this study is an appropriate Kelvin transform reducing the problems (\ref{2.19}), (\ref{2.27}), (\ref{2.28}), (\ref{2.20}) to ones in bounded domains. 
This transform is defined as
\begin{equation}\label{3.79}
\tilde{\xi}:=\frac{\rA^{-\frac{1}{2}}(M_k^\e)(\xi-y_{k,\e})} {|\rA^{-\frac{1}{2}}(M_k^\e)(\xi-y_{k,\e})|^2},
\end{equation}
with the points $y_{k,\e}$ introduced in Assumption~\ref{A1}. By
$\tilde{\om}_{k,\e}$ we denote the images of the domains $\om_{k,\e}$ arising while passing to the variable $\tilde{\xi}$.
It is clear that the domains $\tilde{\om}_{k,\e}$ are unbounded, namely, $\tilde{\om}_{k,\e}\supset \mathds{R}^d\setminus B_{R_1^{-1}}(0)$,  the boundaries of these domains are smooth and the domains $\mathds{R}^d\setminus \tilde{\om}_{k,\e}$ are bounded.
It is also easy to see that the origin does not belong  to $\tilde{\om}_{k,\e}$ and moreover, $B_{(2R_2)^{-1}}(0)\cap\tilde{\om}_{k,\e}=\emptyset$.

We seek a solution to problem (\ref{2.19}), (\ref{2.27}), (\ref{2.28}), (\ref{2.20}) as
\begin{equation}\label{3.81}
X_{k,\e}(\xi)=\frac{1}{|\rA^{-\frac{1}{2}}(M_k^\e)(\xi-y_{k,\e})|^{n-2}} \tilde{X}_{k,\e} \left(
\frac{\rA^{-\frac{1}{2}}(M_k^\e)(\xi-y_{k,\e})} {|\rA^{-\frac{1}{2}}(M_k^\e)(\xi-y_{k,\e})|^2}\right) + \left\{
\begin{aligned}
&1 &&\text{as}\quad n\geqslant 3,
\\
\ln|\rA^{-\frac{1}{2}}(M_k^\e&)(\xi-y_{k,\e})|  &&\text{as}\quad n=2.
\end{aligned}
\right.
\end{equation}
Then for the functions $\tilde{X}_{k,\e}$ we obtain the boundary values problems
\begin{equation}\label{3.85}
\begin{aligned}
&\D_{\tilde{\xi}} \tilde{X}_{k,\e}=0\quad\text{in}\quad \mathds{R}^d\setminus \tilde{\om}_{k,\e},\qquad k\in\mathds{M}^\e,
\\
&\tilde{X}_{k,\e}=\tilde{\phi}_{k,\e}(\xi)\quad\text{on}\quad \p\tilde{\om}_{k,\e},\qquad  k\in\mathds{M}_D^\e\cup \mathds{M}_{R,1}^\e,
\\
\bigg(&\frac{\p\ }{\p\tilde{\nu}} + \tilde{b}_{k,\e}(\tilde{\xi})
\bigg)\tilde{X}_{k,\e}=\tilde{\phi}_{k,\e}(\xi)\quad\text{on}\quad \p\tilde{\om}_{k,\e},\qquad  k\in \mathds{M}_{R,2}^\e,
\end{aligned}
\end{equation}
where $\tilde{b}_{k,\e}$
and $\tilde{\phi}_{k,\e}$ are some complex-valued functions. These functions are the elements of the following spaces:
\begin{equation*}
\tilde{\phi}_{k,\e}\in C^2(\p\tilde{\om}_{k,\e}),\quad k\in\mathds{M}_D^\e\cup M_{R,1}^\e,\qquad
\tilde{\phi}_{k,\e}, \tilde{b}_{k,\e}\in C^1(\p\tilde{\om}_{k,\e}),\quad k\in\mathds{M}_{R,2}^\e,
\end{equation*}
and they are bounded uniformly in $k$ and $\e$ in the norms of these spaces. The functions $\tilde{b}_{k,\e}$ also satisfy the estimate
\begin{equation}\label{3.10}
\RE \tilde{b}_{k,\e}\geqslant \tilde{c}_2,
\end{equation}
where $\tilde{c}_2$ is some fixed positive constant independent of $k$ and $\e$, while $\tilde{\nu}$ is the  normal to $\p\tilde{\om}_{k,\e}$ directed inside $\tilde{\om}_{k,\e}$.

\begin{lemma}\label{lm6.1}
Under Assumption~\ref{A1} and conditions (\ref{2.5}), (\ref{2.5a}), (\ref{2.26}) problems (\ref{2.19}), (\ref{2.27}), (\ref{2.28}), (\ref{2.20})
are uniquely solvable in $C^\infty(\mathds{R}^d\setminus\overline{\om_{k,\e}})\cap W_2^2(B_{R_3}(0)\setminus \om_{k,\e})$.
\end{lemma}

\begin{proof}
We treat problems (\ref{3.85}) in the generalized sense seeking their solutions in $W_2^1(\mathds{R}^d\setminus \tilde{\om}_{k,\e})$.  We consider homogeneous problems (\ref{3.85}) with $\tilde{\phi}_{k,\e}=0$,
write the corresponding integral identities and use inequality (\ref{3.10}) for $k\in\mathds{M}_{R,2}^\e$. Then we see easily that these homogeneous problems can have only trivial solutions. Hence, problems (\ref{3.85}) are uniquely solvable in $W_2^1(\mathds{R}^d\setminus \tilde{\om}_{k,\e})$. By standard smoothness improving estimates we immediately conclude that these functions are the elements of $W_2^2(\mathds{R}^d\setminus \tilde{\om}_{k,\e})$ and are infinitely differentiable in $\mathds{R}^d\setminus \overline{\tilde{\om}_{k,\e}}$.  Since the functions $\tilde{X}_{k,\e}$ are infinitely differentiable in the vicinity of the origin and are represented there by its Taylor series; in particular,
\begin{equation*}
\tilde{X}_{k,\e}(\tilde{\xi})=K_{k,\e}+O(|\tilde{\xi}|),\quad \tilde{\xi}\to0,
\end{equation*}
where $K_{k,\e}$ are some constants. Recovering then functions $X_{k,\e}$ by formulae (\ref{3.81}), we complete the proof.
\end{proof}

\begin{lemma}\label{lm6.5}
For all $k\in\mathds{M}^\e$ the functions $\tilde{X}_{k,\e}$
belong to $L_\infty(\mathds{R}^d\setminus\tilde{\om}_{k,\e})$ and satisfy the uniform  estimates
\begin{equation}\label{3.86}
\|\tilde{X}_{k,\e}\|_{L_\infty(\mathds{R}^d\setminus\tilde{\om}_{k,\e})} \leqslant C,
\end{equation}
where $C$ is some constant independent of $k$ and $\e$.
\end{lemma}

\begin{proof}
As $k\in\mathds{M}_D^\e\cup \mathds{M}_{R,1}^\e$, by the weak maximum principle \cite[Ch. 8, Sect. 8.1, Thm. 8.1]{GT} applied to the real and imaginary parts of the function $\tilde{X}_{k,\e}$ and by the uniform boundedness of the functions $\tilde{\phi}_{k,\e}$ in $C(\p\om_{k,\e})$ we immediately get the statement of the lemma for such $k$.

The case $k\in\mathds{M}_{R,2}^\e$ requires a more detailed study. We first state that for all $u\in W_2^1(\mathds{R}^d\setminus \tilde{\om}_{k,\e})$ and all $k\in\mathds{M}^\e$ the estimate holds
\begin{equation}\label{3.15}
\|u\|_{L_2(\p\tilde{\om}_{k,\e})}\leqslant C \|u\|_{W_2^1(\mathds{R}^d\setminus \tilde{\om}_{k,\e})},
\end{equation}
where $C$ is some constant independent of $u$, $\e$ and $k\in\mathds{M}^\e$. This is implied by a similar estimate for $u\in B_{R_3}(0)\setminus\om_{k,\e}$ established in the proof of Lemma~3.5 in \cite{BK}. We write the integral identity corresponding to problem (\ref{3.85}) with $\tilde{X}_{k,\e}$ as the test function and take then the real part of this identity. In view of (\ref{3.10}), (\ref{3.15}) and the uniform boundedness of $\tilde{\phi}_{k,\e}$ this gives the estimate
\begin{equation}\label{3.14}
\|\nabla_{\tilde{\xi}}\tilde{X}_{k,\e} \|_{L_2(\mathds{R}^d\setminus\tilde{\om}_{k,\e})}^2 + \|\tilde{X}_{k,\e}\|_{L_2(\p\tilde{\om}_{k,\e})}^2\leqslant C;
\end{equation}
hereinafter till the end of the proof by $C$ we denote inessential constants independent of $k$ and $\e$.

Now we use the technique from \cite[Ch. I\!I\!I, Sect. 13]{Ld}. We choose an arbitrary $\vr>0$ and write the integral identity for problem (\ref{3.85}) with the test function $\tilde{X}_{k,\e}$, $Q_{\min,\vr}:=\min\{|\tilde{X}_{k,\e}|^2,\vr\}$. Taking then the real part of the obtained identity and using (\ref{3.10}) and the uniform boundedness of $\tilde{\phi}_{k,\e}$, after some simple arithmetical calculations we get:
\begin{align*}
\int\limits_{\mathds{R}^d\setminus\tilde{\om}_{k,\e}} \Bigg(|& \nabla_{\tilde{\xi}}\tilde{X}_{k,\e}|^2 Q_{\min,\vr} +\frac{1}{2}|\nabla_{\tilde{\xi}}|\tilde{X}_{k,\e}||^2\Bigg)\,dx + \tilde{c}_2\int\limits_{\p\tilde{\om}_{k,\e}} |\tilde{X}_{k,\e}|^2 Q_{\min,\vr}\,ds
\\
&\leqslant  C\int\limits_{\p\tilde{\om}_{k,\e}} |\tilde{X}_{k,\e}|Q_{\min,\vr}\,ds
\leqslant  \frac{\tilde{c}_2}{2} \int\limits_{\p\tilde{\om}_{k,\e}} |\tilde{X}_{k,\e}|^2 Q_{\min,\vr}\,ds + C \int\limits_{\p\tilde{\om}_{k,\e}} Q_{\min,\vr}\,ds,
\end{align*}
where the constants $C$ are independent of $\vr$,  $\tilde{X}_{k,\e}$, $k$ and $\e$. Passing to the limit as $\vr\to+\infty$ and using (\ref{3.14}),
we get
\begin{equation*}
\int\limits_{\mathds{R}^d\setminus\tilde{\om}_{k,\e}} \Bigg(|\nabla_{\tilde{\xi}}\tilde{X}_{k,\e}|^2 |\tilde{X}_{k,\e}|^2 +|\nabla_{\tilde{\xi}}|\tilde{X}_{k,\e}|^2|^2\Bigg)\,dx + \int\limits_{\p\tilde{\om}_{k,\e}} |\tilde{X}_{k,\e}|^4\,ds\leqslant C,
\end{equation*}
where $C$ is a constant independent of $\e$ and $k$.

We once again choose an arbitrary $\vr>0$ and write the integral identity
for problem (\ref{3.85}) with the test function $\tilde{X}_{k,\e} Q_{\max,\vr}$, $Q_{\max,\vr}:=\max\{|\tilde{X}_{k,\e}|^2-\vr,0\}$, taking then the real part of the obtained identity. After simple estimates in the integrals over $\p\tilde{\om}_{k,\e}$ in this identity, in view of the uniform boundedness of $\tilde{\phi}_{k,\e}$ and (\ref{3.10}) we get:
\begin{align*}
\int\limits_{\{\xi:\, |\tilde{X}_{k,\e}|^2\geqslant \vr\}} \bigg( |\nabla_{\tilde{\xi}}\tilde{X}_{k,\e}|^2 Q_{\max,\vr} &+ \frac{1}{2} |\nabla_{\tilde{\xi}}|\tilde{X}_{k,\e}|^2|^2 \bigg)\,dx + \tilde{c}_2\int\limits_{\p \tilde{\om}_{k,\e}} |\tilde{X}_{k,\e}|^2 Q_{\max,\vr} \,ds
 \\
 &\leqslant C \int\limits_{\p \tilde{\om}_{k,\e}} |\tilde{X}_{k,\e}| Q_{\max,\vr} \,ds \leqslant
 \frac{\tilde{c}_2}{2} \int\limits_{\p \tilde{\om}_{k,\e}} |\tilde{X}_{k,\e}|^2 Q_{\max,\vr} \,ds + C \int\limits_{\p\tilde{\om}_{k,\e}} Q_{\max,\vr}\,ds
\end{align*}
and hence,
\begin{equation*}
\int\limits_{\p \tilde{\om}_{k,\e}} \bigg(C-\frac{\tilde{c}_2}{2}|\tilde{X}_{k,\e}|^2\bigg) Q_{\max,\vr} \,ds \geqslant \frac{1}{2} \int\limits_{\{\xi:\, |\tilde{X}_{k,\e}|^2\geqslant \vr\}}  |\nabla_{\tilde{\xi}}|\tilde{X}_{k,\e}|^2|^2\,dx
\end{equation*}
Since the first integral in the left hand side of the above inequality is non-negative and $|\tilde{X}_{k,\e}|^2\geqslant \vr$ on $\supp Q_{\max,\vr}$, we conclude that
\begin{equation*}
\bigg(C-\frac{\tilde{c}_2}{2}\vr\bigg) \int\limits_{\p\om_{k,\e}\cap\{\xi:\ |\tilde{X}_{k,\e}|^2\geqslant \vr\}}  Q_{\max,\vr} \,ds \geqslant \frac{1}{2} \int\limits_{\{\xi:\, |\tilde{X}_{k,\e}|^2\geqslant \vr\}}  |\nabla_{\tilde{\xi}}|\tilde{X}_{k,\e}|^2|^2\,dx.
\end{equation*}
As $\vr>\frac{2 C}{\tilde{c}_2}$, the above inequality is possible only if $\mes_n \{\xi:\, |\tilde{X}_{k,\e}|^2\geqslant \vr\}=0$. Hence, the function $\tilde{X}_{k,\e}$ is belongs to $L_\infty(\mathds{R}^d\setminus\tilde{\om}_{k,\e})$ and satisfies (\ref{3.86}). The proof is complete.
\end{proof}

\begin{lemma}\label{lm4.4}
The functions $\tilde{X}_{k,\e}$ belong to $W_{2n}^2(\mathds{R}^d\setminus\tilde{\om}_{k,\e})$ and satisfy the estimates
\begin{equation*}
\|\tilde{X}_{k,\e}\|_{W_{2n+1}^2(\mathds{R}^d\setminus\tilde{\om}_{k,\e})} \leqslant C,
\end{equation*}
where $C$ is a constant independent of $\e$ and $k$.
\end{lemma}

\begin{proof}
For $k\in\mathds{M}_D^\e\cup \mathds{M}_{R,1}^\e$ by \cite[Ch. I\!I\!I, Sect. 15, Thm. 15.1]{Ld} we conclude that $\tilde{X}_{k,\e}\in W_2^2(\mathds{R}^d\setminus\tilde{\om}_{k,\e})$. Then we observe that in view of definition (\ref{3.79}) of the Kelvin transform and Assumption~\ref{A1}, the boundaries of the domains $\tilde{\om}_{k,\e}$ have the same regularity as described in this assumption. This allows us to reproduce the proof of the apriori estimate from  \cite[Ch. 15, Thm. 15]{ADN1} controlling at the same time the dependence of the constants on the boundaries; this is done while making a standard unity partition. Then, in view of Lemma~\ref{lm6.5} and the uniform boundedness of $\tilde{\phi}_{k,\e}$ we have:
\begin{equation*}
\|\tilde{X}_{k,\e}\|_{W_{2n+1}^2(\mathds{R}^d\setminus\tilde{\om}_{k,\e})} \leqslant C
\big(\|\tilde{\phi}_{k,\e}\|_{C^2(\p\tilde{\om}_{k,\e})} +\|\tilde{X}_{k,\e}\|_{L_{2n+1}(\mathds{R}^d\setminus\tilde{\om}_{k,\e})}\big)\leqslant C.
\end{equation*}
For $k\in\mathds{M}_{R,2}^\e$ the above apriori estimate also holds true; one just should use the norm $\|\tilde{\phi}_{k,\e}\|_{C^1(\p\tilde{\om}_{k,\e})}$. This is why, to complete the proof, we need to show that $\tilde{X}_{k,\e}$ is an element of $W_{2n}^2(\mathds{R}^d\setminus\tilde{\om}_{k,\e})$. This can be done by using an approximation technique from the proof of Theorem~8.34 in \cite[Ch. 8, Sect. 8.11]{GT}. Namely, the domains $\tilde{\om}_{k,\e}$ are to be approximated by a sequence of domains $\tilde{\om}_{k,\e,m}$ with $C^3$-boundaries; we can simply assume that the boundaries are described by the equations $\tilde{\tau}=\a_m(\tilde{s})$, where $\tilde{\tau}$ is the distance along the normal vector $\tilde{\nu}$ to $\p\tilde{\om}_{k,\e}$, the symbol $\tilde{s}$ denotes local variables on $\p\tilde{\om}_{k,\e}$ and $\a_m$ is a sequence of some functions converging to zero  in $C^2(\p\tilde{\om}_{k,\e})$ as $m\to\infty$. The functions $\tilde{b}_{k,\e}$ then are also extended to the surfaces $\tilde{\tau}=\a_m(\tilde{s})$ just by assuming that they are independent of $\tilde{\tau}$ and on each such surface the function $\tilde{b}_{k,\e}$ is approximating by a $C^2$-function $\tilde{b}_{k,\e,m}$ which, in the sense of the above translation along $\tilde{\tau}$, converges to $\tilde{b}_{k,\e}$ in $C^1$-norm as $m\to\infty$. The functions $\tilde{\phi}_{k,\e}$ are also approximated in the same way by a sequence of $C^2$-functions $\tilde{\phi}_{k,\e,m}$ converging to $\tilde{\phi}_{k,\e}$ in $C^1(\p\tilde{\om}_{k,\e})$. Then we consider problems similar to (\ref{3.85}) for $k\in\mathds{M}_{R,2}^\e$ and these solutions are uniquely solvable in $C^2(\mathds{R}^d\setminus\tilde{\om}_{k,\e})$ due to the standard Schauder estimates. We can then map the domains $\tilde{\om}_{k,\e,m}$ onto $\tilde{\om}_{k,\e}$ and this also transforms the approximating problems; we denote their solutions (after the mapping onto $\tilde{\om}_{k,\e}$) by $\tilde{X}_{k,\e,m}$.  These solutions are the elements of the space $W_{2n}^2(\mathds{R}^d\setminus\tilde{\om}_{k,\e})$ and they also satisfy uniform bounds
\begin{equation*}
\|\tilde{X}_{k,\e,m}\|_{W_{2n+1}^2(\mathds{R}^d\setminus \tilde{\om}_{k,\e})}\leqslant C
\end{equation*}
with some constant $C$ independent of $m$. Then, as in the proof of Theorem~8.34 in \cite[Ch. 8, Sect. 8.11]{GT}, we easily show that the sequence $\tilde{X}_{k,\e,m}$ contains a subsequence weakly converging in $W_{2n}^2(\mathds{R}^d\setminus \tilde{\om}_{k,\e})$ and its weak limit coincides with $\tilde{X}_{k,\e}$. The proof is complete.
\end{proof}

\begin{lemma}\label{lm4.2}
There exists  a fixed positive constant $C_0>0$ independent of $\e$ and $k\in\mathds{M}^\e$ such that for $|\xi|\geqslant C_0$ the functions $X_{k,\e}$ satisfy the representations
\begin{equation}\label{3.29}
\begin{aligned}
&X_{k,\e}(\xi)=1+\frac{K_{k,\e}}{|\rA^{-\frac{1}{2}}(M_k^\e)\xi|^{n-2}} + \sum\limits_{j=1}^{n} K_{k,\e}^{(j)} \frac{\big(\rA^{-\frac{1}{2}}(M_k^\e)\xi\big)_j}{|\rA^{-\frac{1}{2}}(M_k^\e)\xi|^{n}} + \mathsf{X}_{k,\e}(\xi), && \xi\to\infty\quad \text{as}\quad n\geqslant 3,
\\
&X_{k,\e}(\xi)=\ln |\rA^{-\frac{1}{2}}(M_k^\e)\xi|+K_{k,\e} + \sum\limits_{j=1}^{2} K_{k,\e}^{(j)}
 \frac{\big(\rA^{-\frac{1}{2}}(M_k^\e)\xi\big)_j}{|\rA^{-\frac{1}{2}}(M_k^\e)\xi|^{2}} + \mathsf{X}_{k,\e}(\xi), && \xi\to\infty\quad \text{as}\quad n=2,
\end{aligned}
\end{equation}
where $K_{k,\e}^{(j)}$ are some constants, $\big(\rA^{-\frac{1}{2}}(M_k^\e)\xi\big)_j$ is the $j$th component of the vector $\rA^{-\frac{1}{2}}(M_k^\e)\xi$, while  $\mathsf{X}_{k,\e}$ are some infinitely differentiable functions such that
\begin{equation}\label{3.30}
\left|\frac{\p^\vt\mathsf{X}_{k,\e}}{\p\xi^\vt}(\xi)\right|\leqslant C_1 |\xi|^{n-|\vt|},\quad |\xi|\geqslant C_0,
\end{equation}
where $C_1$ is some constant independent of $k$, $\e$,  $\xi$ and $\vt\in\mathds{Z}_+^n$ is a multi-index, $|\vt|\leqslant 3$. The estimates hold:
\begin{align}\label{3.50a}
&|K_{k,\e}|\leqslant C,\qquad |K_{k,\e}^{(j)}|\leqslant C,\quad j=1,\ldots,n,
\\
& |X_{k,\e}(\xi)|\leqslant C,\qquad n\geqslant 3,\qquad |X_{k,\e}(\xi) -  \ln|\rA^{-\frac{1}{2}}(M_k^\e)\xi||\leqslant C, \quad \quad n=2,
\label{3.51}
\\
&|\nabla_\xi X_{k,\e}(\xi)|\leqslant C,  \label{3.90}
\end{align}
where $C$
are some  constants independent of $k$, $\e$,  $\xi$.
\end{lemma}

\begin{proof}
By Lemma~\ref{lm4.4} and the definition of the Kelvin transform in (\ref{3.79}), the functions $X_{k,\e}$ are the elements of the space $W_{2n}^2(B_{R_3}(0)\setminus\om_{k,\e})$ and are bounded in this space uniformly in $k$ and $\e$. Using the regularity of the boundaries $\p\om_{k,\e}$ postulated in Assumption~\ref{A1}, we continue the function $X_{k,\e}$ into $\om_{k,\e}$ as follows:
\begin{equation*}
X_{k,\e}(\tau,s)=X_{k,\e}(-\tau,s)\chi_{1}(\tau),
\end{equation*}
where $\chi_{1}=\chi_{1}(\tau)$ is an infinitely differentiable function vanishing as $|\tau|>\frac{2\tau_0}{3}$ and equalling to one as $|\tau|<\frac{\tau_0}{3}$. In the same way we continue each derivative $\frac{\p X_{k,\e}}{\p x_j}$, $j=1,\ldots,n$.  It is clear that after such continuation the obtained functions   are elements of $W_{2n+1}^1(B_{R_3}(0))$ bounded in this space uniformly in $\e$ and $k$. Applying then Sobolev theorem \cite[Ch. I\!I, Sect. 2, Thm. 2.2]{Ld} we see that $X_{k,\e}, \frac{\p X_{k,\e}}{\p x_j}\in C^\frac{1}{2}(\overline{B_{R_2}(0)\setminus\om_{k,\e}})$ and   the estimate holds:
\begin{equation}\label{3.21}
\|X_{k,\e}\|_{C^\frac{3}{2}(\overline{B_{R_2}(0)\setminus\om_{k,\e}})}\leqslant C,
\end{equation}
where $C$ is a constant independent of $k$ and $\e$.

We again use the Kelvin transform introduced in the proof of Lemma~\ref{lm6.1} and even for $k\in\mathds{M}_{R,2}^\e$ we treat $\tilde{X}_{k,\e}$ as solutions to the Dirichlet problem with appropriate $\tilde{\phi}_{k,\e}$. Estimate (\ref{3.21}) allows us to say that these functions $\tilde{\phi}_{k,\e}$ are bounded uniformly in $k$, $\e$ and $\tilde{\xi}$. The well-known estimates for the derivatives of the harmonic function, see, for instance, \cite[Ch. I\!V, Sect. 3.2, Lm. 3]{Mi}, then yield
\begin{equation*}
\left|\frac{\p^\vt \tilde{X}_{k,\e}}{\p\tilde{\xi}^\vt}(\tilde{\xi})\right|
\leqslant  C, \qquad \tilde{\xi}\in B_{(3R_2)^{-1}}(0),\qquad \vt\in\mathds{Z}_+^n,\qquad |\kappa|\leqslant 3,
\end{equation*}
where $C$ are some   constants independent of $k$, $\e$ and $\tilde{\xi}$.
 Using Lemma~\ref{lm6.5}, writing the Taylor series for $\tilde{X}_{k,\e}$ at zero and returning back to the function $X_{k,\e}$ by formula (\ref{3.81}), we   get (\ref{3.29}), (\ref{3.30}), (\ref{3.50a}), (\ref{3.51}), (\ref{3.90}). The proof is complete.
\end{proof}

The next lemma states the unique solvability of problem (\ref{2.34}), (\ref{2.35}), (\ref{2.27}), (\ref{2.28}) and provides certain estimates for its solution.

\begin{lemma}\label{lm3.2}
Problems (\ref{2.34}), (\ref{2.35}), (\ref{2.27}), (\ref{2.28}) are uniquely solvable in $W_2^2(E_{k,\e}\setminus\om_{k,\e})\cap C^1(\overline{E_{k,\e}\setminus\om_{k,\e}})$ and satisfy  the estimates
\begin{equation}
|X_{k,\e}-Z_{k,\e}|\leqslant C\eta^{n-1},
\qquad
|\nabla_\xi X_{k,\e}-\nabla_\xi Z_{k,\e}|\leqslant C\eta^{n-1}  \quad\text{in}\quad \overline{E_{k,\e}\setminus\om_{k,\e}},
\label{3.82}
\end{equation}
where $C$ is some constant independent of $\e$, $\eta$, $k$ and $\xi$.
\end{lemma}

\begin{proof}
The unique solvability in $W_2^2(E_{k,\e}\setminus\om_{k,\e})$ is easily checked in the same way how a similar fact was established in the proof of Lemma~\ref{lm6.1}; here we even do not need to make the Kelvin transform since the domains $E_{k,\e}\setminus\om_{k,\e}$ are bounded. By the standard smoothness improving theorems we also see that $Z_{k,\e}\in C^\infty(\overline{E_{k,\e}}\setminus\overline{\om_{k,\e}})$.

We let
\begin{equation*}
\psi_{k,\e}^+(\xi):=X_{k,\e}(\xi)-\sum\limits_{j=1}^{n} \frac{K_{k,\e}^{(j)}\eta^n}{R_4^n} \big(\rA^{-\frac{1}{2}}(M_k^\e)\xi\big)_j-\left\{
\begin{aligned}
&1+K_{k,\e} |\rA^{-\frac{1}{2}}(M_k^\e)\xi|^{-n+2} &&\text{as}\quad n\geqslant 3,
\\
&\ln |\rA^{-\frac{1}{2}}(M_k^\e)\xi| + K_{k,\e} &&\text{as}\quad n=2,
\end{aligned}\right.
\end{equation*}
and by (\ref{3.29}), (\ref{3.30}) we see that
\begin{equation}\label{3.24}
\|\psi_{k,\e}^+\|_{C^3(\p E_{k,\e})}\leqslant C\eta^n;
\end{equation}
hereinafter in the proof by $C$ we denote inessential constants independent $\e$, $k$ and $\eta$.  By $T_{k,\e}^0(\z)$ we denote the solution to the problem
\begin{equation}\label{3.20}
\D_\z T_{k,\e}^0=0\quad\text{in}\quad B_{R_4\eta^{-1}}(0),\qquad T_{k,\e}^0=\psi_{k,\e}^+\big(\rA^{\frac{1}{2}}(M_k^\e)\z\big) \quad\text{on}\quad\p B_{R_4\eta^{-1}}(0).
\end{equation}
This problem is uniquely solvable. We reproduce the proof of the Schauder estimate \cite[Ch. I\!I\!I, Sect. 1, 2]{Ld} for problem (\ref{3.20}) covering $B_{R_4\eta^{-1}}(0)$ by balls of a fixed radius and in view of (\ref{3.24}) we conclude that $T_{k,\e}^0\in C^2(\overline{B_{R_4\eta^{-1}}(0)})$  and
\begin{equation}\label{3.19}
\|T_{k,\e}^0\|_{C^2(\overline{B_{R_4\eta^{-1}}(0)})}\leqslant C \eta^n.
\end{equation}
The functions
\begin{equation}\label{3.38}
T_{k,\e}(\xi):=X_{k,\e}(\xi)-Z_{k,\e}(\xi) - \sum\limits_{j=1}^{n} \frac{K_{k,\e}^{(j)}\eta^n}{R_4^n} \big(\rA^{-\frac{1}{2}}(M_k^\e)\xi\big)_j - T_{k,\e}^{0}\big(\rA^{-\frac{1}{2}}(M_k^\e)\xi\big)
\end{equation}
solve the boundary value problems
\begin{equation}
\begin{gathered}
\dvr_\xi \rA(M_k^\e) \nabla_\xi T_{k,\e}=0\quad\text{in}\quad E_{k,\e}\setminus\om_{k,\e},
\qquad T_{k,\e}=0\quad \text{on}\quad \p E_{k,\e}
\\
T_{k,\e}=\psi_{k,\e}^-,\quad k\in\mathds{M}_D^\e\cup\mathds{M}_{R,1}^\e,\qquad
 \nu_\xi\cdot\rA(M_k^\e)\nabla_\xi T_{k,\e}+b_k T_{k,\e}=\psi_{k,\e}^-,\quad k\in \mathds{M}_{R,2}^\e\quad\text{on}\quad\p\om_{k,\e},
\end{gathered}\label{3.34}
\end{equation}
where
\begin{align*}
&\psi_{k,\e}^-(\xi):=-\sum\limits_{j=1}^{n} \frac{K_{k,\e}^{(j)}\eta^n}{R_4^n}\big(\rA^{-\frac{1}{2}}(M_k^\e)\xi\big)_j-T_{k,\e}^{0},\quad k\in \mathds{M}_D^\e\cup\mathds{M}_{R,1}^\e,
\\
&\psi_{k,\e}^-(\xi):=-\Big(\nu_\xi\cdot\rA(M_k^\e)\nabla_\xi T_{k,\e}+b_k T_{k,\e}\Big)\Bigg(\sum\limits_{j=1}^{n} \frac{K_{k,\e}^{(j)}\eta^n}{R_4^n}\big(\rA^{-\frac{1}{2}}(M_k^\e)\xi\big)_j+T_{k,\e}^0\Bigg),
\quad k\in \mathds{M}_{R,2}^\e.
\end{align*}
It follows from (\ref{3.19}) that
\begin{equation*}
\|\psi_{k,\e}^-\|_{C^2(\p \om_{k,\e})}\leqslant C\eta^n,\quad k\in\mathds{M}_D^\e\cup\mathds{M}_{R,1}^\e,\qquad \|\psi_{k,\e}^-\|_{C^1(\p \om_{k,\e})}\leqslant C\eta^n, \quad  k\in \mathds{M}_{R,2}^\e.
\end{equation*}
Writing then integral identities associated with   problem (\ref{3.34}), we easily obtain:
\begin{equation*}
\|\nabla_\xi T_{k,\e}\|_{L_2(E_{k,\e}\setminus \om_{k,\e})}+\|T_{k,\e}\|_{L_2(\p\om_{k,\e})} \leqslant C \eta^n,
\end{equation*}
where the second term in the left hand side obviously vanishes for $k\in\mathds{M}_D^\e\cup \mathds{M}_{R,1}^\e$. By the standard smoothness improving theorems we then get:
\begin{equation}\label{3.26}
\|T_{k,\e}\|_{C^3(\overline{B_{R_3}(0)\setminus B_{(R_2+R_3)/2}(0)})}\leqslant C \eta^n.
\end{equation}

Let $\chi_{2}=\chi_{2}(\xi)$ be an infinitely differentiable cut-off function equalling to one as $|\xi|<(R_2+R_3)/2$ and vanishing as $|\xi|>R_3$. Then the functions $T_{k,\e}\chi_{2}$ solve the boundary value problems
\begin{gather*}
\dvr_\xi \rA(M_k^\e) \nabla_\xi T_{k,\e}\chi_{2}=2\nabla_\xi \chi_{2} \cdot \rA(M_k^\e) \nabla_\xi T_{k,\e}+ T_{k,\e}\dvr_\xi \rA(M_k^\e) \nabla_\xi\chi_{2} \quad \text{in}\quad B_{R_3}(0)\setminus \om_{k,\e},
\\
T_{k,\e}\chi_{2}=0\quad\text{on}\quad\p B_{R_3}(0)
\end{gather*}
with boundary conditions (\ref{2.27}), (\ref{2.28}).  These problems can be studied following the lines of the proofs of Lemmata~\ref{lm6.5},~\ref{lm4.4} by employing also estimate (\ref{3.26}) and the inequality \cite[Lm. 3.1]{BK}:
\begin{equation*}
\|u\|_{L_2(B_{R_3}(0)\setminus\om_{k,\e})}\leqslant C \|\nabla_\xi u\|_{L_2(B_{R_3}(0)\setminus\om_{k,\e})}\quad \text{for all}\quad u\in \Ho_2^1(B_{R_3}(0)\setminus\om_{k,\e},\p B_{R_3}(0))
\end{equation*}
with a constant $C$ independent of $k$, $\e$ and $u$. As a result we obtain:
\begin{equation}\label{3.28}
\|T_{k,\e}\|_{C^1(\overline{B_{R_3}(0)\setminus\om_{k,\e}})}\leqslant C\eta^n.
\end{equation}
We also conclude that $T_{k,\e}\in C^1(\overline{E_{k,\e}\setminus\om_{k,\e}})$. Since the function $T_{k,\e}(\rA^\frac{1}{2}(M_k^\e)\xi)$ is obviously harmonic,  by the classical maximum principle for the harmonic functions and the first estimate in (\ref{3.24}) we immediately obtain
\begin{equation}\label{3.36}
\|T_{k,\e}\|_{C(\overline{E_{k,\e}\setminus\om_{k,\e}})}\leqslant C \eta^n.
\end{equation}

For $k\in\mathds{M}_{R,2}^\e$ we also have an appropriate maximum principle. Namely,
the real and imaginary parts of   $T_{k,\e}(\rA^\frac{1}{2}(M_k^\e)\xi)$ are harmonic functions and by the mean value theorem
\begin{equation*}
T_{k,\e}(\rA^\frac{1}{2}(M_k^\e)\xi)=\frac{1}{\mes_n B_\d(\xi)} \int\limits_{B_\d(\xi)} T_{k,\e}(\rA^\frac{1}{2}(M_k^\e)y)\,dy
\end{equation*}
for each $\xi\in E_{k,\e}\setminus\om_{k,\e}$ and each ball $B_\d(\xi)$ such that
$\big\{\rA^\frac{1}{2}(M_k^\e)y:\ y\in B_\d(\xi)\big\}\subset E_{k,\e}\setminus\om_{k,\e}$.
This identity implies
\begin{equation*}
|T_{k,\e}(\rA^\frac{1}{2}(M_k^\e)\xi)|\leqslant \max\limits_{y\in B_\d(\xi)} |T_{k,\e}(\rA^\frac{1}{2}(M_k^\e) y)|.
\end{equation*}
If $\xi\in E_{k,\e}\setminus\overline{\om_{k,\e}}$ is a point of the global maximum of $|T_{k,\e}|$, then the above inequality implies that $|T_{k,\e}|$ is constant in $B_\d(\xi)$. Hence, the function $|T_{k,\e}|$  attains its global maximum on the boundary $\p\om_{k,\e}$ or on $\p E_{k,\e}$. It follows from the boundary condition for $T_{k,\e}$ that
\begin{equation*}
\frac{\p|T_{k,\e}|^2}{\p \nu} +\RE b_{k,\e} |T_{k,\e}|^2=0.
\end{equation*}
If  a point of the global maximum of $T_{k,\e}$ is located on $\p\om_{k,\e}$, then $\frac{\p|T_{k,\e}|^2}{\p \nu}\geqslant 0$ and the above identity due to the positivity of $\RE b_k$, see (\ref{2.26}), implies that $T_{k,\e}=0$. Hence, the function $|T_{k,\e}|$ attains its maximum on $\p E_{k,\e}$ and this gives estimate (\ref{3.36}) for $k\in \mathds{M}_{R,2}^\e$.

We consider the function $T_{k,\e}$ as a solution of the equation from (\ref{3.34}) but on $E_{k,\e}\setminus B_{R_3}(0)$ subject to the boundary condition on $\p E_{k,\e}$ from (\ref{3.34}). Taking into consideration then (\ref{3.26}) and reproducing again the proof of the Schauder estimate with covering by balls of fixed radius, we obtain
\begin{equation*}
\|T_{k,\e}\|_{C^2(\overline{E_{k,\e}\setminus B_{R_3}(0)})} \leqslant C \eta^n.
\end{equation*}
This estimate and (\ref{3.28}) yield
\begin{equation*}
\|T_{k,\e}\|_{C^1(\overline{E_{k,\e}\setminus \om_{k,\e}})} \leqslant C \eta^n.
\end{equation*}
Returning back then to the function $X_{k,\e}-Z_{k,\e}$ by formula (\ref{3.38}) and using estimates (\ref{3.19}), (\ref{3.50a}), we arrive at (\ref{3.82}). The proof is complete.
\end{proof}

The above lemmata implies several properties of the function $\Xi_\e$.
We first  observe that
estimates (\ref{3.29}), (\ref{3.30}),  (\ref{3.82}) imply
\begin{equation}\label{3.66}
|\Xi_\e(x)-1|\leqslant C\quad\text{in}\quad B_{\e R_3}(M_k^\e)\setminus\om_k^\e.
\end{equation}
Employing (\ref{3.29}), (\ref{3.30}), by straightforward calculations we find:
\begin{align*}
&\rA(M_k^\e)  \nabla_\xi X_{k,\e} \cdot \nu=
(2-n)K_{k,\e}
\frac{\eta^{n-2}}{R_4^{n-2}} |\rA^{-1}(M_k^\e)\xi|^{-1}+O(\eta^n)  && \hspace{-1.7 true cm}\text{as}\qquad n\geqslant 3,
\\
&\rA(M_k^\e)  \nabla_\xi X_{k,\e} \cdot \nu=   
|\rA^{-1}(M_k^\e)\xi|^{-1}+O(\eta^2) &&\hspace{-1.7 true cm} \text{as}\qquad n=2,
\end{align*}
on $\p E_{k,\e}$, where $\nu$ is the outward normal to $\p E_{k,\e}$ and the $O$-terms are uniform in $\e$, $k$ and $\eta$. Using  then (\ref{2.14}), (\ref{3.82}) and the definition of $\Xi_\e$  we obtain:
\begin{equation}\label{4.54}
\begin{aligned}
&\rA(M_k^\e)  \nabla \Xi_\e\cdot \nu=(2-n) K_{k,\e}\frac{\eta^{n-2}}{R_4^{n-2}} |\rA^{-1}(M_k^\e)(x-M_k^\e)|^{-1}+O\big(\eta^{n-1}\e^{-1}\big)  && \quad \text{as}\qquad n\geqslant 3,
\\
&\rA(M_k^\e)  \nabla\Xi_\e\cdot \nu=  |\ln\eta|^{-1}|\rA^{-1}(M_k^\e)(x-M_k^\e)|^{-1} + O\big(\e^{-1}\ln^{-2}\eta
\big)
&& \quad\text{as}\qquad n=2,
\end{aligned}
\end{equation}
on $\p E_k^\e$, where the $O$-terms are uniform in $\e$, $k$ and $\eta$.

\begin{lemma}\label{lm4.3}
The estimates hold
\begin{equation}\label{3.63}
 \|(\Xi_\e-1)u\|_{L_2(B_{\e R_3}(M_k^\e)\setminus\om_k^\e)} \leqslant  C \e \|u\|_{W_2^1(B_{\e R_3}(M_k^\e)\setminus\om_k^\e)}
\end{equation}
for all $u\in W_2^1(B_{\e R_3}(M_k^\e)\setminus\om_k^\e)$ and
\begin{align}
 &\|u \nabla \Xi_\e\|_{L_2(B_{\e R_3}(M_k^\e)\setminus\om_k^\e)} \leqslant  C \big(\eta^{\frac{n}{2}-1}\e^{-1}\vk^{-\frac{1}{2}} +\e^{\frac{1}{2}}\eta^{\frac{1}{2}}
 +\e\big) \|u\|_{W_2^2(B_{\e R_3}(M_k^\e)\setminus\om_k^\e)}, \label{3.64}
 \\
 &\|(\Xi_\e-1)u\|_{L_2(B_{\e R_3}(M_k^\e)\setminus\om_k^\e)} \leqslant  C (\e^2+\e \eta) \|u\|_{W_2^2(B_{\e R_3}(M_k^\e)\setminus\om_k^\e)}
 \label{3.63a}
\end{align}
for all $u\in W_2^2(B_{\e R_3}(M_k^\e)\setminus\om_k^\e)$,
where $C$ are some constants  independent of the parameters $\e$, $k$ and the function $u$.
\end{lemma}

\begin{proof}
We fix $k\in\mathds{M}^\e$ and for a given $u\in   W_2^1(B_{\e R_3}(M_k^\e)\setminus\om_k^\e)$ we denote
\begin{equation}\label{4.56}
\la u \ra:=\frac{1}{\mes_n B_{\e R_3}(M_k^\e)\setminus\om_k^\e} \int\limits_{B_{\e R_3}(M_k)\setminus\om_k^\e} u(x)\,dx,\qquad u^\bot:=u-\la u \ra.
\end{equation}
Then, in view of the fact that $\Xi_\e$ is identically one in $B_{\e R_3}(M_k^\e)\setminus E_k^\e$,
\begin{equation}\label{3.65}
\|(\Xi_\e-1)u\|_{L_2(B_{\e R_3}(M_k^\e)\setminus\om_k^\e)}^2\leqslant
2|\la u\ra|^2 \|(\Xi_\e-1)\|_{L_2(E_k^\e\setminus\om_k^\e)}^2
+2\|(\Xi_\e-1)u^\bot\|_{L_2(B_{\e R_3}(M_k^\e)\setminus\om_k^\e)}^2.
\end{equation}
Since the function $u^\bot$ obeys condition (\ref{3.1}), by Lemma~\ref{lm3.3} it satisfies estimate (\ref{3.2a}) and by  (\ref{3.66})  we immediately get
\begin{equation}\label{3.67}
\|(\Xi_\e-1)u^\bot\|_{L_2(B_{\e R_3}(M_k^\e)\setminus\om_k^\e)}^2 \leqslant C\e^2\|u\|_{W_2^1(B_{\e R_3}(M_k^\e)\setminus\om_k^\e)}^2;
\end{equation}
hereinafter in the proof by $C$ we denote various inessential constants independent of $\e$, $k$ and $u$.
Passing then to the variables $\xi=(x-M_k^\e)\e^{-1}\eta^{-1}$, by
(\ref{3.29}), (\ref{3.30}), (\ref{3.50a}), (\ref{3.82}), (\ref{3.66}) for $n\geqslant 3$ we find:
\begin{align*}
\|\Xi_\e-1\|_{L_2(E_k^\e\setminus\om_k^\e)}^2 \leqslant & C\e^n\eta^n \|Z_{k,\e}-1-K_{k,\e}R_4^{-n+2}\eta^{n-2}\|_{L_2(E_{k,\e}\setminus\om_{k,\e})}^2
 \\
 \leqslant & C \e^n \eta^{2n-2} +C\e^n\eta^n \|X_{k,\e}-1-K_{k,\e}R_4^{-n+2}\eta^{n-2}\|_{L_2(E_{k,\e}\setminus\om_{k,\e})}^2
\\
\leqslant & C \e^n\eta^n + C\e^n\eta^n \int\limits_{C_0}^{R_4\eta^{-1}} \big((r^{-n+2}-R_4^{-n+2}\eta^{n-2})^2+r^{-2n+2}\big)r^{n-1}\,dr
\\
\leqslant & C\e^n\eta^n \left\{
\begin{aligned}
&1+\eta^{n-4}, && n\ne 3,
\\
&|\ln\eta|, && n=4.
\end{aligned}\right.
\end{align*}
As $n=2$, we estimate along the same lines:
\begin{align*}
\|\Xi_\e-1\|_{L_2(E_k^\e\setminus\om_k^\e)}^2 \leqslant & C\e^2\eta^2 \ln^{-2}\eta\|Z_{k,\e}-\ln R_4\eta^{-1}1-K_{k,\e}\|_{L_2(E_{k,\e}\setminus\om_{k,\e})}^2
\\
\leqslant & C \e^2\eta^2\ln^{-2}\eta + C\e^2\eta^2 \ln^{-2}\eta\|X_{k,\e}-\ln R_4\eta^{-1}1-K_{k,\e}\|_{L_2(E_{k,\e}\setminus\om_{k,\e})}^2
\\
\leqslant & C \e^2\eta^2\ln^{-2}\eta + C\e^2\eta^2 \ln^{-2}\eta  \int\limits_{C_0}^{R_4\eta^{-1}}
\left(\ln^2 \frac{r}{R_4\eta^{-1}} + r^{-2}\right)r\,dr \leqslant C \e^2\eta^2.
\end{align*}
Hence, in view of convergence  (\ref{2.14}),
\begin{equation}\label{3.69}
\|\Xi_\e-1\|_{L_2(E_k^\e\setminus\om_k^\e)}^2 \leqslant
C\e^{n+2}.
\end{equation}
We also see easily that
\begin{equation}\label{3.69a}
|\la u\ra|^2\leqslant C\e^{-n}\|u\|_{L_2(B_{\e R_3}(M_k^\e)\setminus\om_k^\e)}^2.
\end{equation}
Employing this estimate and (\ref{3.65}), (\ref{3.67}), (\ref{3.69}), (\ref{2.14}), we obtain (\ref{3.63}). In the same way we also prove easily the estimate
\begin{equation}\label{3.76}
\||\Xi_\e-1|^\frac{1}{2}u\|_{L_2(B_{\e R_3}(M_k^\e)\setminus\om_k^\e)}
\leqslant C\e\|u\|_{W_2^1(B_{\e R_3}(M_k^\e)\setminus\om_k^\e)}.
\end{equation}

We proceed to proving (\ref{3.64}) and (\ref{3.63a}). We begin with simple relations:
\begin{equation}\label{3.70}
\|u \nabla \Xi_\e\|_{L_2(B_{\e R_3}(M_k^\e)\setminus\om_k^\e)}^2= \|u \nabla \Xi_\e\|_{L_2(E_k^\e\setminus\om_k^\e)}^2
\leqslant  C\big(\rA(M_k^\e) u \nabla \Xi_\e, u \nabla \Xi_\e\big)_{L_2(B_{\e R_3}(M_k^\e)\setminus\om_k^\e)}.
\end{equation}
Then we integrate by parts employing the definition of the function $\Xi_\e$:
\begin{equation}\label{3.13}
\begin{aligned}
\big(\rA(M_k^\e) u \nabla  \Xi_\e, u \nabla \Xi_\e\big)_{L_2(B_{\e R_3}(M_k^\e)\setminus\om_k^\e)} = &\RE\int\limits_{\p E_k^\e} |u|^2 \rA(M_k^\e) \nabla \Xi_\e \cdot\nu\,ds  + \RE\int\limits_{\p\om_k^\e} \overline{\Xi_\e}
|u|^2 \rA(M_k^\e) \nabla \Xi_\e \cdot\nu\,ds
\\
&- \RE\int\limits_{E_k^\e\setminus\om_k^\e} \overline{\Xi_\e} \dvr \rA(M_k^\e) |u|^2 \nabla \Xi_\e \,dx
\\
= & \RE\int\limits_{\p E_k^\e} |u|^2 \rA(M_k^\e) \nabla \Xi_\e \cdot\nu\,ds + \RE\int\limits_{\p\om_k^\e} \overline{\Xi_\e}
|u|^2 \rA(M_k^\e) \nabla \Xi_\e \cdot\nu\,ds
\\
 &- \frac{1}{2} \int\limits_{E_k^\e\setminus\om_k^\e}\rA(M_k^\e)\nabla (|\Xi_\e|^2-1) \cdot\nabla |u|^2\,dx
\\
=& \RE\int\limits_{\p E_k^\e} |u|^2 \rA(M_k^\e) \nabla \Xi_\e \cdot\nu\,ds + \RE\int\limits_{\p\om_k^\e} \overline{\Xi_\e}
|u|^2 \rA(M_k^\e) \nabla \Xi_\e \cdot\nu\,ds
\\
&-
\frac{1}{2} \int\limits_{\p \om_k^\e} (|\Xi_\e|^2-1) \rA(M_k^\e) \nabla |u|^2\cdot\nu\,ds
\\
&+ \frac{1}{2}
\int\limits_{E_k^\e\setminus\om_k^\e} (|\Xi_\e|^2-1) \dvr \rA(M_k^\e) \nabla |u|^2 \,dx,
\end{aligned}
\end{equation}
where $\nu$ denotes the unit normal to $\p E_k^\e$ directed outside $E_k^\e$ and also  the unit normal to $\p\om_k^\e$ directed inside  $\om_k^\e$. In view of the boundary conditions for $X_{k,\e}$ in (\ref{2.27}), (\ref{2.28}) and the inequality for $b_k$ in (\ref{2.26}) we also have:
\begin{align*}
& \RE\int\limits_{\p\om_k^\e} \overline{\Xi_\e}
|u|^2 \rA(M_k^\e) \nabla \Xi_\e \cdot\nu^\e\,ds=0, && k\in\mathds{M}_D^\e\cup\mathds{M}_{R,1}^\e,
\\
& \RE\int\limits_{\p\om_k^\e} \overline{\Xi_\e}
|u|^2 \rA(M_k^\e) \nabla \Xi_\e \cdot\nu^\e\,ds=-\e^{-1}\eta^{-1}\RE b_k(\e)\|u\|_{L_2(\p\om_k^\e)}^2\leqslant 0, && k\in\mathds{M}_{R,2}^\e.
\end{align*}
Hence, in view of relations (\ref{3.66}),
(\ref{3.76}), (\ref{3.70}), (\ref{3.22}),  (\ref{2.14}), (\ref{4.54})
  and the definition of the function $\Xi_\e$, by (\ref{3.13}) we get:
\begin{equation}\label{3.77}
\begin{aligned}
\|u \nabla \Xi_\e\|_{L_2(B_{\e R_3}(M_k^\e)\setminus\om_k^\e)}^2 \leqslant & C   \big(\eta^{n-2}\e^{-1}\vk^{-1}\|u\|_{L_2(\p E_k^\e)}^2
+\|\nabla u\|_{L_2(\p \om_k^\e)}\|u\|_{L_2(\p \om_k^\e)}\big)
\\
&+ C\|u\|_{W_2^2(E_k^\e\setminus\om_k^\e)}\| (\Xi_\e-1) u\|_{L_2(E_k^\e\setminus\om_k^\e)}
+ C\big\||\Xi_\e-1|^{\frac{1}{2}}\nabla u \big\|_{L_2(E_k^\e\setminus\om_k^\e)}^2
\\
\leqslant &
C(\e^2+\e\eta\vk)\|u\|_{W_2^2(B_{\e R_3}(M_k^\e)\setminus\om_k^\e)}^2
+C \eta^{n-2}\e^{-1}\vk^{-1}\|u\|_{L_2(\p E_k^\e)}^2
\\
&+
C\|u\|_{W_2^2(E_k^\e\setminus\om_k^\e)}\| (\Xi_\e-1) u\|_{L_2(E_k^\e\setminus\om_k^\e)}.
\end{aligned}
\end{equation}
This inequality gives an opportunity to improve (\ref{3.67})
for $u\in W_2^2(B_{\e R_3}(M_k^\e)\setminus\om_k^\e)$. Namely, we integrate by parts and estimate then using (\ref{3.77}) with $u=u^\bot$:
\begin{align*}
\|(\Xi_\e-1)u^\bot \|_{L_2(B_{\e R_3}(M_k^\e)\setminus\om_k^\e)}^2 =& \frac{1}{n}\int\limits_{L_2(B_{\e R_3}(M_k^\e)\setminus\om_k^\e)} (\Xi_\e-1)^2|u^\bot|^2\dvr x\,dx
\\
= & \frac{1}{n}\int\limits_{\p\om_k^\e} |u^\bot|^2 x\cdot\nu\,ds
 -\frac{2}{n}\int\limits_{B_{\e R_3}(M_k^\e)} |u^\bot|(\Xi_\e-1) x\cdot\nabla |u^\bot|(\Xi_\e-1)
\,dx
\\
\leqslant &C\e\eta\|u^\bot\|_{L_2(\p\om_k^\e)}^2+C\e  \|(\Xi_\e-1)u^\bot \|_{L_2(B_{\e R_3}(M_k^\e)\setminus\om_k^\e)}
\|u^\bot \nabla \Xi_\e\|_{L_2(B_{\e R_3}(M_k^\e)\setminus\om_k^\e)}
  \\
  &+C\e \|(\Xi_\e-1)u^\bot \|_{L_2(B_{\e R_3}(M_k^\e)\setminus\om_k^\e)}
 \|(\Xi_\e-1)\nabla u\|_{L_2(B_{\e R_3}(M_k^\e)\setminus\om_k^\e)}
\\
\leqslant & C\e\eta \|u^\bot\|_{L_2(\p\om_k^\e)}^2
 +\frac{1}{4} \|(\Xi_\e-1)u^\bot \|_{L_2(B_{\e R_3}(M_k^\e)\setminus\om_k^\e)}^2
 \\
 &+C\e^2 \left( \|u^\bot \nabla \Xi_\e\|_{L_2(B_{\e R_3}(M_k^\e)\setminus\om_k^\e)}^2 + \|(\Xi_\e-1)\nabla u\|_{L_2(B_{\e R_3}(M_k^\e)\setminus\om_k^\e)}^2
 \right)
 \\
 \leqslant & C\e\eta \|u^\bot\|_{L_2(\p\om_k^\e)}^2
 +\frac{1}{2} \|(\Xi_\e-1)u^\bot \|_{L_2(B_{\e R_3}(M_k^\e)\setminus\om_k^\e)}^2
 \\
 &+C\e^2 \|(\Xi_\e-1)\nabla u\|_{L_2(B_{\e R_3}(M_k^\e)\setminus\om_k^\e)}^2
 + C(\e^4+\e^3\eta\vk) \|u\|_{W_2^2(E_k^\e\setminus\om_k^\e)}^2
 \\
 &+ C \eta^{n-2}\e \vk^{-1} \|u^\bot\|_{L_2(\p E_k^\e)}^2.
\end{align*}
It also follows from (\ref{3.6}) with $\eta=1$, $\p\om_k^\e=\p E_k^\e$, $u=u^\bot$ and (\ref{3.2a}) that
\begin{equation*}
 \|u^\bot\|_{L_2(\p E_k^\e)}^2\leqslant C \e \|\nabla u\|_{L_2(B_{\e R_3}(M_k^\e)\setminus\om_{k,\e})}^2.
\end{equation*}
Therefore,
due to (\ref{3.22}), (\ref{3.63}), (\ref{3.23}),
\begin{align*}
\|(\Xi_\e-1)u^\bot \|_{L_2(B_{\e R_3}(M_k^\e)\setminus\om_k^\e)}^2 \leqslant & C(\e^2\eta^2\vk+\e^4+\e^3\eta\vk) \|u\|_{W_2^2(B_{\e R_3}(M_k^\e)\setminus\om_k^\e)}^2
\\
\leqslant & C(\e^4+\e^2\eta^2)
\|u\|_{W_2^2(B_{\e R_3}(M_k^\e)\setminus\om_k^\e)}^2.
\end{align*}
This estimate and (\ref{3.65}), (\ref{3.69}), (\ref{3.69a}) yield (\ref{3.63a}). Substituting then   (\ref{3.63a}) into the right hand side of (\ref{3.77}) and using (\ref{3.6}) with $\eta=1$, $\p\om_k^\e=\p E_k^\e$,  we arrive at (\ref{3.64}). The proof is complete.
\end{proof}

For each $r>0$ and $x\in\Om$ we denote
\begin{equation*}
\sE_r(x):=\big\{\xi\in\mathds{R}^d:\ |\rA^{-\frac{1}{2}}(x)\xi|<r\big\}.
\end{equation*}

\begin{lemma}\label{lm5.3}
For all $x\in \Om$ the identity holds:
\begin{equation*}
\int\limits_{\p \sE_1(x)}\frac{ds}{|\rA^{-\frac{1}{2}}(x)\xi|}=\Ups(x).
\end{equation*}
\end{lemma}

\begin{proof}
We first consider the case $n\geqslant 3$. For a given $x\in\Om$ we choose  a sufficiently large $r>0$ and denoting by $\nu$ the unit outward normal to $\p \sE_r(x)$, we integrate once by parts as follows:
\begin{align*}
0=& (2-n)^{-1}\int\limits_{\sE_r(x)\setminus\sE_1(x)} |\rA^{-\frac{1}{2}}(x)\xi|^{-n+2} \dvr_\xi \rA(x)\nabla_\xi |\rA^{-\frac{1}{2}}(x)\xi|^{-n+2}\,d\xi
\\
=&(2-n)^{-1}
\int\limits_{\p\sE_r(x)} |\rA^{-\frac{1}{2}}(x)\xi|^{-n+2} \rA(x)\nabla_\xi |\rA^{-\frac{1}{2}}(x)\xi|^{-n+2}\cdot \nu\,d\xi
\\
&- (2-n)^{-1}\int\limits_{\p\sE_1(x)} |\rA^{-\frac{1}{2}}(x)\xi|^{-n+2} \rA(x)\nabla_\xi |\rA^{-\frac{1}{2}}(x)\xi|^{-n+2}\cdot \nu\,d\xi
\\
&-(2-n)^{-1}\int\limits_{\sE_r(x)\setminus\sE_1(x)} \rA(x)\nabla_\xi|\rA^{-\frac{1}{2}}(x)\xi|^{-n+2} \cdot\nabla_\xi |\rA^{-\frac{1}{2}}(x)\xi|^{-n+2}\,d\xi
\\
=&
\int\limits_{\p\sE_r(x)}\frac{ds}{r^{2n-4}|\rA^{-1}(x)\xi|}
-\int\limits_{\sE_r(x)\setminus\sE_1(x)}\frac{ds}{|\rA^{-1}(x)\xi|}- (2-n)\int\limits_{\sE_r(x)\setminus \sE_1(x)} \frac{d\xi}{|\rA^{-\frac{1}{2}}(x)\xi|^{2n-2}}.
\end{align*}
Making the change of the variables $y=\rA^{-\frac{1}{2}}(x)\xi$ in the integral over $\sE_r(x)\setminus\sE_1(x)$ and passing then to the limit as $r\to+\infty$, we obtain:
\begin{equation*}
\int\limits_{\sE_r(x)\setminus\sE_1(x)}\frac{ds}{|\rA^{-1}(x)\xi|}=(n-2)\sqrt{\det \rA(x)}\int\limits_{\mathds{R}^n\setminus B_1(0)} \frac{dy}{|y|^{2n-2}}=\Ups(x)
\end{equation*}
and this proves the needed formula for $n\geqslant 3$.  For $n=2$ we integrate in a similar way:
\begin{align*}
0=& \int\limits_{\sE_r(x)\setminus \sE_1(x)} \ln |\rA^{-\frac{1}{2}}(x)\xi|\dvr_\xi\rA(x)\nabla_\xi \ln |\rA^{-\frac{1}{2}}(x)\xi|\,d\xi
=\ln r\int\limits_{\sE_r(x)} \frac{ds}{|\rA^{-1}(x)\xi|} -
\int\limits_{\sE_r(x)\setminus\sE_1(x)} \frac{d\xi}{|\rA^{-\frac{1}{2}}(x)\xi|^2}
\\
=&\ln r\int\limits_{\sE_1(x)} \frac{ds}{|\rA^{-1}(x)\xi|} -\sqrt{\det \rA(x)}
\int\limits_{B_r(0)\setminus B_1(0)} \frac{d y}{|y|^2}=\ln r\int\limits_{\sE_1(x)} \frac{ds}{|\rA^{-1}(x)\xi|} -2\pi\sqrt{\det \rA(x)}\ln r
\end{align*}
and we arrive at the statement of the lemma for $n=2$. The proof is complete.
\end{proof}

Estimates (\ref{3.50a}) allow us to prove one more auxiliary lemma, which will be used then in the proof of our main theorems.

\begin{lemma}\label{lm5.2} The family $\b_\e$ is uniformly bounded in $L_\infty(\Om)$.
Under Assumption~\ref{A7}, the limit $\b$ of the family $\b_\e$ in the space $\mathfrak{M}$ is an element of $L_\infty(\Om)$.
\end{lemma}

\begin{proof} The family $\b_\e$ is uniformly bounded in $L_\infty(\Om)$ due to its definition (\ref{2.9}) and estimates (\ref{3.50a}). Since the space $L_\infty(\Om)$ is dual to $L_1(\Om)$, there exists a sequence $\e'$ such that $\b_{\e'}$ converges weakly in $L_\infty(\Om)$ to some limit $\tilde{\b}\in L_\infty(\Om)$. Hence,
\begin{equation*}
(\b_{\e'}u,v)_{L_2(\Om)}=\int\limits_{\Om} \b_{\e'} u\overline{v}\,dx\to \int\limits_{\Om} \tilde{\b} u\overline{v}\,dx,\qquad \e'\to0,
\end{equation*}
for all $u,v\in C^\infty(\overline{\Om})$ vanishing on $\p\Om$. At the same time, it follows from definition (\ref{2.3}) of the norm in $\mathfrak{M}$ that $(\b_{\e'}u,v)_{L_2(\Om)}\to\la \b u,v\ra$, $\e'\to0$. Hence, $\la\b u,v\ra= (\tilde{\b} u,v)_{L_2(\Om)}$, and due to the density of the functions from $C^\infty(\overline{\Om})$ vanishing on $\p\Om$ in $\Ho_2^1(\Om)$ and $\Ho_2^1(\Om)\cap W_2^2(\Om)$. Therefore, $\b=\tilde{\b}$ and this proves the lemma.
\end{proof}

\section{Operator estimates}

In this section we prove Theorems~\ref{th1},~\ref{th2}. The proofs consist of three main steps. At the first step we prove estimates (\ref{2.16}), (\ref{2.21}). At the second step we establish estimates (\ref{2.17}), (\ref{2.23}). And at third step  we show the order sharpness of the certain terms in the estimates.

\subsection{$W_2^1$-estimates: general case}\label{sec4.1}

In this subsection we prove estimate (\ref{2.16}). We choose an arbitrary $f\in L_2(\Om)$ and by means of the solutions of problems (\ref{2.7}), (\ref{2.11}) we define $v_\e:=u_\e- u_0 \Xi_\e$. In view of the definition  of the function $\Xi_\e$, the function $v_\e$ belongs to $W_2^1(\Om^\e)$ and satisfies the boundary conditions
\begin{equation}\label{4.25}
v_\e=0\quad\text{on}\quad\p\tht_D^\e,\qquad v_\e=u_\e\quad\text{on}\quad \p\tht_R^\e.
\end{equation}
We write integral identity (\ref{2.18}) choosing $v_\e$ as the test function:
\begin{equation}\label{4.12}
\hf(u_\e,v_\e)-\l(u_\e,v_\e)_{L_2(\Om^\e)}+ \big(a^\e(\,\cdot\,,u_\e),v_\e\big)_{L_2(\p\tht_R^\e)} =(f,v_\e)_{L_2(\Om^\e)}.
\end{equation}
Then we multiply the equation in (\ref{2.11}) by $\overline{v_\e}\Xi_\e$  and integrate once by parts over $\Om^\e$:
\begin{equation}\label{4.13}
\hf(u_0,v_\e\overline{\Xi_\e})-\left(\frac{\p u_0}{\p\boldsymbol{\nu}},v_\e\overline{\Xi_\e}\right)_{L_2(\p\tht_\e^R)} +\g(\b\Ups u_0\Xi_\e,v_\e)_{L_2(\Om^\e)}-\l(u_0\Xi_\e,v_\e)_{L_2(\Om^\e)} =(f,v_\e\overline{\Xi_\e})_{L_2(\Om^\e)}.
\end{equation}
It follows from the definition of the form $\hf$ that
\begin{equation}\label{4.14}
\hf(u_0,v_\e\overline{\Xi_\e})=
\hf(u_0\Xi_\e,v_\e)+
\big(\rA\nabla u_0,v_\e\nabla \overline{\Xi_\e}\big)_{L_2(\Om^\e)}
-
\big(\rA u_0\nabla\Xi_\e,\nabla v_\e\big)_{L_2(\Om^\e)}
-\sum\limits_{j=1}^{n} \left(A_j u_0 \frac{\p \Xi_\e}{\p x_j},v_\e\right)_{L_2(\Om^\e)}.
\end{equation}
Having this identity and (\ref{4.25}) in mind, we calculate the difference of identities (\ref{4.12}), (\ref{4.13}) and take then the real part of the resulting relation:
\begin{equation}\label{4.15}
\RE\hf(v_\e,v_\e)-\RE\l\|v_\e\|_{L_2(\Om^\e)}^2+ \RE \big(a^\e(\,\cdot\,,u_\e)-a^\e(\,\cdot\,,u_0\Xi_\e),v_\e\big)_{L_2(\p\tht_R^\e)} =\RE h_\e,
\end{equation}
where
\begin{align*}
& h_\e:=h_\e^{(1)}+h_\e^{(2)}+h_\e^{(3)}+h_\e^{(4)}, 
\\
&h_\e^{(1)}:=  (f,(1-\overline{\Xi_\e})v_\e)_{L_2(\Om^\e)} + \g(\b \Ups u_0,(\overline{\Xi_\e}-1)v_\e)_{L_2(\Om^\e)},\nonumber
\\
&h_\e^{(2)}:= -\sum\limits_{k\in\mathds{M}^\e}
\big((\rA-\rA(M_k^\e))u_0\nabla \Xi_\e,
\nabla v_\e
\big)_{L_2(E_k^\e\setminus\om_k^\e)},\nonumber
\\
&
h_\e^{(3)}:=  \sum\limits_{k\in\mathds{M}^\e} h_{\e,k}^{(3)},
\qquad h_\e^{(4)}:=  \sum\limits_{k\in\mathds{M}^\e} h_{\e,k}^{(4)}+\g(\b \Ups u_0,v_\e)_{L_2(\Om^\e)},
\\
&
\begin{aligned}
h_{\e,k}^{(3)}:=&-\sum\limits_{j=1}^{n}\left(A_j u_0 \frac{\p \Xi_\e}{\p x_j},v_\e\right)_{L_2(E_k^\e\setminus \om_k^\e)}+
\big((\rA +\rA(M_k^\e))\nabla u_0,v_\e\nabla \overline{\Xi_\e}
\big)_{L_2(E_k^\e\setminus\om_k^\e)}
\\
&
-
(\rA(M_k^\e)\nabla u_0\cdot\nu^\e,v_\e\overline{\Xi_\e})_{L_2(\p\om_k^\e)},
\end{aligned}
\\
&
h_{\e,k}^{(4)}:=-\big(\rA(M_k^\e) u_0\nabla\Xi_\e, \nabla v_\e\big)_{L_2(E_k^\e\setminus\om_k^\e)} - \big(\rA(M_k^\e) \nabla u_0,v_\e\nabla\overline{\Xi_\e}\big)_{L_2(E_k^\e\setminus\om_k^\e)}-
\big(a^\e(\,\cdot\,,u_0\Xi_\e),v_\e\big)_{L_2(\p\om_k^\e)};
\end{align*}
we recall that $\nu^\e$ is the unit normal to $\p\om_k^\e$ directed inside $\om_k^\e$.

Letting $\tht_{R,i}^\e:=\bigcup\limits_{k\in\mathds{M}_{R,i}^\e} \om_k^\e$, $i=1,2$, we observe that the function $\Xi_\e$ vanishes on $\p\tht_{R,1}^\e$. By Assumption~\ref{A6} and (\ref{2.6}) we find:
\begin{align*}
\RE\big(a^\e(\,\cdot\,,u_\e)-a^\e(\,\cdot\,,u_0\Xi_\e), v_\e\big)_{L_2(\p\tht_R^\e)} =&\RE\big(a^\e(\,\cdot\,,v_\e),v_\e\big)_{L_2(\p\tht_{R,1}^\e)}
+\RE\big(a^\e(\,\cdot\,,u_\e)-a^\e(\,\cdot\,,u_0\Xi_\e), v_\e\big)_{L_2(\p\tht_{R,2}^\e)}
\\
\geqslant  & \mu_1 \|v_\e\|_{L_2(\p\tht_{R,1}^\e)}^2 - \mu_0 \|v_\e\|_{L_2(\p\tht_{R,2}^\e)}^2.
\end{align*}
Using then inequalities (\ref{4.3}), (\ref{3.22}) and Assumption~\ref{A6}, we find a lower bound for the left hand side of identity (\ref{4.15}):
\begin{equation}\label{4.20}
\RE\hf(v_\e,v_\e)-\RE\l\|v_\e\|_{L_2(\Om^\e)}^2+ \RE \big(a^\e(\,\cdot\,,u_\e)-a^\e(\,\cdot\,,u_0\Xi_\e),v_\e\big)_{L_2(\tht_R^\e)} \geqslant C \|v_\e\|_{W_2^1(\Om^\e)}^2 + \mu_1
\|v_\e\|_{L_2(\p\tht_{R,1}^\e)}^2.
\end{equation}
Hereinafter in this section by $C$ we denote inessential constants independent of $\e$, $x$, $k$, $u_\e$, $u_0$, $v_\e$, $f$  but, generally speaking, depending on $\l$. Our next key step is to estimate the right hand side in (\ref{4.15}) and to get in this way a bound for $\|v_\e\|_{W_2^1(\Om^\e)}$.

Taking into consideration (\ref{3.63}) and Lemma~\ref{lm5.2}, we can estimate  the function $h_\e^{(1)}$:
\begin{equation}\label{4.21}
\begin{aligned}
|h_\e^{(1)}|\leqslant &\sum\limits_{k\in\mathds{M}^\e} \big| (f,(1-\overline{\Xi_\e})v_\e)_{L_2(E_k^\e\setminus\om_k^\e)} \big|+\g
\sum\limits_{k\in\mathds{M}^\e}  \big| (\b \Ups u_0,(1-\overline{\Xi_\e})v_\e)_{L_2(E_k^\e\setminus\om_k^\e)}
\big|
\\
\leqslant &
\sum\limits_{k\in\mathds{M}^\e}
\|f\|_{L_2(E_k^\e \setminus\om_k^\e)} \|(\Xi_\e-1)v_\e\|_{L_2(E_k^\e \setminus\om_k^\e)} + C \sum\limits_{k\in\mathds{M}^\e} \|u_0\|_{L_2(E_k^\e \setminus\om_k^\e)} \|(\Xi_\e-1)v_\e\|_{L_2(E_k^\e \setminus\om_k^\e)}
\\
 \leqslant & C \e\sum\limits_{k\in\mathds{M}^\e}
\big(\|f\|_{L_2(E_k^\e\setminus\om_k^\e)}+\|u_0\|_{L_2(E_k^\e \setminus\om_k^\e)}\big)\|v_\e\|_{W_2^1(B_{\e R_3}(M_k^\e)\setminus\om_k^\e)}
\leqslant  C \e\|f\|_{L_2(\Om^\e)}   \|v_\e\|_{W_2^1(\Om^\e)}.
\end{aligned}
\end{equation}

The assumed smoothness of the functions $A_{ij}$ implies the inequality
\begin{equation*}
|\rA(x)-\rA(M_k^\e)|\leqslant C|x|\leqslant C\e \quad\text{a.e. in}\quad E_k^\e,
\end{equation*}
and this is why by (\ref{3.64}), (\ref{2.14}) we get:
\begin{equation}\label{4.49}
|h_\e^{(2)}|\leqslant C   \big(\eta^{\frac{n}{2}-1}\vk^{-\frac{1}{2}} +\e^{\frac{3}{2}}\eta^{\frac{1}{2}}
 +\e^2\big) \|u_0\|_{W_2^2(\Om)}\|\nabla v_\e\|_{L_2(\Om^\e)}.
\end{equation}
In order to estimate $h_\e^{(3)}$, we first integrate by parts using the properties of the function $\Xi_\e$:
\begin{align}\label{4.50}
&
\begin{aligned}
h_{\e,k}^{(3)}=& -\sum\limits_{j=1}^{n}\left(A_j u_0 \frac{\p (\Xi_\e-1)}{\p x_j},v_\e\right)_{L_2(E_k^\e\setminus \om_k^\e)}+
\big((\rA+\rA(M_k^\e))\nabla u_0,v_\e\nabla (\overline{\Xi_\e}-1)\big)_{L_2(E_k^\e\setminus\om_k^\e)}
\\
& -
(\rA(M_k^\e)\nabla u_0\cdot\nu,v_\e\overline{\Xi_\e})_{L_2(\p\om_k^\e)}
=h_{\e,k}^{(5)}+h_{\e,k}^{(6)},
\end{aligned}
\\
&\nonumber
\begin{aligned}
h_{\e,k}^{(5)}:=&-\sum\limits_{j=1}^{n}\left((\Xi_\e-1)A_j \nu_j u_0, v_\e\right)_{L_2(\p \om_k^\e)} +\big(((\Xi_\e-1)\rA-\rA(M_k^\e))\nabla u_0\cdot\nu,
v_\e\big)_{L_2(\p \om_k^\e)},
\end{aligned}
\\
&\nonumber
\begin{aligned}
h_{\e,k}^{(6)}:=&\sum\limits_{j=1}^{n}\left(A_j (\Xi_\e-1) \frac{\p u_0}{\p x_j},v_\e\right)_{L_2(E_k^\e\setminus \om_k^\e)}+ \sum\limits_{j=1}^{n}\left(A_j (\Xi_\e-1) u_0,\frac{\p v_\e}{\p x_j}\right)_{L_2(E_k^\e\setminus \om_k^\e)}
\\
&-
\big((\Xi_\e-1)\dvr(\rA+\rA(M_k^\e))\nabla u_0,
v_\e\big)_{L_2(E_k^\e\setminus\om_k^\e)} -
\big((\Xi_\e-1) (\rA+\rA(M_k^\e))\nabla u_0,
 \nabla v_\e
 \big)_{L_2(E_k^\e\setminus\om_k^\e)},
\end{aligned}
\end{align}
where $\nu_j$ are the components of the unit normal $\nu$.  Inequality (\ref{3.63}) applied with $v=u_0$ and $v=v_\e$ allows us to estimate $h_{\e,k}^{(6)}$:
\begin{equation}\label{4.51}
|h_{\e,k}^{(6)}|\leqslant C \e \|u_0\|_{W_2^2(B_{\e R_3}(M_k^\e)\setminus \om_k^\e)} \|v_\e\|_{W_2^1(B_{\e R_3}(M_k^\e)\setminus \om_k^\e)}.
\end{equation}
Inequalities (\ref{3.22}) and (\ref{3.66}) give rise to a similar estimate for $h_{\e,k}^{(5)}$:
\begin{equation*}
|h_{\e,k}^{(5)}|\leqslant C \e\eta\vk\|u_0\|_{W_2^2(B_{\e R_3}(M_k^\e)\setminus \om_k^\e)} \|v_\e\|_{W_2^1(B_{\e R_3}(M_k^\e)\setminus \om_k^\e)}.
\end{equation*}
This estimate and (\ref{4.51}), (\ref{4.50})  imply:
\begin{equation}\label{4.52}
|h_\e^{(3)}|\leqslant C\e\|u_0\|_{W_2^2(\Om^\e)} \|v_\e\|_{W_2^1(\Om^\e)}.
\end{equation}

We proceed to estimating the function $h_{\e,k}^{(4)}$, which is one of the most non-trivial steps in the proof. As above, we first integrate by parts in $h_{\e,k}^{(4)}$ taking into consideration the definition of $\Xi_\e$ and the  equation for $X_{k,\e}$:
\begin{equation}\label{4.53}
h_{\e,k}^{(4)}=-(\rA(M_k^\e) u_0\nabla\Xi_\e\cdot \nu, v_\e)_{L_2(\p E_k^\e)} - (\rA(M_k^\e) u_0\nabla\Xi_\e\cdot \nu^\e, v_\e)_{L_2(\p \om_k^\e)} -
\big(a^\e(\,\cdot\,,u_0\Xi_\e),v_\e\big)_{L_2(\p\om_k^\e)},
\end{equation}
where $\nu$ stands for the unit outward normal to $\p E_k^\e$. We fix $k\in\mathds{M}_R^\e$ and represent the functions $u_0$ and $v_\e$ as
\begin{equation}\label{4.58}
u_0=\la u_0\ra  + u_0^\bot,\qquad v_\e=\la v_\e\ra  + v_\e^\bot,
\end{equation}
where the operations $\la\cdot\ra$ and $\cdot^\bot$ were defined in (\ref{4.56}). Then by inequality (\ref{3.6}) with $\p\om_k^\e=\p E_k^\e$ and $\eta=1$ and by Lemma~\ref{lm3.3} the functions $u_0$, $u_0^\bot$, $v_\e^\bot$ satisfy the estimates
\begin{equation}\label{4.57}
\begin{gathered}
\|u_0\|_{L_2(\p E_k^\e)}\leqslant C\e^{-\frac{1}{2}} \|u_0\|_{W_2^1(B_{\e R_3}(M_k^\e)\setminus\om_k^\e)},\qquad
\|u_0^\bot\|_{L_2(\p E_k^\e)}\leqslant C\e^{\frac{1}{2}} \|\nabla u_0\|_{L_2(B_{\e R_3}(M_k^\e)\setminus\om_k^\e)},
\\
\|v_\e^\bot\|_{L_2(\p E_k^\e)}\leqslant C\e^{\frac{1}{2}} \|\nabla v_\e\|_{L_2(B_{\e R_3}(M_k^\e)\setminus\om_k^\e)}.
\end{gathered}
\end{equation}
By identities (\ref{4.58}), we rewrite the first term in formula (\ref{4.53}) as
\begin{equation}\label{4.59}
\begin{aligned}
(\rA(M_k^\e) u_0\nabla\Xi_\e\cdot \nu, v_\e)_{L_2(\p E_k^\e)}=&\la u_0\ra \la \overline{v_\e}\ra \int\limits_{\p E_k^\e} \rA(M_k^\e)\nabla\Xi_\e\cdot \nu\,ds +
(\rA(M_k^\e) u_0\nabla \Xi_\e\cdot \nu, v_\e^\bot)_{L_2(\p E_k^\e)}
\\
&+
\la \overline{v_\e}\ra \int\limits_{\p E_k^\e} \rA(M_k^\e) u_0^\bot \Xi_\e\cdot \nu\,ds.
\end{aligned}
\end{equation}
The second and the third term in the right hand side  can be estimated by means of  (\ref{3.66}),
 (\ref{4.54}),  (\ref{4.57}), (\ref{3.23}):
\begin{equation}\label{4.60}
\begin{aligned}
& \big|(\rA(M_k^\e) u_0\nabla \Xi_\e\cdot \nu, v_\e^\bot)_{L_2(\p E_k^\e)}\big|\leqslant C \e \|u_0\|_{W_2^1(B_{\e R_3}(M_k^\e)\setminus\om_k^\e)} \|\nabla v_\e\|_{L_2(B_{\e R_3}(M_k^\e)\setminus\om_k^\e)},
\\
&\Bigg|\la \overline{v_\e}\ra \int\limits_{\p E_k^\e} \rA(M_k^\e) u_0^\bot \Xi_\e\cdot \nu\,ds\Bigg|\leqslant C \e^{\frac{n+2}{2}} \big|\la \overline{v_\e}\ra \big| \|\nabla u_0\|_{L_2(B_{\e R_3}(M_k^\e)\setminus\om_k^\e)}
\\
&\hphantom{\Bigg|\la \overline{v_\e}\ra \int\limits_{\p E_k^\e} u_0^\bot \rA(M_k^\e)\nabla \Xi_\e\cdot \nu\,ds\Bigg|} \leqslant C\e \|\nabla u_0\|_{L_2(B_{\e R_3}(M_k^\e)\setminus\om_k^\e)} \|v_\e\|_{L_2(B_{\e R_3}(M_k^\e)\setminus\om_k^\e)}.
\end{aligned}
\end{equation}
We calculate the integral in first term in the right hand side in (\ref{4.59}) by using (\ref{4.54}), (\ref{3.23}):
\begin{align*}
&\int\limits_{\p E_k^\e}  \rA(M_k^\e)\nabla \Xi_\e\cdot \nu\,ds = \frac{(2-n)K_{k,\e}}{R_4^{n-2}} \eta^{n-2}  \int\limits_{\p E_k^\e} |\rA^{-1}(M_k^\e)(x-M_k^\e)|^{-1}\,ds+O(\e^{n+2}\eta^2) && \text{as}\qquad n\geqslant 3,
\\
&\int\limits_{\p E_k^\e}  \rA(M_k^\e)\nabla \Xi_\e\cdot \nu\,ds = -\frac{1}{\ln\eta} \int\limits_{\p E_k^\e} |\rA^{-1}(M_k^\e)(x-M_k^\e)|^{-1}\,ds+O(\e^4\eta)&& \text{as}\qquad n=2,
\end{align*}
where the $O$-terms are uniform in  $k$.
In the integrals in the right hand sides of the above identities we make the change of variables
$y=\e^{-1}R_4^{-1}(x-M_k^\e)$ and then we get:
\begin{equation*}
\frac{1}{R_4^{n-2}}\int\limits_{\p E_k^\e} |\rA^{-1}(M_k^\e)(x-M_k^\e)|^{-1}\,ds=\e^n\Ups(M_k^\e), \qquad n\geqslant 2.
\end{equation*}
Hence,
\begin{equation}\label{4.65}
\begin{aligned}
&\int\limits_{\p E_k^\e}  \rA(M_k^\e)\nabla \Xi_\e\cdot \nu\,ds = \e^n\g K_{k,\e}(2-n)\Ups(M_k^\e)+O\big(\e^n|\eta^{n-2}\e^{-2}-\g|+\e^{n+2}\eta^2\big) && \text{as}\quad n\geqslant 3,
\\
&\int\limits_{\p E_k^\e}  \rA(M_k^\e)\nabla \Xi_\e\cdot \nu\,ds = \e^2\g\Ups(M_k^\e)+O\big(\e^2|\e^{-2}\ln^{-1}\eta-\g|+\e^4\eta\big) && \text{as}\quad  n=2.
\end{aligned}
\end{equation}
It is clear that
\begin{align*}
&(\Ups u_0,v_\e)_{L_2(B_{\e R_3}(M_k^\e)\setminus\om_k^\e)}= \la u_0\ra \la \overline{v_\e}\ra \int\limits_{B_{\e R_3}(M_k^\e)\setminus\om_k^\e} \Ups\,dx+ \la \overline{v_\e}\ra \int\limits_{B_{\e R_3}(M_k^\e)\setminus\om_k^\e} \Ups u_0^\bot\,dx
\\
&\hphantom{(\Ups u_0,v_\e)_{L_2(B_{\e R_3}(M_k^\e)\setminus\om_k^\e)}= }+ (\Ups u_0,v_\e^\bot)_{L_2(B_{\e R_3}(M_k^\e)\setminus\om_k^\e)},
\\
&\Bigg|\int\limits_{B_{\e R_3}(M_k^\e)\setminus\om_k^\e} \Ups\,dx-  \Ups(M_k^\e)\mes_n B_{\e R_3}(M_k^\e)
\Bigg|\leqslant C\e^{n+1}.
\end{align*}
Then by Lemma~\ref{lm3.3}   the estimates
\begin{align*}
& \big|(\Ups u_0,v_\e^\bot)_{L_2(B_{\e R_3}(M_k^\e)\setminus\om_k^\e)}\big| \leqslant C \e \|u_0\|_{L_2(B_{\e R_3}(M_k^\e)\setminus\om_k^\e)} \|\nabla v_\e\|_{L_2(B_{\e R_3}(M_k^\e)\setminus\om_k^\e)},
\\
& \Bigg| \la \overline{v_\e}\ra \int\limits_{B_{\e R_3}(M_k^\e)\setminus\om_k^\e} \Ups u_0^\bot\,dx\Bigg| \leqslant C \|u_0^\bot\|_{L_2(B_{\e R_3}(M_k^\e)\setminus\om_k^\e)} \|v_\e\|_{L_2(B_{\e R_3}(M_k^\e)\setminus\om_k^\e)}
\\
& \hphantom{\Bigg| \la \overline{v_\e}\ra \int\limits_{B_{\e R_3}(M_k^\e)\setminus\om_k^\e} u_0^\bot.dx\leqslant}
\leqslant C \e \|\nabla u_0\|_{L_2(B_{\e R_3}(M_k^\e)\setminus\om_k^\e)}
\|v_\e\|_{L_2(B_{\e R_3}(M_k^\e)\setminus\om_k^\e)},
\end{align*}
hold and therefore,
\begin{equation*}
\Big| (\Ups u_0,v_\e)_{L_2(B_{\e R_3}(M_k^\e)\setminus\om_k^\e)}- \la u_0\ra \la \overline{v_\e}\ra \Ups(M_k^\e)\mes_n B_{\e R_3}(M_k^\e)
\Big|\leqslant C \e \|u_0\|_{W_2^1(B_{\e R_3}(M_k^\e)\setminus\om_k^\e)}\|v_\e\|_{W_2^1(B_{\e R_3}(M_k^\e)\setminus\om_k^\e)}.
\end{equation*}
The above estimate and  (\ref{4.59}), (\ref{4.60}), (\ref{4.65}) imply:
\begin{equation}\label{4.62}
(\rA(M_k^\e) u_0\nabla\Xi_\e\cdot \nu, v_\e)_{L_2(\p E_k^\e)}=h_{\e,k}^{(7)}+h_{\e,k}^{(8)},
\end{equation}
where
\begin{align*}
&h_{\e,k}^{(7)}:=-\frac{K_{k,\e}(2-n)\g}{R_3^n\mes_n B_1(0)} (\Ups u_0,v_\e)_{L_2(B_{\e  R_3}(M_k^\e)\setminus\om_k^\e)}
&&\hspace{-2.5 true cm}\text{as}\quad n\geqslant 3,
\\
&h_{\e,k}^{(7)}:=-\frac{\g}{R_3^2\mes_n B_1(0)}(\Ups u_0,v_\e)_{L_2(B_{\e  R_3}(M_k^\e)\setminus\om_k^\e)}
 &&\hspace{-2.5 true cm}\text{as}\quad n=2,
\end{align*}
while the functions $h_{\e,k}^{(8)}$ obey the estimates
\begin{equation}\label{4.64}
\begin{aligned}
|h_{\e,k}^{(8)}|\leqslant& C\e\|u_0\|_{W_2^1(B_{\e R_3}(M_k^\e)\setminus\om_k^\e)}\|v_\e\|_{W_2^1(B_{\e R_3}(M_k^\e)\setminus\om_k^\e)}
\\
&+C\big|\eta^{n-2}\e^{-2}\vk^{-1}-\g\big| \|u_0\|_{L_2(B_{\e R_3}(M_k^\e)\setminus\om_k^\e)}\|v_\e\|_{L_2(B_{\e R_3}(M_k^\e)\setminus\om_k^\e)}.
\end{aligned}
\end{equation}
We observe that
\begin{equation}\label{4.66}
h_{\e,k}^{(7)}=-\g(\Ups\b_\e u_0,v_\e)_{L_2(B_{\e  R_3}(M_k^\e)\setminus\om_k^\e)}.
\end{equation}

We proceed to estimating two other terms in the right hand side of (\ref{4.53}). We first of all note  that as $k\in\mathds{M}_D^\e$, this term vanishes since $v_\e=0$ on $\p\om_k^\e$ for such $k$. This is why we need to estimate it only for $k\in\mathds{M}_R^\e$. We first consider the case $k\in\mathds{M}_{R,1}^\e$. The function $\Xi_\e$ vanishes on $\p\om_k^\e$ and hence, in view of the identity in (\ref{2.6}), the function $a^\e(\,\cdot\,,\Xi_\e v_\e)$ vanishes and the same is true for the third term in  the right hand side of (\ref{4.53}).
 Since $Z_k^\e=0$ on $\p\om_k^\e$ and $Z_k^\e\in C^1(\overline{B_{\e R_3(M_k^\e)}\setminus \om_k^\e})$, we have
\begin{equation*}
\frac{\p Z_k^\e}{\p x_j}=\frac{\p\tau}{\p x_j} \frac{\p Z_k^\e}{\p\tau}\quad\text{on}\quad\p\om_k^\e;
\end{equation*}
where $\tau$ is the distance measured along the unit normal to $\p\om_k^\e$, see Assumption~\ref{A1}. According to this assumption, the derivatives of $\tau$ in $x_j$ are bounded uniformly in $k$, $\e$ and the spatial variables. By  (\ref{3.90}), (\ref{3.82})  we then obtain:
\begin{equation*}
\bigg|\frac{\p Z_k^\e}{\p x_j}\bigg|\leqslant C (\e\eta)^{-1}\quad\text{on}\quad \p\om_k^\e,\qquad k\in\mathds{M}_R^\e.
\end{equation*}
Hence, by (\ref{3.22}), (\ref{2.6}),
\begin{equation}\label{4.70}
\begin{aligned}
&\big|(\rA(M_k^\e) u_0\nabla\Xi_\e\cdot \nu, v_\e)_{L_2(\p \om_k^\e)}\big| \leqslant C (\e\eta\vk)^{-\frac{1}{2}} \|u_0\|_{W_2^1(B_{\e R_3}(M_k^\e)\setminus\om_k^\e)} \|v_\e\|_{L_2(\p\om_k^\e)},
\\
&\big(a^\e(\,\cdot\,,u_0\Xi_\e),v_\e\big)_{L_2(\p\om_k^\e)} =\big(a^\e(\,\cdot\,,0),v_\e\big)_{L_2(\p\om_k^\e)}=0,
\end{aligned}
\end{equation}
as $k\in\mathds{M}_{R,1}^\e$. For $k\in\mathds{M}_{R,2}^\e$, by boundary condition (\ref{2.20}) and the definition of $\Xi_\e$ we see that
\begin{equation*}
\rA(M_k^\e) \nabla \Xi_\e \cdot \nu =\e^{-1}\eta^{-1} b_k\Xi_\e\quad\text{on}\quad \p\om_k^\e,
\end{equation*}
and therefore, in view of (\ref{2.22}),
\begin{equation*}
- (\rA(M_k^\e) u_0\nabla\Xi_\e\cdot \nu, v_\e)_{L_2(\p \om_k^\e)} -
\big(a^\e(\,\cdot\,,u_0\Xi_\e),v_\e\big)_{L_2(\p\om_k^\e)} = - \big(\tilde{a}_k^\e(\,\cdot\,,u_0\Xi_\e),v_\e\big)_{L_2(\p\om_k^\e)}.
\end{equation*}
Using then the estimate for $\tilde{a}_k^\e$ in (\ref{2.26}) as well as (\ref{3.22}), (\ref{3.66}), we get:
\begin{equation}\label{4.85}
\begin{aligned}
\big| (\rA(M_k^\e) u_0\nabla\Xi_\e\cdot \nu, v_\e)_{L_2(\p \om_k^\e)} +\big(a^\e(\,\cdot\,,u_0\Xi_\e),v_\e\big)_{L_2(\p\om_k^\e)}\big| \leqslant & C \e\eta\vk \mu_2 \|u_0\|_{W_2^1(B_{\e R_3}(M_k^\e)\setminus\om_{k,\e})}
\\
&\cdot\|v_\e\|_{W_2^1(B_{\e R_3}(M_k^\e)\setminus\om_{k,\e})}.
\end{aligned}
\end{equation}
Summing up the above estimates over $k\in\mathds{M}_{R,2}^\e$,
relations (\ref{4.70}) over $k\in\mathds{M}_{R,1}^\e$ and identities (\ref{4.62}), (\ref{4.66}) and inequalities (\ref{4.64}) over $k\in\mathds{M}^\e$, by (\ref{3.71}), we finally obtain:
\begin{align}\label{4.69}
&h_\e^{(4)}= \g((\b-\b_\e) u_0,\Ups v_\e)_{L_2(\Om^\e)} + h_\e^{(9)},
\\
&
\begin{aligned}
|h_\e^{(9)}|\leqslant  &C(\e+\e\eta\vk\mu_2)\|u_0\|_{W_2^1(\Om^\e)}\|v_\e\|_{W_2^1(\Om^\e)} + C (\e\eta\vk)^{-\frac{1}{2}} \|u_0\|_{W_2^1(\Om^\e)}   \|v_\e\|_{L_2(\p\tht_{R,1}^\e)}
\\
&+C\big|\eta^{n-2}\e^{-2}\vk^{-1}-\g\big| \|u_0\|_{L_2(\Om^\e)}\|v_\e\|_{L_2(\Om^\e)}.
\end{aligned}\label{4.68}
\end{align}

If $\g\ne 0$, then the first term in the left hand side of (\ref{4.69}) is, generally speaking, non-zero and we need Assumption~\ref{A7}  to estimate it. Under this assumption,
we continue the function $v_\e$ inside $\om_k^\e$ as follows. We first let
\begin{equation}\label{4.71}
v_\e:=0\quad\text{in}\quad \om_k^\e\quad\text{for}\quad k\in\mathds{M}_D^\e.
\end{equation}
For $k\in\mathds{M}_R^\e$ we introduce the quantities $\la v_\e\ra_k$ and the functions $v_{\e,k}$ by formulae (\ref{3.78}).
By Lemma~\ref{lm3.3} with $\e$ replaced by $\e\eta$ we have
\begin{equation}\label{4.84}
\|v_{\e,k}\|_{L_2(B_{\e\eta R_2}(M_k^\e)\setminus\om_k^\e)}^2
\leqslant C \e^2\eta^2 \|\nabla v_\e\|_{L_2(B_{\e\eta R_2}(M_k^\e)\setminus\om_k^\e)}^2.
\end{equation}
Then for $k\in\mathds{M}_R^\e$ we define the continuation of the function $v_\e$ inside $\om_k^\e$ in terms of the local variables $(\tau,s)$ as follows:
\begin{equation}\label{4.72}
\begin{aligned}
&v_\e(\tau,s):=\la v_\e\ra_k +v_{\e,k}(-\tau,s)\chi_{1}(\tau\e^{-1}\eta^{-1}) &&\text{for}\quad x\in\om_k^\e,\quad \dist(x,\p\om_k^\e)\leqslant\e\eta\tau_0,
\\
&v_\e(\tau,s):=\la v_\e\ra_k  && \text{for}\quad x\in\om_k^\e,\quad \dist(x,\p\om_k^\e)>\e\eta\tau_0,
\end{aligned}
\end{equation}
where $\chi_{1}$ is the cut-off function introduced in the proof of Lemma~\ref{lm4.2}. It is obvious that this continuation gives a function in $W_2^1(B_{\e R_3}(M_k^\e))$, which in view of  estimates (\ref{3.72}), (\ref{4.84}) satisfies the inequalities
\begin{align}\label{4.73}
&
\begin{aligned}
\|v_\e\|_{L_2(\om_k^\e)}^2\leqslant & C \left(|\la v_\e\ra_k|^2\mes_n \om_k^\e
+\|v_{\e,k}\|_{L_2(B_{\e\eta R_2}(M_k^\e)\setminus\om_k^\e)}^2
\right)
\\
\leqslant &C \|v_\e\|_{L_2(B_{\e\eta R_2}(M_k^\e)\setminus\om_k^\e)}^2 \leqslant   C  \e^2\eta^2\vk
\|v_\e\|_{W_2^1(B_{\e R_3}(M_k^\e)\setminus\om_k^\e)}^2,
\end{aligned}
\\
&\label{4.74}
\begin{aligned}
\|\nabla v_\e\|_{L_2(\om_k^\e)}^2\leqslant & C \Big(\|\nabla v_\e\|_{L_2(B_{\e\eta R_2}(M_k^\e)\setminus\om_k^\e)}^2+ \e^{-2}\eta^{-2}\|v_{\e,k}\|_{L_2(B_{\e\eta R_2}(M_k^\e)\setminus\om_k^\e)}^2\Big)
\\
\leqslant & C
 \|v_\e\|_{W_2^1(B_{\e R_3}(M_k^\e)\setminus\om_k^\e)}^2.
\end{aligned}
\end{align}
Due to these inequalities and identity (\ref{4.71}), the continued function $v_\e$, regarded as defined on the entire domain $\Om$, is an element of $\Ho_2^1(\Om)$ and
\begin{equation}\label{4.75}
\|v_\e\|_{L_2(\tht^\e)}^2\leqslant C\e^2\eta^2\vk
\|v_\e\|_{\Om^\e}^2,
\qquad
\|\nabla v_\e\|_{L_2(\tht^\e)}^2\leqslant C \|\nabla v_\e\|_{L_2(\Om^\e)}^2.
\end{equation}
These inequalities and (\ref{3.9}) allow us to rewrite the scalar product in the right hand side of (\ref{4.69}) as
\begin{equation}\label{4.76}
((\b-\b_\e) u_0,\Ups v_\e)_{L_2(\Om^\e)}=((\b-\b_\e) u_0,\Ups v_\e)_{L_2(\Om)}+h_\e^{(10)},
\end{equation}
where $h_\e^{(10)}$ is a function satisfying the estimate
\begin{equation}\label{4.77}
|h_\e^{(10)}|\leqslant C\|u_0\|_{L_2(\tht^\e)} \|v_\e\|_{L_2(\tht^\e)} \leqslant C
\e^2\eta^2\vk \|u_0\|_{W_2^1(\Om)}\|v_\e\|_{W_2^1(\Om^\e)}.
\end{equation}
In the scalar product in the right hand of (\ref{4.76}) the function $\Ups v_\e$ is an element of $\Ho_2^1(\Om)$ and then the function $(\b-\b_\e) u_0$ can be regarded as a functional on this space. Then by formula (\ref{2.3}), Assumption~\ref{A7} and inequalities (\ref{4.75}) we can estimate this scalar product as follows:
\begin{equation}\label{4.78}
\Big|((\b-\b_\e) u_0,\Ups v_\e)_{L_2(\Om)}\Big| \leqslant  C\|\b_\e-\b\|_{\mathfrak{M}} \|u_0\|_{W_2^1(\Om)} \|\Ups v_\e\|_{W_2^1(\Om)}
\leqslant   C\|\b_\e-\b\|_{\mathfrak{M}} \|u_0\|_{W_2^2(\Om)}\|v_\e\|_{W_2^1(\Om)}.
\end{equation}
The above estimate and (\ref{4.21}), (\ref{4.49}), (\ref{4.52}), (\ref{4.69}), (\ref{4.68}), (\ref{4.76}), (\ref{4.77}), (\ref{2.15}) yield a final estimate for the right hand side in (\ref{4.15}):
\begin{equation}\label{4.79}
\begin{aligned}
|\RE h_\e|\leqslant& |h_\e|\leqslant C \|f\|_{L_2(\Om)}\Big((\e+\e\eta\vk\mu_2)\|v_\e\|_{W_2^1(\Om^\e)} +  \big|\eta^{n-2}\e^{-2}\vk^{-1}(\e)-\g\big| \|v_\e\|_{L_2(\Om^\e)}
\\
&\hphantom{|h_\e|\leqslant C \|f\|_{L_2(\Om)}\Big(}+
\g\|\b_\e-\b\|_{\mathfrak{M}} \|v_\e\|_{W_2^1(\Om^\e)}+(\e\eta\vk)^{-\frac{1}{2}}  \|v_\e\|_{L_2(\p\tht_{R,1}^\e)}  \Big)
\\
\leqslant &\d \big(\|v_\e\|_{W_2^1(\Om^\e)}^2+\mu_1 \|v_\e\|_{L_2(\p\tht_{R,1}^\e)}^2\big)
\\
&+
C(\d)\big(\e^2+(\e\eta\vk\mu_2)^2 +\big|\eta^{n-2}\e^{-2}\vk^{-1} -\g\big|^2+(\e\eta\vk\mu_1)^{-1} +
\g^2\|\b_\e-\b\|_{\mathfrak{M}}^2
\big)\|f\|_{L_2(\Om)}^2,
\end{aligned}
\end{equation}
where $\d>0$ is arbitrary but fixed, while $C(\d)$ is a constant independent of $\e$, $f$, $u_0$ and $v_\e$. Substituting the above estimate with a sufficiently small $\d$ into the left hand side of (\ref{4.15}) and employing then (\ref{4.20}), we obtain:
\begin{equation}\label{4.80}
\|v_\e\|_{W_2^1(\Om^\e)}^2+\mu_1 \|v_\e\|_{L_2(\p\tht_{R,1}^\e)}^2 \leqslant  C\big(\e^2+(\e\eta\vk\mu_2)^2
+\big|\eta^{n-2}\e^{-2}\vk^{-1} -\g\big|^2 + (\e\eta\vk\mu_1)^{-1} +
\g^2\|\b_\e-\b\|_{\mathfrak{M}}^2
\big)\|f\|_{L_2(\Om)}^2.
\end{equation}
This estimate implies (\ref{2.16}).

If the set $\mathds{M}_R$ is empty, then the second term in (\ref{4.53}) is zero for $k\in\mathds{M}_R$ just because the function $v_\e$ vanishes on $\p\om_k^\e$. Then   estimates (\ref{4.70}), (\ref{4.85}) are no longer needed.   Estimate (\ref{4.68}) also simplifies:
\begin{equation*}
|h_\e^{(9)}|\leqslant  C\e\|u_0\|_{W_2^1(\Om^\e)} \|v_\e\|_{W_2^1(\Om^\e)}  +C\big|\eta^{n-2}\e^{-2}\vk^{-1} - \g\big| \|u_0\|_{L_2(\Om^\e)}\|v_\e\|_{L_2(\Om^\e)}.
\end{equation*}
We also do not need continuation (\ref{4.72}) and, hence, Assumption~\ref{A6}. Estimates (\ref{4.75}) are also omitted in the considered case. Estimates (\ref{4.78}), (\ref{4.79}) then become
\begin{align*}
&
\Big|((\b-\b_\e) u_0,\Ups v_\e)_{L_2(\Om)}\Big| \leqslant  C\|\b_\e-\b\|_{\mathfrak{M}} \|u_0\|_{W_2^2(\Om)} \|v_\e\|_{W_2^1(\Om^\e)},
\\
& |\RE h_\e|
\leqslant  \d \|v_\e\|_{W_2^1(\Om^\e)}^2
+
C(\d)\Big(\e^2+\big|\eta^{n-2}\e^{-2}\vk^{-1} -\g\big|^2 +
\|\b_\e-\b\|_{\mathfrak{M}}^2
\Big)\|f\|_{L_2(\Om)}^2.
\end{align*}
All other above arguing remain the same and we  arrive at estimate (\ref{2.16}) without the terms $(\e\eta\mu_1\vk)^{-\frac{1}{2}}$ and $\e\eta\vk\mu_2$.

\subsection{$W_2^1$-estimates: $\g=0$}

Here we prove estimate (\ref{2.21}). Assume that $\g=0$ and only Assumption~\ref{A1} holds.
In this case, the function $h_\e$ can be written as $h_\e=g_\e[v_\e]$, where
\begin{equation}\label{5.22}
\begin{aligned}
g_\e[v]:=&
(f,(1-\overline{\Xi_\e})v)_{L_2(\Om^\e)}
-\sum\limits_{k\in\mathds{M}^\e}
\Big(\big(\rA u_0\nabla \Xi_\e,
\nabla v
\big)_{L_2(E_k^\e\setminus\om_k^\e)}
- \big(\rA\nabla u_0, v
\nabla \overline{\Xi_\e}
\big)_{L_2(E_k^\e\setminus\om_k^\e)}
\Big)
\\
&
-\sum\limits_{k\in\mathds{M}^\e}(\rA \nabla u_0\cdot\nu^\e,v_\e\overline{\Xi_\e})_{L_2(\p\om_k^\e)}
-\sum\limits_{k\in\mathds{M}^\e}
\sum\limits_{j=1}^{n}\left(A_j u_0 \frac{\p \Xi_\e}{\p x_j},v\right)_{L_2(E_k^\e\setminus \om_k^\e)} .
\end{aligned}
\end{equation}
We also integrate by parts:
\begin{equation*}
\sum\limits_{k\in\mathds{M}^\e}
\big(\rA\nabla u_0, v
\nabla \overline{\Xi_\e}
\big)_{L_2(E_k^\e\setminus\om_k^\e)}= -\sum\limits_{k\in\mathds{M}^\e}
\int\limits_{k\in\mathds{M}^\e} (\Xi_\e-1) \dvr \rA \overline{v} \nabla u_0\,dx + \sum\limits_{k\in\mathds{M}^\e}
(\rA \nabla u_0\cdot\nu^\e,(\overline{\Xi_\e}-1)v_\e)_{L_2(\p\om_k^\e)}
\end{equation*}
and this allows us to rewrite (\ref{5.22}):
\begin{align*}
g_\e[v]:=&
(f,(1-\overline{\Xi_\e})v)_{L_2(\Om^\e)} -\sum\limits_{k\in\mathds{M}^\e}
\sum\limits_{j=1}^{n}\left(A_j u_0 \frac{\p \Xi_\e}{\p x_j},v\right)_{L_2(E_k^\e\setminus \om_k^\e)} + \sum\limits_{k\in\mathds{M}^\e} \big(\rA \nabla u_0\cdot\nu^\e,(\overline{\Xi_\e}-1)v_\e\big)_{L_2(\p\om_k^\e)}
\\
&
-\sum\limits_{k\in\mathds{M}^\e}
\big(\rA u_0\nabla \Xi_\e,
\nabla v
\big)_{L_2(E_k^\e\setminus\om_k^\e)}
-\sum\limits_{k\in\mathds{M}^\e}
\int\limits_{k\in\mathds{M}^\e} (\Xi_\e-1) \dvr \rA \overline{v} \nabla u_0\,dx.
\end{align*}
Then the function $h_\e$ can be directly estimated by means of inequality inequalities (\ref{2.15}), (\ref{3.66}), (\ref{3.6}) and  Lemma~\ref{lm4.3}:
\begin{align*}
|h_\e|=&|g_\e[v_\e]|\leqslant C \Big(\e\|f\|_{L_2(\Om)}+\big(\eta^{\frac{n}{2}-1}\e^{-1}\vk^{-\frac{1}{2}}
+\e^\frac{1}{2}\eta^\frac{1}{2}+\e\big)\|u_0\|_{W_2^2(\Om)}\Big)
\|v_\e\|_{W_2^1(\Om^\e)}
\\
\leqslant&  C\big( \eta^{\frac{n}{2}-1}\e^{-1}\vk^{-\frac{1}{2}} +\e^\frac{1}{2}\eta^\frac{1}{2}+\e\big)\|f\|_{L_2(\Om)}
\|v_\e\|_{W_2^1(\Om^\e)}.
\end{align*}
Substituting this estimate into the right hand side of (\ref{4.15})
 and using (\ref{4.20}), we obtain
\begin{equation}\label{4.83}
\|v_\e\|_{W_2^1(\Om^\e)}  \leqslant  C\big( \eta^{\frac{n}{2}-1}\e^{-1}\vk^{-\frac{1}{2}} +\e^\frac{1}{2}\eta^\frac{1}{2}+\e \big) \|f\|_{L_2(\Om)}.
\end{equation}
By Lemma~\ref{lm4.3} and inequality (\ref{2.15}) we also have:
\begin{align}
&\|(\Xi_\e-1)u_0\|_{L_2(\Om^\e)} \leqslant C (\e^2+\e\eta) \|u_0\|_{W_2^1(\Om^\e)} \leqslant C\e \|f\|_{L_2(\Om)},\label{4.1}
\\
&\|\nabla (\Xi_\e-1)u_0\|_{L_2(\Om^\e)} \leqslant
\| (\Xi_\e-1)\nabla u_0\|_{L_2(\Om^\e)} + \|u_0\nabla \Xi_\e\|_{L_2(\Om^\e)}\nonumber
\\
&\hphantom{\|\nabla (\Xi_\e-1)u_0\|_{L_2(\Om^\e)}}\leqslant C \big( \eta^{\frac{n}{2}-1}\e^{-1}\vk^{-\frac{1}{2}} +\e^\frac{1}{2}\eta^\frac{1}{2}+\e\big) \|f\|_{L_2(\Om)}.\nonumber
\end{align}
These estimates and (\ref{4.83}) and an obvious identity
\begin{equation}\label{4.82}
u_\e-u_0=v_\e + (1-\Xi_\e) u_0
\end{equation}
 prove (\ref{2.21}).

\subsection{$L_2$-estimates}

Here we prove inequalities (\ref{2.17}) and (\ref{2.23}). The former is implied immediately by  identity (\ref{4.82}) and estimates (\ref{3.63}), (\ref{2.15}), (\ref{2.16}).

In the proof of (\ref{2.23}) we follow an approach proposed recently in \cite{PMA22-2}, which is a modification of the technique used in \cite{Pas1}, \cite{Pas2}, \cite{Pas3}, \cite{Sen1}, \cite{Sen2}. Namely, we first
introduce a  differential expression
\begin{equation*}
\cL^*:=-\sum\limits_{i,j=1}^{n} \frac{\p\ }{\p x_i} A_{ij}^\e \frac{\p\ }{\p x_j}  - \sum\limits_{j=1}^{n}\frac{\p\ }{\p x_j}\overline{A_j}  + \overline{A_0}
\end{equation*}
and consider an auxiliary boundary value problem
\begin{equation}\label{5.3}
(\cL^*-\overline{\l})w=f_\e
\quad\text{in}\quad \Om,\qquad w=0\quad\text{on} \quad \p\Om,
\end{equation}
where $f_\e=v_\e$ in $\Om^\e$ and $f_\e=0$ in $\tht^\e$; here we use the notations from Subsection~\ref{sec4.1}. Since $A_j\in W_\infty^1(\Om)$, this problem is of the same nature as (\ref{2.11}). This is why it is solvable for $\RE\l\leqslant \l_0$ and its solution belongs to $W_2^2(\Om)$ and satisfies the estimate
\begin{equation}\label{5.11}
\|w\|_{W_2^2(\Om)}\leqslant C\|v_\e\|_{L_2(\Om^\e)}.
\end{equation}
Hereinafter by $C$ we denote inessential constants independent of $\e$, $k$, $f$, $v_\e$ and $w$.

In what follows the function $v_\e$ is supposed to be continued inside $\tht^\e$ in accordance with (\ref{4.71}), (\ref{4.72}) and thus is  regarded as an element of $W_2^1(\Om)$.
We then write an integral identity associated with problem (\ref{5.3}) choosing $v_\e$ as a test function:
\begin{equation}\label{5.38}
\begin{aligned}
(v_\e,f_\e)_{L_2(\Om)}=\|v_\e\|_{L_2(\Om^\e)}^2 =&\hf(w,v_\e)-\l(w,v_\e)_{L_2(\Om^\e)} +
(\rA \nabla w,\nabla v_\e)_{L_2(\tht^\e)}
\\
 &+ \sum\limits_{j=1}^{n} \left(A_j \frac{\p w}{\p x_j}, v_\e\right)_{L_2(\tht^\e)}
  + (A_0 w,v_\e)_{L_2(\tht^\e)} - \l(w,v_\e)_{L_2(\tht^\e)}.
\end{aligned}
\end{equation}
By   estimates (\ref{4.73}), (\ref{4.74}), (\ref{4.83}), (\ref{5.11})  we then immediately obtain:
\begin{equation}\label{5.39}
\begin{aligned}
\bigg|\sum\limits_{j=1}^{n} \left(A_j \frac{\p w}{\p x_j}, v_\e\right)_{L_2(\tht^\e)}
& + (A_0 w,v_\e)_{L_2(\tht^\e)}- \l(w,v_\e)_{L_2(\tht^\e)}
\bigg|\leqslant C\e\eta\vk^{\frac{1}{2}}  \|w\|_{W_2^2(\Om)}\|v_\e\|_{W_2^1(\Om^\e)}
\\
&\leqslant C  \e\eta\vk^{\frac{1}{2}} \big( \eta^{\frac{n}{2}-1}\e^{-1}\vk^{-\frac{1}{2}}+ \e^\frac{1}{2}\eta^\frac{1}{2}+\e \big) \|f\|_{L_2(\Om)}
\|v_\e\|_{L_2(\Om^\e)}.
\end{aligned}
\end{equation}

By $\Pi_\e$ we denote a  particular case of function $\Xi_\e$ in
the case when on the boundaries of all cavities
the Dirichlet condition is imposed. In other words,   only the functions $X_{k,\e}$ satisfying Dirichlet condition (\ref{2.27}) are used in (\ref{2.31}) for all $k\in\mathds{M}^\e$ while defining $\Pi_\e$. The function $\Pi_\e$ vanishes on $\p\tht^\e$ and is real.

We write identities (\ref{4.12}), (\ref{4.13}), (\ref{4.14}) replacing there $v_\e$ by $\Pi_\e w$ and then we take the difference of the obtained analogues of (\ref{4.12}), (\ref{4.13}). This gives:
\begin{equation}\label{5.40}
\begin{aligned}
\hf(v_\e,\Pi_\e w)-\l(v_\e,\Pi_\e w)_{L_2(\Om^\e)}=& \big(f,(1-\overline{\Xi_\e})\Pi_\e w\big)_{L_2(\Om^\e)} + (\rA\nabla u_0,\Pi_\e w\nabla\overline{\Xi_\e})_{L_2(\Om^\e)}
\\
&- \big(\rA u_0\nabla \Xi_\e,\nabla(\Pi_\e w)\big)_{L_2(\Om^\e)} - \sum\limits_{j=1}^{n} \left(A_j u_0\frac{\p\Xi_\e}{\p x_j},\Pi_\e w\right)_{L_2(\Om^\e)}.
\end{aligned}
\end{equation}
Let us estimate the right hand side of this identity.

Since the function $\Pi_\e$ is a particular case of $\Xi_\e$, it possesses the same properties, namely, relations (\ref{3.66}), (\ref{4.54})  and Lemma~\ref{lm4.3} hold true for $\Pi_\e$. Then by (\ref{3.63}), (\ref{5.11}) we obtain:
\begin{equation}\label{5.41}
\big|\big(f,(1-\overline{\Xi_\e})\Pi_\e w\big)_{L_2(\Om^\e)}\big|\leqslant
 C\|f\|_{L_2(\Om)}\|(1-\Xi_\e) w\|_{L_2(\Om^\e)}
\leqslant C\e\|f\|_{L_2(\Om)} \|v_\e\|_{L_2(\Om^\e)}.
\end{equation}
Using the definition of the functions $\Xi_\e$ and $\Pi_\e$, we integrate by parts as follows:
\begin{align*}
\sum\limits_{j=1}^{n} \left(A_j u_0\frac{\p\Xi_\e}{\p x_j},\Pi_\e w\right)_{L_2(\Om^\e)}=\sum\limits_{j=1}^{n} \left(A_j u_0\frac{\p(\Xi_\e-1)}{\p x_j},\Pi_\e w\right)_{L_2(\Om^\e)}= \sum\limits_{j=1}^{n} \int\limits_{\Om^\e}(\Xi_\e-1)\frac{\p\ }{\p x_j} A_j u_0\Pi_\e \overline{w}\,dx.
\end{align*}
Hence, by estimates (\ref{3.66}), (\ref{3.63}), (\ref{3.64})
for the functions $\Xi_\e$ and $\Pi_\e$ and by estimates (\ref{2.15}), (\ref{5.11}) we obtain:
\begin{equation}\label{5.42}
\begin{aligned}
\bigg|\sum\limits_{j=1}^{n} \left(A_j u_0\frac{\p\Xi_\e}{\p x_j},\Pi_\e w\right)_{L_2(\Om^\e)}\bigg|\leqslant & C \|(\Xi_\e-1)u_0\|_{L_2(\Om^\e)}\big(\|w\|_{W_2^1(\Om^\e)}
+\|w \nabla \Pi_\e\|_{L_2(\Om^\e)} \big)
\\
&+C\|(\Xi_\e-1)w\|_{L_2(\Om^\e)} \|u_0\|_{W_2^1(\Om)}
\\
\leqslant & C (\e^2+\e\eta) \|u_0\|_{W_2^2(\Om)}\|w\|_{W_2^2(\Om)}
\leqslant  C (\e^2+\e\eta)\|f\|_{L_2(\Om)}\|v_\e\|_{L_2(\Om^\e)}.
\end{aligned}
\end{equation}

We rewrite two remaining terms in the right hand side of (\ref{5.40}) as \begin{equation}\label{5.43}
\begin{aligned}
(\rA\nabla u_0,\Pi_\e w\nabla\overline{\Xi_\e})_{L_2(\Om^\e)}
&- \big(\rA u_0\nabla \Xi_\e,\nabla(\Pi_\e w)\big)_{L_2(\Om^\e)}
\\
=&
(\rA\nabla u_0,  w\nabla\overline{\Xi_\e})_{L_2(\Om^\e)}- \big(\rA u_0\nabla \Xi_\e,\nabla w \big)_{L_2(\Om^\e)}
\\
&+
(\rA\nabla u_0,(\Pi_\e-1) w\nabla\overline{\Xi_\e})_{L_2(\Om^\e)}
- \big(\rA u_0\nabla \Xi_\e,(\Pi_\e-1)\nabla w\big)_{L_2(\Om^\e)}
\\
&- \big(\rA u_0\nabla \Xi_\e,w\nabla \Pi_\e \big)_{L_2(\Om^\e)}.
\end{aligned}
\end{equation}
By Lemma~\ref{lm4.3} and  estimates (\ref{2.15}), (\ref{5.11}) we see that
\begin{equation}\label{5.44}
\begin{aligned}
\Big|(\rA\nabla u_0,(\Pi_\e-1) & w\nabla\overline{\Xi_\e})_{L_2(\Om^\e)}
- \big(\rA u_0\nabla \Xi_\e,(\Pi_\e-1)\nabla w\big)_{L_2(\Om^\e)}
- \big(\rA u_0\nabla \Xi_\e,w\nabla \Pi_\e \big)_{L_2(\Om^\e)}
\Big|
\\
\leqslant & C \Big( \|(\Pi_\e-1) \nabla u_0\|_{L_2(\Om^\e)} \|w\nabla \Xi_\e\|_{L_2(\Om^\e)}
\\
&\hphantom{C \Big(}+ \|(\Pi_\e-1) \nabla w\|_{L_2(\Om^\e)} \|u_0\nabla \Xi_\e\|_{L_2(\Om^\e)}
+\|u_0\nabla \Xi_\e\|_{L_2(\Om^\e)}
\|w\nabla \Pi_\e\|_{L_2(\Om^\e)} \Big)
\\
\leqslant & C(\eta^{n-2}\e^{-2}\vk^{-1}+\e^2+\e\eta)
 \|u_0\|_{W_2^2(\Om)} \|w\|_{W_2^2(\Om)}.
\end{aligned}
\end{equation}
In the other two terms in the right hand side of (\ref{5.43}) we integrate by parts using the definition of the functions $\Xi_\e$ and $\Pi_\e$:
\begin{align*}
(\rA\nabla u_0,  w\nabla\overline{\Xi_\e})_{L_2(\Om^\e)}- \big(\rA u_0\nabla \Xi_\e,\nabla w \big)_{L_2(\Om^\e)}=&
(\rA\nabla u_0,  w\nabla\overline{\Xi_\e})_{L_2(\Om)}- \big(\rA u_0\nabla \Xi_\e,\nabla w \big)_{L_2(\Om)}
\\
=&
\int\limits_{\Om} (\Xi_\e-1)\big( \dvr \rA u_0 \nabla \overline{w} -    \dvr \rA  \overline{w}\nabla u_0
\big)\,dx
\\
=&\int\limits_{\Om} (\Xi_\e-1)\big(u_0\dvr \rA  \nabla \overline{w} -    \overline{w}\dvr \rA   \nabla u_0
\big)\,dx.
\end{align*}
By   estimates (\ref{2.15}), (\ref{3.66}), (\ref{3.63}), (\ref{5.11}) and Lemma~\ref{lm3.6} we then obtain:
\begin{align*}
\Big|(\rA\nabla u_0,  w\nabla\overline{\Xi_\e})_{L_2(\Om^\e)}&- \big(\rA u_0\nabla \Xi_\e,\nabla w \big)_{L_2(\Om^\e)}
\Big|
\\
\leqslant & C  \|(\Xi_\e-1)u_0\|_{L_2(\Om)}\|w\|_{W_2^2(\Om)}
 + C
\|(\Xi_\e-1)w\|_{L_2(\Om)}\|u_0\|_{W_2^2(\Om)}
\\
\leqslant &
C \Big(\|(\Xi_\e-1)u_0\|_{L_2(\Om^\e)}+
\|u_0\|_{L_2(\tht^\e)}\Big)
\|w\|_{W_2^2(\Om)}
\\
& +C\Big(
\|(\Xi_\e-1)w\|_{L_2(\Om^\e)}+ \|(\Xi_\e-1)w\|_{L_2(\tht^\e)}
\Big)\|u_0\|_{W_2^2(\Om)}
\\
\leqslant & C(\e^2+\e\eta)
 \|u_0\|_{W_2^2(\Om)}\|w\|_{W_2^2(\Om)}
\\
\leqslant & C  (\e^2+\e\eta)
\|f\|_{L_2(\Om)}\|v_\e\|_{L_2(\Om^\e)}.
\end{align*}
These estimates and (\ref{5.40}), (\ref{5.41}), (\ref{5.42}), (\ref{5.43}), (\ref{5.44}) yield
\begin{equation}\label{5.45}
\big|\hf(v_\e,\Pi_\e w)-\l(v_\e,\Pi_\e w)_{L_2(\Om^\e)}\big|
\leqslant C\big(\eta^{n-2}\e^{-2}\vk^{-1}+\e^2+\e\eta\big) \|f\|_{L_2(\Om)}\|v_\e\|_{L_2(\Om^\e)}.
\end{equation}
It also follows from Lemma~\ref{lm4.3} and estimate (\ref{3.66}) for the function $\Pi_\e$ and from (\ref{4.83}) that
\begin{align*}
\big|\hf(v_\e,(1-\Pi_\e) w)-\l(v_\e,(1-\Pi_\e) w)_{L_2(\Om^\e)}\big|
\leqslant & C \|v_\e\|_{W_2^1(\Om^\e)} \big( \|(1-\Pi_\e) \nabla w\|_{L_2(\Om^\e)} + \|w\nabla \Pi_\e\|_{L_2(\Om^\e)}
\big)
\\
\leqslant & C\big(\eta^{n-2}\e^{-2}\vk^{-1}+\e^2+\e\eta
\big) \|f\|_{L_2(\Om)} \|w\|_{W_2^2(\Om)}
\\
\leqslant & C\big(\eta^{n-2}\e^{-2}\vk^{-1}+\e^2+\e\eta
\big) \|f\|_{L_2(\Om)} \|v_\e\|_{L_2(\Om^\e)}.
\end{align*}
These inequalities and (\ref{5.45}),  (\ref{5.39}) allow us to estimate the left hand side in (\ref{5.38}):
\begin{equation*}
\|v_\e\|_{L_2(\Om^\e)}^2\leqslant C\big(\eta^{n-2}\e^{-2}\vk^{-1} +\e^2+\e\eta \big)\|f\|_{L_2(\Om)} \|v_\e\|_{L_2(\Om^\e)}
\end{equation*}
and hence,
\begin{equation*}
\|v_\e\|_{L_2(\Om^\e)} \leqslant C\big(\eta^{n-2}\e^{-2}\vk^{-1} +\e^2+\e\eta \big)\|f\|_{L_2(\Om)}.
\end{equation*}
Employing now estimate (\ref{3.63}) with $u=u_0$ and identity (\ref{4.82}), we arrive at (\ref{2.23}).

\subsection{Sharpness of estimates}

In this subsection we study the sharpness of the terms in the right hand sides of inequalities (\ref{2.16}), (\ref{2.17}), (\ref{2.21}), (\ref{2.23}) are order sharp. First let us show that the term $\|\b_\e-\b\|_{\mathfrak{M}}$ in (\ref{2.16}), (\ref{2.17}) is order sharp.

We choose  $\cL:=-\D+1$ and we impose only the Dirichlet condition on the boundaries of the cavities, that is, $\mathds{M}_R^\e=\emptyset$. In this case, the function $\b_\e$ is non-negative. For $n=2$ this fact is implied by  definition (\ref{2.9}), while for $n\geqslant 3$ it follows from a simple integration by parts:
\begin{align*}
0=&\lim\limits_{R\to+\infty} \int\limits_{\{\xi:\, |\rA(M_k^\e)\xi|<R,\ \xi\notin\om_{k,\e}\}} X_{k,\e} \dvr_\xi \rA(M_k^\e)\nabla_\xi X_{k,\e}\,d\xi
\\
=&\lim\limits_{R\to+\infty} \int\limits_{\{\xi:\, |\rA(M_k^\e)\xi|=R\}}X_{k,\e}\rA(M_k^\e)\nabla_\xi X_{k,\e}\cdot\nu\,ds-
\int\limits_{\mathds{R}^n\setminus\om_{k,\e}} \rA(M_k^\e)\nabla_\xi X_{k,\e}\cdot\nabla_\xi X_{k,\e}
\,d\xi
\\
=&(2-n) K_{k,\e}\mes_{n-1} \p B_1(0){\det}^\frac{1}{2}\rA(M_k^\e)
-
\int\limits_{\mathds{R}^n\setminus\om_{k,\e}} \rA(M_k^\e)\nabla_\xi X_{k,\e}\cdot\nabla_\xi X_{k,\e}
\,d\xi.
\end{align*}
Definition (\ref{2.3}) of the norm in $\mathfrak{M}$ yields that
\begin{equation*}
\int\limits_{\Om} \b_\e \phi\,dx\to \int\limits_{\Om} \b\phi\,dx,\quad \e\to+0,\quad\text{for all}\quad \phi\in C^\infty_0(\Om).
\end{equation*}
Choosing then non-negative functions $\phi$, we see that $\b\geqslant 0$ almost everywhere in $\Om$.

 We also assume that $\eta$ is so that $\eta^{n-2}\e^{-2}\vk^{-1}=\g>0$ for all $\e$ and the domain $\Om$ is bounded.
  Then  the choice $\l=\l_0=0$   ensures the solvability of both perturbed and limiting problems (\ref{2.7}), (\ref{2.11}) as well as of the following auxiliary boundary problem:
\begin{equation}\label{5.48}
(\cL+\g\b_\e)\tilde{u}_0=f\quad\text{in}\quad\Om,\qquad \tilde{u}_0=0\quad\text{on}\quad\p\Om.
\end{equation}
Since the function $\b_\e$ is piece-wise constant, non-negative and is uniformly bounded due its definition and (\ref{3.50a}), the above problem is solvable in $W_2^2(\Om)$ and its solution also satisfies estimate (\ref{2.15}). Then we can replace the function $u_0$ and problem (\ref{2.11}) by $\tilde{u}_0$ and problem (\ref{5.48}) and reproduce all calculations in  Subsection~\ref{sec4.1}  up to (\ref{4.69}), (\ref{4.68}) taking into consideration that $\mathds{M}_R^\e=\emptyset$.  In the right hand side of identity (\ref{4.69}) then the first term vanishes and this removes the term $\|\b_\e-\b\|_{\mathfrak{M}}$ from (\ref{4.79}), (\ref{4.80}). Using then (\ref{4.1}), (\ref{4.82}),  we get modifications of estimates (\ref{2.16}), (\ref{2.17}):
\begin{equation*}
\|u_\e-\Xi_\e \tilde{u}_0\|_{W_2^1(\Om^\e)}\leqslant C \e \|f\|_{L_2(\Om)},\qquad \|u_\e- \tilde{u}_0\|_{L_2(\Om^\e)}\leqslant C \e \|f\|_{L_2(\Om)}.
\end{equation*}
Therefore, it is sufficient to prove that
\begin{equation}\label{5.50}
 \|\tilde{u}_0-u_0\|_{L_2(\Om^\e)}\geqslant C
\|\b_\e-\b\|_{\mathfrak{M}}   \|f\|_{L_2(\Om)},
\end{equation}
for  some $f$ to show the sharpness of the term $\|\b_\e-\b\|_{\mathfrak{M}}$ in (\ref{2.16}), (\ref{2.17}). If $\|\b_\e-\b\|_{\mathfrak{M}}=0$, then the above inequalities are obvious and this is why in what follows we assume that $\|\b_\e-\b\|_{\mathfrak{M}}\ne 0$.

Assume that $\|\b_\e-\b\|_{\mathfrak{M}}\ne0$. We rewrite definition   (\ref{2.3}) of the norm $\|\b_\e-\b\|_{\mathfrak{M}}$ as
\begin{equation}
\|\b_\e-\b\|_{\mathfrak{M}}=\sup\limits_{u\in \Ho_2^1(\Om)\cap W_2^2(\Om)}\frac{1}{\|u\|_{W_2^2(\Om)}}
\sup\limits_{
v\in \Ho_2^1(\Om)}  \frac{|((\b_\e-\b) u,v)_{L_2(\Om)}|
}{\|v\|_{W_2^1(\Om)}}
\end{equation}
and conclude that there exists a non-zero function $u^\e\in W_2^2(\Om)\cap \Ho_2^1(\Om)$  such that
\begin{equation}\label{5.51}
 \big|\big((\b_\e-\b)u^\e,v\big)_{L_2(\Om)}\big|
 \geqslant
 \frac{1}{2} \|\b_\e-\b\|_{\mathfrak{M}} \|u^\e\|_{W_2^2(\Om)}\|v\|_{W_2^1(\Om)}\quad \text{for all}
 \quad v\in\Ho_2^1(\Om).
\end{equation}
We choose $f:=(\cL+\g\b_\e)u^\e$ and we see that $\tilde{u}_0=u^\e$ solves problem (\ref{5.48}) and
\begin{equation}\label{5.12}
\|f\|_{L_2(\Om)}\leqslant C\|u^\e\|_{W_2^2(\Om)}
\end{equation}
with some fixed constant $C$ independent of $\e$ and $u^\e$.  Using the corresponding solution $u_0$ of problem (\ref{2.11}) with $\l=0$, we define $\phi^\e:=u^\e-u_0$. The latter function solves the  boundary value problem
\begin{equation}\label{5.4}
(\cL+\g \b) \phi^\e+\g (\b_\e-\b)u^\e=0\quad\text{in}\quad\Om,\qquad \phi^\e=0\quad\text{on}\quad\p\Om.
\end{equation}
Let $\L$ and $\Phi$ be a positive eigenvalue and an associated normalized in $L_2(\Om)$ eigenfunction of a self-adjoint operator in $L_2(\Om)$ with the differential expression $\cL+\g\b$ and the Dirichlet condition on $\p\Om$. We write integral identity corresponding to (\ref{5.4}) with $\Phi$ as a test function to obtain:
\begin{equation*}
(\phi^\e,\Phi)_{L_2(\Om)}=-\frac{1}{\L}\big((\b_\e-\b) u^\e,\Phi\big)_{L_2(\Om)}.
\end{equation*}
Hence, by (\ref{5.51}), (\ref{5.12})
\begin{align*}
\|\phi^\e\|_{L_2(\Om)}\geqslant & \big|(\phi^\e,\Phi)_{L_2(\Om)}\big|\|\Phi\|_{L_2(\Om)}=\frac{1}{\L}\big|\big((\b_\e-\b) u^\e,\Phi\big)_{L_2(\Om)}|
\geqslant
 \frac{1}{2\L} \|\b_\e-\b\|_{\mathfrak{M}} \|u^\e\|_{W_2^2(\Om)}\|\Phi\|_{W_2^1(\Om)}
 \\
 \geqslant & \frac{1}{2\L} \|\b_\e-\b\|_{\mathfrak{M}} \|u^\e\|_{W_2^2(\Om)}\|\Phi\|_{L_2(\Om)}\geqslant
 \frac{1}{2\L} \|\b_\e-\b\|_{\mathfrak{M}} \|f\|_{L_2(\Om)}
\end{align*}
and this proves (\ref{5.50}).

In order to prove the sharpness of the term $\big|\eta^{n-2}\e^{-2}\vk^{-1}-\g\big|$ in estimates (\ref{2.16}), (\ref{2.17}), we proceed in a similar way. Namely, assuming that $\eta^{n-2}\e^{-2}\vk^{-1}$ does not coincide with its limit $\g$, now we define $\tilde{u}_0$ as a solution to the problem
\begin{equation}\label{5.68}
(-\D+\mu_3\b_0)\tilde{u}_0=f\quad\text{in}\quad\Om,\qquad \tilde{u}_0=0\quad\text{on}\quad\p\Om
\end{equation}
with $\mu_3:=\eta^{n-2}\e^{-2}\vk^{-1}$. Then we can again reproduce the calculations from Section~\ref{sec4.1} skipping just identities (\ref{4.65}) and replacing $\g$ by
$\eta^{n-2}\e^{-2}\vk^{-1}$ in all relations after (\ref{4.65}). This gives the estimates
\begin{align*}
&\|u_\e-\Xi_\e \tilde{u}_0\|_{W_2^1(\Om^\e)}\leqslant C\big(\e
+
\|\b_\e-\b\|_{\mathfrak{M}}
+(\e\eta\mu_1\vk)^{-1}  + \e\eta\vk\mu_2
\big)\|f\|_{L_2(\Om)},
\\
&  \|u_\e- \tilde{u}_0\|_{L_2(\Om^\e)}\leqslant C\big(\e
+
\|\b_\e-\b\|_{\mathfrak{M}}
+(\e\eta\mu_1\vk)^{-1} + \e\eta\vk\mu_2
\big)\|f\|_{L_2(\Om)}.
\end{align*}
At the same time, it is easy to see that the solution  of problem (\ref{5.68}) is analytic in $\mu_3$ and for $\mu_3=\g$, this solution coincides with the solution $u_0$ to homogenized problem (\ref{2.11}). Hence, in the general situation, the next-to-leading term in the Taylor expansion of $\tilde{u}_0$ in $\mu_3-\g$ is non-zero and the estimates
\begin{equation*}
\|u_0-\tilde{u}_0\|_{L_2(\Om)}\leqslant C|\mu_3-\g|\|f\|_{L_2(\Om)},\qquad
\|u_0-\tilde{u}_0\|_{W_2^1(\Om)}\leqslant C|\mu_3-\g|\|f\|_{L_2(\Om)},
\end{equation*}
are order sharp. In particular, we can calculate the norms in their left hand sides over the set $\{x\in\Om:\, \Xi_\e(x)=1\}$ and this proves that the term $\big|\eta^{n-2}\e^{-2}\vk^{-1}-\g\big|$ in estimates (\ref{2.16}), (\ref{2.17}) is order sharp.

We proceed to checking the sharpness of the other terms in (\ref{2.16}) and (\ref{2.21}). We are going to do this by adducing an appropriate example.  We let $\square:=[-2,2)^n$, $\om:=B_1(0)$  and   $\Om=\mathds{R}^n$. The points $M_k^\e$ are defined as $M_k^\e:=\e k$, $k\in\mathds{M}^\e:=4\mathds{Z}^n$. In this case, we deal with a periodic perforation in $\Om$. Each cavity is ball of the radius $\e\eta$ centered at a point $\e k$, $k\in\mathds{Z}^n$, and hence,
\begin{equation*}
\tht^\e=\bigcup\limits_{k\in\mathds{Z}^n} B_{\e\eta}(k),\qquad \Om^\e=\mathds{R}^n\setminus \tht^\e.
\end{equation*}
The differential expression is chosen to be the negative Laplacian, $\cL:=-\D+1$. It is clear that  we can take $\l=\l_0=0$. We choose $R_1=1$, $R_2:=\frac{7}{6}$, $R_3:=\frac{3}{2}$,  $R_4:=\frac{4}{3}$.

We consider the solution to the problem
\begin{equation*}
(\cL+ \eta^{n-2}\e^{-2}\vk^{-1}\b_\e) \tilde{u}_0=f\quad\text{in}\quad\mathds{R}^d
\end{equation*}
and hence, as in the first part of section, by reproducing the arguing from Section~\ref{sec4.1}, the solution $u_\e$ to the corresponding  equation in (\ref{2.7}) satisfies the modified versions of estimates (\ref{2.16}), (\ref{2.17}):
\begin{equation}\label{5.74}
\begin{aligned}
&\|u_\e-\tilde{u}_0\Xi_\e\|_{W_2^1(\mathds{R}^n\setminus\tht^\e)}  \leqslant  C\big(\e +
 (\e\eta\mu_1\vk)^{-\frac{1}{2}}  + \e\eta\vk\mu_2\big)\|f\|_{L_2(\Om)},
 \\
 &\|u_\e-\tilde{u}_0\|_{L_2(\mathds{R}^n\setminus\tht^\e)}  \leqslant  C \big(\e +
 (\e\eta\mu_1\vk)^{-\frac{1}{2}} + \e\eta\vk\mu_2\big)
\|f\|_{L_2(\Om)}.
\end{aligned}
\end{equation}
This is why, to confirm the sharpness of the other terms in (\ref{2.16}), (\ref{2.17}), we need to estimate from below the norms in the left hand sides of the above inequalities.

We first consider the case of only Dirichlet conditions on the boundaries of the cavities, that is,
  $\mathds{M}_R^\e=\emptyset$ and $\mathds{M}_D^\e=\mathds{M}^\e$. Such choice of the boundary conditions on $\p\tht^\e$ removes the terms $ (\e\eta\mu_1\vk)^{-\frac{1}{2}}$ and $ \e\eta\vk\mu_2$ from inequalities (\ref{5.74}).
The functions $X_{k,\e}$ and $Z_{k,\e}$ for the considered model can be found explicitly
\begin{equation}\label{5.24}
\begin{aligned}
&
Z_{k,\e}(\xi)=X_{k,\e}(\xi)=\left\{
\begin{aligned}
&1-|\xi|^{-n+2},\quad && n\geqslant 3,
\\
&\ln|\xi|, && n=2,
\end{aligned}\right.
\end{aligned}
\end{equation}
The corresponding function $\b_\e$ given by (\ref{2.9}) then is $\e\square$-periodic and reads as $\b_\e(x)=\b_0(x\e^{-1})$, where
\begin{equation}\label{5.81}
\begin{aligned}
&\b_0(\z):=\frac{(n-2)}{R_3^n \mes_n B_1(0)} \quad&& \text{on}\quad B_{R_3}(k),\quad k\in 4\mathds{Z}^n,\quad n\geqslant 3,
\\
&\b_0(\z):=\frac{1}{R_3^2 \mes_n B_1(0)} &&  \text{on}\quad B_{R_3}(k),\quad k\in 4\mathds{Z}^n, \quad n=2,
\\
& \b_0(\z):=0 &&\text{on}\quad \mathds{R}^n\setminus \bigcup\limits_{k\in 4\mathds{Z}^n}
B_{R_3}(k).
\end{aligned}
\end{equation}

By $Y_0=Y_0(\z)$ we denote the $\square$-periodic solution to the following boundary value problem
\begin{equation}
\label{5.25}
\begin{gathered}
-\D_\z Y_0=\b_0\quad\text{in}\quad \mathds{R}^n\setminus 4\mathds{Z}^n,
\\
\begin{aligned}
&Y_0(\z)=-\frac{1}{|\z-k|^{n-2}}+O(|\z-k|), \quad&& \z\to k,\quad n\geqslant 3,
\\
& Y_0(\z)=\ln|\z-k|+O(|\z-k|), &&\z\to k, \quad n=2.
\end{aligned}
\end{gathered}
\end{equation}
It is easy to see that  this problem satisfies the standard solvability condition; the uniqueness of the solution  is ensured by the order of the error terms in the prescribed asymptotics.

Given an arbitrary infinitely differentiable function $\tilde{u}_0=\tilde{u}_0(x)$,
we denote
\begin{equation*}
 f:=(\cL+ \eta^{n-2}\e^{-2}\vk^{-1}\b_\e) \tilde{u}_0\in C_0^\infty(\mathds{R}^d),
\qquad
U_\e(x):=
\left(1+\eta^{n-2}(\e)\vk^{-1}(\e) Y_0\left(\frac{x}{\e}\right)\right) \tilde{u}_0(x).
\end{equation*}
It is clear that
\begin{equation}\label{5.1}
\|f\|_{L_2(\mathds{R}^d)}\leqslant C\|\tilde{u}_0\|_{W_2^2(\mathds{R}^d)},
\end{equation}
where a constant $C$ is independent of $\e$, $\eta$ and $\tilde{u}_0$.
It also follows from (\ref{5.25})
that the function $U_\e$ solves the boundary value problem
\begin{align}\label{5.64}
 &\cL U_\e=f+f_\e\quad\text{in}\quad\mathds{R}^n\setminus\tht^\e, \qquad U_\e=\eta^{n-2}\vk^{-1}\vp_D^\e \quad\text{on}\quad \tht^\e,
\\
&\nonumber
\begin{aligned}
&f_\e(x):=\eta^{n-2}(\e)\vk^{-1}(\e)Y_0\left(\frac{x}{\e}\right)\cL\tilde{u}_0 -2\eta^{n-2}(\e)\vk^{-1}(\e) \nabla Y_0\left(\frac{x}{\e}\right) \cdot \nabla \tilde{u}_0(x),
\\
&\vp_D^\e(x):=\left(Y_0\left(\frac{x}{\e}\right)+\eta^{-n+2}(\e)\vk(\e) \right)
\tilde{u}_0(x).
\end{aligned}
\end{align}
By straightforward calculations we confirm that
\begin{align*}
&f_\e(x)=\eta^{n-2}(\e)\vk^{-1}(\e) f_{\e,0}(x)-2 \eta^{n-2}(\e)\vk^{-1}(\e) \sum\limits_{j=1}^{n} \frac{\p f_{\e,j}}{\p x_j}(x),
\\
&f_{\e,0}(x):=Y_0\left(\frac{x}{\e}\right)(-\cL+2)\tilde{u}_0(x),\qquad
f_{\e,j}(x):=Y_0\left(\frac{x}{\e}\right) \frac{\p\tilde{u}_0}{\p x_j}(x).
\end{align*}

By $\chi_3=\chi_3(\z,\eta)$ we denote an infinitely differentiable $\square$-periodic cut-off function equalling to one as $|\z-k|\leqslant 2\eta$ and vanishing as $|\z-k|\geqslant 3\eta$ for $k\in 4\mathds{Z}^n$ and obeying the uniform estimate $|\nabla_\z \chi_3(\z,\eta)|\leqslant C\eta^{-1}$ with a constant $C$ independent of $\z$ and $\eta$. The function
\begin{align*}
&\tilde{\vp}_D^\e(x):= \tilde{u}_0(x) \chi_3\left(\frac{x}{\e},\eta\right)  Y_1\left(\frac{x}{\e}\right)\quad \text{on}\quad B_{3\e\eta}(k),
\qquad \tilde{\vp}_D^\e(x):=0\quad \text{outside}\quad B_{3\e\eta}(k),\qquad k\in 4\mathds{Z}^n,
\\
&Y_1(\xi):=Y_0(\xi)+\left\{
\begin{aligned}
&|\xi-k|^{-n+2} &&\text{as}\quad n\geqslant 3,\qquad k\in 4\mathds{Z}^n,
\\
-&\ln|\xi-k| &&\text{as}\quad n=2,\qquad k\in 4\mathds{Z}^n,
\end{aligned}\right.
\end{align*}
multiplied by $\eta^{n-2}\vk^{-1}$ satisfies the boundary condition in (\ref{5.64}). Then we consider the solution  $u_\e$ to problem (\ref{2.7}) with the introduced function $f$ and in a standard way we get the estimate:
\begin{equation}\label{5.67}
\begin{aligned}
\|U_\e-u_\e\|_{W_2^1(\mathds{R}^n\setminus\tht^\e)}\leqslant C\bigg(&\eta^{n-2}\vk^{-1} \|\tilde{\vp}_D^\e\|_{W_2^1(\mathds{R}^n\setminus\tht^\e)} + \eta^{n-2}\vk^{-1} \|f_{\e,0}\|_{L_2(\mathds{R}^n\setminus\tht^\e)}
\\
&+ \eta^{n-2}\vk^{-1}\e
\sum\limits_{j=1}^{n}\|f_{\e,j}\|_{L_2(\mathds{R}^n\setminus\tht^\e)}
\bigg),
\end{aligned}
\end{equation}
where $C$ is a fixed constant independent of $\e$, $\eta$, $f_\e$, $\tilde{\vp}_D^\e$.

In view of convergence (\ref{2.14}) and the asymptotics for $Y_0$
in (\ref{5.25})
and the smoothness of this function, by routine straightforward calculations we find that
\begin{equation}\label{5.69}
\eta^{n-2}\vk^{-1}\|\tilde{\vp}_D^\e\|_{W_2^1(\mathds{R}^n\setminus\tht^\e)}\leqslant C\e^{-1}\eta^{\frac{3n}{2}-2}\vk^{-1}\leqslant C \e^2\eta\vk,
\end{equation}
where $C$ are  some constants  independent of $\e$ and $\eta$ but depending on the choice of the function $\tilde{u}_0$. In the same way we find that
\begin{equation*}
\sum\limits_{j=0}^{n}\|f_{\e,j}\|_{L_2(\mathds{R}^n\setminus\tht^\e)}\leqslant C \left\{
\begin{aligned}
& 1+\eta^{2-\frac{n}{2}}, && n\ne 4,
\\
& 1+|\ln\eta|^{\frac{1}{2}}, && n=4,
\end{aligned}\right.
\end{equation*}
where $C$ is a constant independent of $\e$ and $\eta$ but depending on the choice of the function $\tilde{u}_0$.  We substitute these estimates and (\ref{5.69}) into (\ref{5.67}) and use convergence (\ref{2.14}) to obtain:
\begin{equation}\label{5.71}
\|U_\e-u_\e\|_{W_2^1(\mathds{R}^n\setminus\tht^\e)}
\leqslant C \e^2+C\e\eta \left\{
\begin{aligned}
&1, && n\ne 4,
\\
& |\ln\eta|^{\frac{1}{2}}+1, && n=4,
\end{aligned}\right.
\end{equation}
where $C$ is some constant independent of $\e$ but depending on $\tilde{u}_0$.

The function $\Xi_\e$ defined by (\ref{2.31}) with the functions $Z_{k,\e}$ from (\ref{5.24}) reads as
\begin{equation*}
\Xi_\e(x):=\left\{
\begin{aligned}
& \frac{1-\left|\frac{x}{\e}-k\right|^{-n+2}\eta^{n-2}} {1-\left(\frac{4\eta}{3}\right)^{n-2}} &&\text{in}\quad B_{\frac{3}{4}\e}(\e k)\setminus B_{\e\eta}(\e k),\quad k\in 4\mathds{Z}^n,
\\
& \hphantom{|\rA^{-1}|}1 &&\text{in}\quad \mathds{R}^n\setminus \bigcup\limits_{k\in\mathds{M}^\e} B_{\frac{3}{4}\e}(\e k),
\end{aligned}
\right.
\end{equation*}
as $n\geqslant 3$, and
\begin{equation*}
\Xi_\e(x):=\left\{
\begin{aligned}
& \frac{\ln \left|\frac{x}{\e}-k\right|-\ln\eta} {\ln \frac{3\eta}{4}} &&\text{in}\quad B_{\frac{3}{4}\e}(\e k)\setminus B_{\e\eta}(\e k),\quad k\in 4\mathds{Z}^n,
\\
& \hphantom{|1-|}1 &&\text{in}\quad \mathds{R}^n\setminus \bigcup\limits_{k\in\mathds{M}^\e} B_{\frac{3}{4}\e}(\e k),
\end{aligned}
\right.
\end{equation*}
as $n=2$. As above, by straightforward calculations we confirm that
\begin{align}
&\|\nabla(U_\e-\Xi_\e\tilde{u}_0)\|_{L_2(\mathds{R}^n\setminus\tht^\e)} \geqslant C \eta^{n-2}\e^{-1}\vk^{-1},\label{5.75}
\\
&\|\nabla(U_\e-\tilde{u}_0)\|_{L_2(\mathds{R}^n\setminus\tht^\e)} \geqslant C
\eta^{\frac{n}{2}-1}\e^{-1}\vk^{-\frac{1}{2}},\label{5.76}
\end{align}
provided $|u_0|$ is uniformly separated from zero on some fixed ball. Assuming that $\eta^{n-2}\e^{-2}\vk^{-1}\equiv \g$,  by (\ref{5.1}), (\ref{5.71}), (\ref{5.75}) we get
\begin{equation*}
\|\nabla(U_\e-\Xi_\e\tilde{u}_0)\|_{L_2(\mathds{R}^n\setminus\tht^\e)} \geqslant C \eta^{n-2}\e^{-1}\vk^{-1} \geqslant C\e \|f\|_{L_2(\mathds{R}^d)}
\end{equation*}
with a constant $C$ independent of $\e$ and this proves the sharpness of the term $\e$ in the right hand side of (\ref{2.16}).

As $\g=0$, the solution $u_0$ to the equation
\begin{equation*}
\cL u_0=f\quad\text{in}\quad\mathds{R}^d
\end{equation*}
obviously satisfies the estimate
\begin{equation*}
\|\tilde{u}_0-u_0\|_{W_2^1(\mathds{R}^d)}\leqslant C \eta^{n-2}\e^{-2}\vk^{-1}
\end{equation*}
and in view of   (\ref{5.71}), (\ref{5.76}),    (\ref{2.14}) we also see that the term $\eta^{\frac{n}{2}-1}\e^{-1}\vk^{-\frac{1}{2}}$ in the right hand side of (\ref{2.21}) is order sharp.

Let us show that the term $\e\eta\vk\mu_2$ in (\ref{2.16}) is order sharp. Here we again consider the above example, but on the boundaries of the cavities we impose the Robin condition
\begin{equation*}
\frac{\p u_\e}{\p\nu}+(\e^{-1}\eta^{-1}+\mu_2) u_\e=0\quad\text{on}\quad\p \tht^\e,\qquad a^\e(x,u)=(\e^{-1}\eta^{-1}+\mu_2) u,\qquad b_k=1.
\end{equation*}
Such choice of the boundary conditions removes the term $ (\e\eta\mu_1\vk)^{-\frac{1}{2}}$ from (\ref{5.74}). The functions $X_{k,\e}$ and $Z_{k,\e}$ can be again found explicitly
\begin{align*}
&
Z_{k,\e}(\xi)=X_{k,\e}(\xi)=\left\{
\begin{aligned}
&1-\frac{|\xi|^{-n+2}}{n-1},\quad && n\geqslant 3,
\\
&\ln|\xi|+1, && n=2.
\end{aligned}\right.
\end{align*}
The corresponding function $\b_\e$   then again reads as $\b_\e(x)=\b_0(x\e^{-1})$, where the function $\b_0$ is defined by the formula
\begin{equation*}
\b_0(\z):=\frac{(n-2)}{(n-1)R_3^n \mes_n B_1(0)}\quad \text{on}\quad B_{R_3}(k),\quad k\in 4\mathds{Z}^n,\quad n\geqslant 3,
\end{equation*}
and by the second and third formulae in (\ref{5.81}). We suppose that $\eta^{n-2}\e^{-2}\vk^{-1}\equiv 1$ for all $\e$.

By $u_\e$ we define the solution to problem (\ref{2.7}) for an arbitrary $\mu_2$, while $u_\e^{(0)}$ is the solution to the same problem for $\mu_2=0$. Both these solutions converge to the same solution $u_0$ of the homogenized problem (\ref{2.11}). The function $u_\e$ satisfies estimate (\ref{2.16}), while the function $u_\e^{(0)}$ satisfies the same estimate  with $\mu_2=0$. This is why, in order to prove that the term $\e\eta\vk\mu_2$ is order sharp in estimate (\ref{2.16}), it is sufficient to find an example of the function $f$ such that
\begin{equation}\label{5.89}
\|u_\e-u_\e^{(0)}\|_{L_2(\mathds{R}^d\setminus\tht^\e)} \geqslant C \e\eta\mu_2 \|f\|_{L_2(\mathds{R}^d)}.
\end{equation}

We choose an arbitrary infinitely differentiable compactly supported in $\mathds{R}^n$ function $u_0=u_0(x)$ and we choose the aforementioned function $f$ as $f:=(\cL+\b)u_0$. Then we define:
\begin{align}
&U_\e(x):=u_\e^{(0)}(x) +  \e \eta^{n-1}(\e)\mu_2(\e) Y_2\left(\frac{x}{\e}\right) u_0(x),\label{5.66}
\\
&Y_2:=\frac{1}{n-1}Y_0,\quad n\geqslant 3,\qquad Y_2:=Y_0-\ln\eta,\quad n=2.\nonumber
\end{align}
This function solves the following boundary value problem:
\begin{equation*}
\cL U_\e=f +  f_\e \quad \text{in}\quad \mathds{R}^n\setminus\tht^\e,
\qquad \frac{\p U_\e}{\p \nu}+(\e^{-1}\eta^{-1}+\mu_2)U_\e=\vp_R^\e\quad \text{on}\quad\p\tht^\e,
\end{equation*}
where
\begin{align*}
&f_\e(x)=f_{\e,0}(x) - 2\sum\limits_{j=1}^{n} \frac{\p f_{\e,j}}{\p x_j}(x),
\\
&f_{\e,0}(x):=\e \eta^{n-1} \mu_2 Y_2\left(\frac{x}{\e}\right)(-\cL+2)u_0(x),\qquad
f_{\e,j}(x):=\e \eta^{n-1} \mu_2 Y_2\left(\frac{x}{\e}\right) \frac{\p u_0}{\p x_j}(x),
\\
&\vp_R^\e(x):=\e\eta^{n-1}Y_2 \left(\frac{x}{\e}\right) \left(\frac{\p\  }{\p\nu}  +\mu_2(\e)\right)u_0(x)+\e\eta^{n-1} u_0(x)\left(\frac{\p\  }{\p\nu}
 +\e^{-1}\eta^{-1}\right)Y_3\left(\frac{x}{\e}\right)
 \\
&\hphantom{\vp_R^\e(x):=} +\mu_2(\e)\big(u_\e^{(0)}(x)-u_0(x)\big),
\\
&Y_3(\z):=\frac{1}{n-1}\big(Y_0(\z)+\eta^{-n+2}\big),\quad n\geqslant 3,\qquad
Y_3(\z):=Y_2(\z),\quad n=2.
\end{align*}
Writing then a problem for $U_\e-u_0$ and an associated integral identity with $U_\e-u_0$ as the test function and using (\ref{3.22}), we obtain an analogue of inequality (\ref{5.67}):
\begin{equation}\label{5.88}
\|U_\e-u_\e\|_{W_2^1(\mathds{R}^n\setminus\tht^\e)}\leqslant C
\bigg(\e^{\frac{1}{2}}\eta^{\frac{1}{2}}\vk^\frac{1}{2} \|\vp_R^\e\|_{L_2(\p\tht^\e)} +  \sum\limits_{j=0}^{n} \|f_{\e,j}\|_{L_2(\mathds{R}^n\setminus\tht^\e)}\bigg),
\end{equation}
where $C$ is some constant independent of $\e$, $\vp_R^\e$ and $f_{\e,j}$, $j=0,\ldots,n$. Using (\ref{3.22}), (\ref{2.15}), the asymptotics for $Y_0$ in (\ref{5.25}) and the definition of the function $Y_3$, we find:
\begin{equation}\label{5.82}
\begin{aligned}
\|\vp_R^\e\|_{L_2(\p\tht^\e)}\leqslant &C \big((\e\eta\vk\mu_2+\eta^{n-1})\|u_0\|_{L_2(\p\tht^\e)} + \e\eta\vk\|\nabla u_0\|_{L_2(\p\tht^\e)}+\mu_2\|u_\e^{(0)}-u_0\|_{L_2(\p\tht^\e)}\big)
\\
\leqslant & C (\e\eta\vk)^{\frac{1}{2}}\big((\e\eta\vk\mu_2+\e\eta\vk
)\|f\|_{L_2(\mathds{R}^n)}+\mu_2 \|u_\e^{(0)}-u_0\|_{W_2^1(\mathds{R}^d\setminus\tht^\e)}\big),
\end{aligned}
\end{equation}
where $C$ is some constant independent of $\e$, $u_0$, $u_\e^{(0)}$. The functions $f_{\e,j}$ can be estimated as follows:
\begin{equation*}
\|f_{\e,j}\|_{L_2(\mathds{R}^n\setminus\tht^\e)}\leqslant C\e\eta\mu_2
(\eta^{n-2}+\vr) \|u_0\|_{C^2(\supp u_0)},\quad j=0,\ldots,n,
\end{equation*}
where $C$ is some constant independent of $\e$ and $u_0$. This estimate and (\ref{5.82}), (\ref{5.88}) yield:
\begin{equation}\label{5.65}
\|U_\e-u_\e\|_{W_2^1(\mathds{R}^n\setminus\tht^\e)}\leqslant C
\big( \big(\e \eta \vk \mu_2(\e\eta\vk+\eta^{n-2}+\vr)+\e^2\eta^2\vk^2
\big)\|u_0\|_{C^2(\supp u_0)} +  \e \eta \vk \mu_2 \|u_\e^{(0)}-u_0\|_{W_2^1(\mathds{R}^d\setminus\tht^\e)}\big),
\end{equation}
where $C$ is some constant independent of $\e$, $u_0$, $u_\e^{(0)}$.
At the same time, it follows from definition (\ref{5.66}) of $U_\e$ that provided $|u_0|$ is uniformly separated from zero on some fixed ball, the estimate
\begin{equation*}
\|U_\e-u_\e^{(0)}\|_{W_2^1(\mathds{R}^d\setminus\tht^\e)}\geqslant
C \eta^{\frac{n}{2}}\vk^{\frac{1}{2}}\mu_2=C\e\eta\vk\mu_2
\end{equation*}
holds with some fixed constant $C$ independent of $\e$, where we have also assumed that $\eta^{\frac{n}{2}-1}\e^{-1}\vk^{-\frac{1}{2}}=1$.
This estimate, (\ref{5.65}) and (\ref{2.16}) for $\|u_\e^{(0)}-u_0\|_{W_2^1(\mathds{R}^d\setminus\tht^\e)}$ prove (\ref{5.89}) and hence, the term $\e\eta\vk\mu_2$ is order sharp in (\ref{2.16}).

\section{Convergence in $\|\,\cdot\,\|_{\mathfrak{M}}$-norm}\label{secAS7}

In this section we discuss the convergence postulated in Assumption~\ref{A7}. As a main tool of checking Assumption~\ref{A7}, we propose the following way. We introduce one more space of multipliers $\tilde{\mathfrak{M}}$, which consists of the functions $F$ defined on $\Om$ such that for each $u\in \Ho_2^1(\Om)$ the function $F u$ is a continuous antilinear functional on $\Ho_2^1(\Om)$.  The norm in $\tilde{\mathfrak{M}}$ is introduced as
\begin{equation*}
\|F\|_{\tilde{\mathfrak{M}}}=\sup\limits_{u,v\in \Ho_2^1(\Om)}  \frac{|\la F u,v \ra|
}{\|u\|_{W_2^1(\Om)}\|v\|_{W_2^1(\Om)}}.
\end{equation*}
It is clear that $\tilde{\mathfrak{M}}\subset\mathfrak{M}$ and
\begin{equation}\label{6.2}
\|F\|_{\mathfrak{M}}\leqslant \|F\|_{\tilde{\mathfrak{M}}}.
\end{equation}
Having this inequality in mind, instead of convergence in the space $\mathfrak{M}$ as it is postulated in Assumption~\ref{A7}, we propose to check  the convergence in the space $\tilde{\mathfrak{M}}$.

The convergence in the sense of the norm $\|\,\cdot\,\|_{\tilde{\mathfrak{M}}}$-norm was studied in details in \cite{B22} and a simple criterion was established. Namely, we choose an arbitrary lattice $\G$ in $\mathds{R}^n$ with a periodicity cell $\square$. Given a function $\rho_1=\rho_1(\e)$, we denote
\begin{equation}\label{6.29}
\G_{\rho_1}:=\big\{z\in\G:\ \rho_1 z+\rho_1\square\subset\Om\big\}.
\end{equation}
The mentioned criterion reads as follows: the function $\b_\e$ converges to some function $\b$ in $\|\,\cdot\,\|_{\tilde{\mathfrak{M}}}$-norm if and only if there exist functions $\rho_1=\rho_1(\e)$, $\rho_2=\rho_2(\e)$ such that
\begin{equation}\label{6.24}
\sup\limits_{z\in\G_{\rho_1(\e)}}
\frac{1}{\rho_1^n(\e)} \Bigg|\int\limits_{\rho_1(\e) z+\rho_1(\e)\square}
\big(\b_\e(x)-\b(x)\big)\,dx\Bigg|\leqslant \rho_2(\e),\qquad \rho_1(\e)\to0,\qquad \rho_2(\e)\to0,\qquad \e\to+0.
\end{equation}
Therefore, Assumption~\ref{A7} can be guaranteed by condition (\ref{6.24}), which is very explicit. If condition (\ref{6.24}) is satisfied, by Theorem~2.4  from \cite{B22} we obtain the estimate
\begin{equation}\label{6.25}
\|\b_\e-\b\|_{\tilde{\mathfrak{M}}}\leqslant C (\rho_2+\rho_1);
\end{equation}
hereinafter in this section by $C$ we denote various constants independent of $\e$ and spatial variables.

In paper \cite{B22}, a way for explicit calculation of function $\b$ for a given $\b_\e$ was provided. Namely, let $\om\subset \mathds{R}^n$ be a fixed domain and assume that the function
\begin{equation*}
\b(x):=\lim\limits_{\e\to+0} \frac{1}{\rho_3(\e)\mes_n \om} \int\limits_{x+\rho_3(\e)\om} \b_\e(y)\,dy
\end{equation*}
is well-defined in $\Om$. If the limit in the above formula is uniform in $x$, namely,
\begin{equation*}
\sup\limits_{x}\Bigg|\b(x)-\frac{1}{\rho_3(\e)\mes_n \om} \int\limits_{x+\rho_3(\e)\om} \b_\e(y)\,dy\Bigg|\leqslant \rho_4(\e),\qquad \rho_4(\e)\to+0,\qquad \e\to+0,
\end{equation*}
where $\rho_3(\e)$, $\rho_4(\e)$ are some functions independent  of $x$ and the supremum is taken over $x\in\Om$ such $x+\rho_3(\e)\om\subset \Om$, then condition (\ref{6.24}) is satisfied and
\begin{equation*}
\|\b_\e-\b\|_{\tilde{\mathfrak{M}}}\leqslant C\big(\rho_4+\rho_3^\frac{1}{2}\big).
\end{equation*}

Paper \cite{B22} provides many particular examples of possible functions obeying condition (\ref{6.24}) and all of them can be adapted  also for our particular function $\b_\e$. We do not reproduce here all these examples but instead we discuss a few close examples.

The first example is a sparsely distributed perforation. Here we assume that there exists a function $\rho_5=\rho_5(\e)$ such that
\begin{equation}\label{6.1}
\e\rho_5^{-1}(\e)\to+0,\quad \e\to+0,\qquad B_{\rho_5}(M_k^\e)\cap B_{\rho_5}(M_j^\e)=\emptyset,\quad k\ne j.
\end{equation}
Then according to the example discussed in Section~3.2 in \cite{B22}, condition (\ref{6.24})  is satisfied with $\b=0$ and
\begin{equation*}
\|\b_\e\|_{\tilde{\mathfrak{M}}}\leqslant C\big(\e^n\rho_5^{-n} +\rho_5^{\frac{1}{2}}\big).
\end{equation*}
This estimate can be even improved for our particular case
as the following lemma shows.

\begin{lemma}\label{lm5.1}
Suppose that condition (\ref{6.1}) holds. Then
\begin{equation*}
\|\b_\e\|_{\tilde{\mathfrak{M}}}\leqslant C(\e^n\rho_5^{-n} +\e^2),\quad n\geqslant 3,\qquad \|\b_\e\|_{\tilde{\mathfrak{M}}}\leqslant C(\e^n\rho_5^{-n} +\e^2|\ln\e|),\quad n=2,
\end{equation*}
\end{lemma}

\begin{proof}
Given an arbitrary function $u\in W_2^1(B_{\rho_5}(M_k^\e))$, by Lemma~\ref{lm3.6} with $\om_k^\e=\emptyset$ and $(\e,\eta)$ replaced by $(\rho_5,\e\rho_5^{-1})$  we obtain:
\begin{align*}
&\|u\|_{L_2(B_{\e R_3}(M_k^\e))}^2\leqslant C(\e^n\rho_5^{-n}+\e^2)\|u\|_{W_2^1(B_{\rho_5}(M_k^\e))}^2, &&\hspace{-1.7 true cm} n\geqslant 3,
\\
& \|u\|_{L_2(B_{\e R_3}(M_k^\e))}^2\leqslant C(\e^n\rho_5^{-n}+\e^2|\ln\e|)\|u\|_{W_2^1(B_{\rho_5}(M_k^\e))}^2, &&\hspace{-1.7 true cm} n=2,
\end{align*}
where $C$ is a constant independent of $\e$, $u$ and $k$.
Hence,   by the uniform boundedness of $\b_\e$ stated in Lemma~\ref{lm5.2},
\begin{equation*}
\big|(\b_\e u,v)_{L_2(B_{\e R_3}(M_k^\e))}\big|\leqslant C(\e^n\rho_5^{-n}+\e^2)\|u\|_{W_2^1(B_{\rho_5}(M_k^\e))}^2,
\end{equation*}
where $C$ is a constant independent of $\e$, $u$ and $k$. Summing up this estimate over $k\in\mathds{M}^\e$, we complete the proof.
\end{proof}

The proven lemma says that if the distances between the cavities are much larger than $\e$, then Assumption~\ref{A7} holds with $\b=0$. In particular, this is the case when we deal with  cavities separated by finite distances.

The second situation describes a perforation, which can be regarded as a general perturbation of a periodically distributed perforations. We choose a fixed lattice $\G$ in $\mathds{R}^n$ with a periodicity cell $\square$.
Then we define the set $\G_\e$ by formula (\ref{6.29}) with $\rho_1(\e)=\e$ and in each rescaled cell $\e k+\e\square$, $k\in \G_\e$, we choose a point $M_k^\e$ such that $B_{\e R_3}(M_k^\e)\subset \e k+\e\square$, $k\in \G_\e$. Then we arbitrary choose the corresponding cavities $\om_{k,\e}$ and in the case $n\geqslant 3$ we additionally assume that the constants $K_{k,\e}$ satisfy the identity
\begin{equation}\label{6.3}
(2-n)K_{k,\e}=\Psi_\e(M_k^\e).
\end{equation}
Here
 $\Psi_\e\in L_\infty(\Om)\cap C(\overline{\Om})$
is some family of functions such that
\begin{equation}\label{6.7}
\rho_6(\e):=\max\limits_{x\in\overline{\Om}} |\Psi_\e(x)-\Psi_0(x)|\to+0,\quad \e\to+0,
\end{equation}
where $\Psi_0\in L_\infty(\Om)\cap C(\overline{\Om})$ is some uniformly continuous in $\overline{\Om}$ function, namely,
\begin{equation}\label{6.8}
\rho_7(\e):=\max\limits_{k\in\mathds{M}^\e} \max\limits_{x,y\in\e\overline{\square}+\e k} |\psi_0(x)-\psi_0(y)| \to +0,\quad \e\to+0.
\end{equation}
We stress that condition (\ref{6.3}) is imposed only on the constants $K_{k,\e}$ and not on the shapes of the corresponding cavities. This means that the cavities corresponding to different $k$ are not necessarily of the same shapes even if the constants $K_{k,\e}$ coincide.
In the case $n=2$ we let $\Psi_\e(x):=\Psi_0(x)\equiv 1$. These conditions ensure (\ref{6.24}) with $\rho_1(\e)=\e$ and
\begin{equation*}
\b=\frac{2-n}{|\square|}\psi_0,\quad n\geqslant 3, \qquad \b=\frac{1}{|\square|},\quad n\geqslant 2.
\end{equation*}
Indeed,
\begin{align*}
 \int\limits_{\e\square+\e k} \big(\b_\e(x)-\b(x)\big)\,dx=&
\e^n \psi_\e(M_k^\e)-\frac{1}{ |\square| }\int\limits_{\e\square+\e k} \psi_0(x) \,dx
= \int\limits_{\e\square+\e k} \frac{\psi_\e(M_k^\e)-\psi_0(x)}{|\square|} \,dx
\\
=& \int\limits_{\e\square+\e k} \frac{\psi_\e(M_k^\e)-\psi_0(M_k^\e)}{|\square|}\,dx
+ \int\limits_{\e\square+\e k} \frac{\psi_0(M_k^\e)-\psi_0(x)}{|\square|}\,dx.
\end{align*}
The right hand side of the above identity
 can be estimated by means of the functions $\rho_6$, $\rho_7$ introduced in (\ref{6.7}), (\ref{6.8}) and this yields
\begin{equation*}
\frac{1}{\e^n}\bigg|\int\limits_{\e\square+\e k}(\b_\e-\b)\,dx\bigg|\leqslant \rho_6 +\rho_7.
\end{equation*}
This is exactly condition (\ref{6.24}) for our case and by (\ref{6.25})
we obtain
\begin{equation*}
\|\b_\e-\b\|_{\tilde{\mathfrak{M}}}\leqslant C  (\e +\rho_6+\rho_7),\quad n\geqslant 3,\qquad \|\b_\e-\b\|_{\tilde{\mathfrak{M}}}\leqslant C  \e,\quad n=2.
\end{equation*}

Our next step is to show how to generate new perforations obeying Assumption~\ref{A7} if we are given one already obeying this assumption. The first way is provided by the following lemma.

\begin{lemma}\label{lm6.3}
Let a perforation described by the points $M_k^\e$ and cavities $\om_{k,\e}$ obey Assumption~\ref{A7}. Let $\tilde{M}_k^\e$, $k\in\mathds{M}^\e$, be another set of points satisfying Assumption~\ref{A1} with the same constants $R_i$, $i=1,\ldots,4$, as for $M_k^\e$, such that
\begin{equation*}
|M_k^\e-\tilde{M}_k^\e|\leqslant C\e
\end{equation*}
with some constant $C$ independent of $k$ and $\e$.
Then the function $\tilde{\b}_\e$ corresponding to the perforation described by the points $\tilde{M}_k^\e$ and the same cavities $\om_{k,\e}$ also obeys Assumption~\ref{A7} with the same function $\b$ and  the estimate holds:
\begin{equation*}
\|\tilde{\b}_\e-\b\|_{\mathfrak{M}}\leqslant \|\b_\e-\b\|_{\mathfrak{M}} + C\e^\frac{1}{2},
\end{equation*}
where $C$ is a constant independent of $\e$.
\end{lemma}

\begin{proof}
It is clear that
\begin{equation*}
\|\tilde{\b}_\e-\b\|_{\mathfrak{M}}\leqslant \|\b_\e-\b\|_{\mathfrak{M}} + \|\b_\e-\tilde{\b}_\e\|_{\mathfrak{M}}
\end{equation*}
and this is why it is sufficient to estimate just the second term in the right hand of this inequality. We are going to do this by means of condition (\ref{6.24}). Namely, we let $\rho_1(\e):=\e^\frac{1}{2}$. Then the integral in (\ref{6.24})
can be rewritten as a sum of the integrals over the balls $B_{\e R_3}(M_k^\e)$ and $B_{\e R_3}(\tilde{M}_k^\e)$. If for some $k$
both these balls are contained in the cell $\e^\frac{1}{2} z + \e^\frac{1}{2}\square$, then their contributions to the total integral cancel out just due to the definitions of the points $\tilde{M}_k^\e$ and of the function $\tilde{\b}_\e$. Hence, only the balls  $B_{\e R_3}(M_k^\e)$ and $B_{\e R_3}(\tilde{M}_k^\e)$ intersecting with the boundary $ \e^\frac{1}{2} z + \e^\frac{1}{2}\p\square $ contribute to the considered integral. Then the total number of such balls is proportional to the measure of this boundary, which is of order $\sim \e^{\frac{n-1}{2}}$ and the total measure of such balls is obviously estimated by $C\e^{\frac{n+1}{2}}$ with some fixed constant $C$. Since the functions $\b_\e$ and $\tilde{\b}_\e$ are uniformly bounded, see Lemma~\ref{lm5.2}, we then get the estimate
\begin{equation*}
\frac{1}{\e^{\frac{n}{2}}} \Bigg|\int\limits_{\e^\frac{1}{2} z+\e^\frac{1}{2}\square}
\big(\b_\e(x)-\b(x)\big)\,dx\Bigg|\leqslant C \e^\frac{1}{2}
\end{equation*}
and we arrive at (\ref{6.24}) with $\rho_2(\e)=\e^\frac{1}{2}$.
Employing then estimate (\ref{6.25}) for $\|\b_\e-\tilde{\b}_\e\|_{\tilde{\mathfrak{M}}}$ and (\ref{6.2}), we complete the proof.
\end{proof}

The proven lemma shows that given a perforation obeying Assumption~\ref{A7}, we can shift the points $M_k^\e$ by the distance of order $O(\e)$   provided the new points satisfy Assumption~\ref{A1}. This gives an easy way to generate many new non-periodic perforations from a given one keeping Assumption~\ref{A7} satisfied.

The second way of generating new perforations obeying Assumption~\ref{A7} is as follows. Suppose that we are given two perforations described by  $M_k^\e$, $\om_{k,\e}$, $k\in\mathds{M}^\e$, and $\hat{M}_k^\e$, $\hat{\om}_{k,\e}$, $k\in\hat{\mathds{M}}^\e$. Let these perforations satisfy Assumption~\ref{A7} respectively with the functions $\b_\e$, $\b$ and
$\hat{\b}_\e$, $\hat{\b}$. Consider then the union of these perforations formed by the unions of the points and cavities $M_k^\e\cup \hat{M}_j^\e$, $\om_{k,\e}$, $\om_{j,\e}$, $k\in\mathds{M}^\e$, $j\in\mathds{M}^\e$, and let this union of the perforations satisfy Assumption~\ref{A1}. Then function (\ref{2.9}) corresponding to this union of the perforations is $\b_\e+\hat{\b}_\e$ and it satisfies Assumption~\ref{A7} with the limiting function $\b+\hat{\b}$ thanks to the following simple estimate:
\begin{equation*}
\|\b_\e+\hat{\b}_\e- \b -\hat{\b} \|_{\mathfrak{M}}\leqslant
\|\b_\e - \b \|_{\mathfrak{M}} +  \| \hat{\b}_\e-  \hat{\b} \|_{\mathfrak{M}}.
\end{equation*}

It is also possible to remove some cavities from a given perforation keeping at the same time Assumption~\ref{A7} and this is our third way of producing new perforations. Namely, given a perforation described by the points and cavities $M_k^\e$, $\om_{k,\e}$, $k\in\mathds{M}^\e$ and obeying Assumption~\ref{A7}, suppose that there is a subset $\hat{\mathds{M}}^\e\subset \mathds{M}^\e$ such that the corresponding perforation
satisfies Assumption~\ref{A7} with $\b=0$; the associated function (\ref{2.9}) is denoted by $\hat{\b}_\e$. Then we consider a difference of perforations corresponding to $\mathds{M}^\e\setminus \hat{\mathds{M}}^\e$ and we see that its function (\ref{2.9}) is $\b_\e-\hat{\b}_\e$. Hence,
\begin{equation*}
\|\b_\e-\hat{\b}_\e-\b\|_{\mathfrak{M}}\leqslant \|\b_\e-\b\|_{\mathfrak{M}} +\|\hat{\b}_\e\|_{\mathfrak{M}},\qquad \|\hat{\b}_\e\|_{\mathfrak{M}}\to+0,\quad \e\to+0,
\end{equation*}
and the introduced difference of perforations also satisfies Assumption~\ref{A7} with the same function $\b$.

The fourth way of producing new perforations is to vary the shapes of the cavities. In the dimension $n=2$ the function $\b_\e$ is independent on the shapes of the cavities and we therefore have a very rich freedom in choosing the shapes of the cavities. As $n\geqslant 3$, the shapes of the cavities are reflected in the constants $K_{k,\e}$. Then, given a perforation obeying Assumption~\ref{A7} with functions $\b_\e$ and $\b$, one can deform slightly the shapes of the cavities so that new constants $\hat{K}_{k,\e}$ differ from $K_{k,\e}$ by a small quantity, namely, $|\hat{K}_{k,\e}-K_{k,\e}|\leqslant \rho_8(\e)$, where $\rho_8(\e)\to+0$ as $\e\to+0$. Then it is clear that the new function $\hat{\b}_\e$ corresponding to the constants $\hat{K}_{k,\e}$ satisfies the estimate
\begin{equation*}
\|\hat{\b}_\e-\b\|_{\mathfrak{M}}\leqslant \|\b_\e-\b\|_{\mathfrak{M}} + \rho_8(\e),
\end{equation*}
which means that Assumption~\ref{A7} holds also for modified perforation.

\section*{Acknowledgment}

The author thanks A.I.~Nazarov for a valuable discussion on the regularity of solutions to elliptic boundary value problems. The authors is grateful to the referee for useful remarks.

\section*{Funding}

The work is supported by  the Czech Science Foundation within the project 22-18739S.

\section*{Conflict of interest}

The author declares that he has no conflicts of interest.

\section*{Data availability}

Not applicable in the manuscript as no datasets were generated or analysed during the current study.

\end{document}